\documentclass[reqno]{amsart}

\usepackage[dvipsnames]{xcolor}

\usepackage{baskervald}
\usepackage[mathscr]{euscript}
\usepackage{etoolbox}
\usepackage{hyperref}
\usepackage{cleveref}

\usepackage{amssymb}
\usepackage{amsmath}
\usepackage{tikz-cd}

\usepackage{epigraph}
\hypersetup{
    colorlinks=true,
    linkcolor=Maroon,
    citecolor=MidnightBlue,
    urlcolor=MidnightBlue
}

\usepackage{todonotes}

\DeclareMathSymbol{A}{\mathalpha}{operators}{`A}
\DeclareMathSymbol{B}{\mathalpha}{operators}{`B}
\DeclareMathSymbol{C}{\mathalpha}{operators}{`C}
\DeclareMathSymbol{D}{\mathalpha}{operators}{`D}
\DeclareMathSymbol{E}{\mathalpha}{operators}{`E}
\DeclareMathSymbol{F}{\mathalpha}{operators}{`F}
\DeclareMathSymbol{G}{\mathalpha}{operators}{`G}
\DeclareMathSymbol{H}{\mathalpha}{operators}{`H}
\DeclareMathSymbol{I}{\mathalpha}{operators}{`I}
\DeclareMathSymbol{J}{\mathalpha}{operators}{`J}
\DeclareMathSymbol{K}{\mathalpha}{operators}{`K}
\DeclareMathSymbol{L}{\mathalpha}{operators}{`L}
\DeclareMathSymbol{M}{\mathalpha}{operators}{`M}
\DeclareMathSymbol{N}{\mathalpha}{operators}{`N}
\DeclareMathSymbol{O}{\mathalpha}{operators}{`O}
\DeclareMathSymbol{P}{\mathalpha}{operators}{`P}
\DeclareMathSymbol{Q}{\mathalpha}{operators}{`Q}
\DeclareMathSymbol{R}{\mathalpha}{operators}{`R}
\DeclareMathSymbol{S}{\mathalpha}{operators}{`S}
\DeclareMathSymbol{T}{\mathalpha}{operators}{`T}
\DeclareMathSymbol{U}{\mathalpha}{operators}{`U}
\DeclareMathSymbol{V}{\mathalpha}{operators}{`V}
\DeclareMathSymbol{W}{\mathalpha}{operators}{`W}
\DeclareMathSymbol{X}{\mathalpha}{operators}{`X}
\DeclareMathSymbol{Y}{\mathalpha}{operators}{`Y}
\DeclareMathSymbol{Z}{\mathalpha}{operators}{`Z}

\newcommand{\Maps}{\mathrm{Maps}}

\newcommand{\Spec}{\operatorname{Spec}}

\newcommand{\End}{\operatorname{End}}
\newcommand{\id}{\mathbf{1}}
\newcommand{\tr}{\operatorname{tr}}

\newcommand{\ind}{\operatorname{ind}}
\newcommand{\res}{\operatorname{res}}

\newcommand{\IND}{\mathbf{ind}}
\newcommand{\RES}{\mathbf{res}}
\newcommand{\OBLV}{\mathbf{oblv}}
\newcommand{\COIND}{\mathbf{coind}}
\newcommand{\CORES}{\mathbf{cores}}
\newcommand{\GAMMA}{\mathbf{\Gamma}}
\newcommand{\LOC}{\mathbf{Loc}}

\newcommand{\Cat}{\mathrm{Cat}}
\newcommand{\Poset}{\mathrm{Poset}}
\newcommand{\Prop}{\mathrm{Prop}}
\newcommand{\Tw}{\operatorname{Tw}}

\newcommand{\DGCAT}{\textbf{dgCat}}
\newcommand{\twoDGCAT}{\textbf{2-dgCat}}
\newcommand{\MOD}{\textbf{-mod}}
\newcommand{\BIMOD}{\textbf{-bimod}}

\newcommand{\too}{\longrightarrow}
\newcommand{\mon}{\hookrightarrow}

\newcommand{\add}{\mathrm{add}}
\newcommand{\coev}{\mathrm{coev}}
\newcommand{\ev}{\mathrm{ev}}

\newcommand{\KD}{\mathrm{KD}}
\newcommand{\mult}{\mathrm{mult}}
\newcommand{\comult}{\mathrm{comult}}
\newcommand{\oblv}{\mathrm{oblv}}
\newcommand{\pr}{\mathrm{pr}}

\newcommand{\unit}{\mathrm{unit}}
\newcommand{\counit}{\mathrm{counit}}

\renewcommand{\AA}{\mathscr{A}}
\newcommand{\bfA}{\mathbf{A}}
\newcommand{\BB}{\mathscr{B}}
\newcommand{\CC}{\mathscr{C}}

\newcommand{\DD}{\mathscr{D}}

\newcommand{\bfE}{\mathbf{E}}
\newcommand{\FF}{\mathscr{F}}

\newcommand{\HH}{\mathbf{H}}

\newcommand{\Cot}{\mathbf{L}}
\newcommand{\MM}{\mathscr{M}}

\newcommand{\NN}{\mathscr{N}}
\newcommand{\bfO}{\mathbf{O}}
\newcommand{\OO}{\mathscr{O}}

\newcommand{\bfP}{\mathbf{P}}

\newcommand{\SSh}{{\mathbf{S}^\Rightarrow}}

\newcommand{\Tan}{\mathbf{T}}
\newcommand{\UU}{\mathscr{U}}

\newcommand{\WW}{\mathscr{W}}
\newcommand{\XX}{\mathscr{X}}
\newcommand{\YY}{\mathscr{Y}}
\newcommand{\ZZ}{\mathscr{Z}}

\newcommand{\Zc}{{\overset{\circ}{Z}}}

\newcommand{\fg}{\mathfrak{g}}

\newcommand{\Gm}{{\mathbf{G}_m}}

\newcommand{\act}{\mathrm{act}}

\newcommand{\cl}{\mathrm{cl}}
\newcommand{\co}{\mathrm{co}}

\newcommand{\der}{\mathrm{der}}

\newcommand{\dR}{\mathrm{dR}}

\newcommand{\Nilp}{\mathscr{N}}

\newcommand{\obj}{\mathrm{obj}}
\newcommand{\op}{\mathrm{op}}

\newcommand{\red}{\mathrm{red}}

\newcommand{\rev}{\mathrm{rev}}

\newcommand{\dgCat}{\operatorname{dgCat}}
\newcommand{\ShCat}{\operatorname{ShCat}}
\newcommand{\Fun}{\operatorname{Fun}}

\newcommand{\Sm}{\operatorname{Sm}}

\newcommand{\Bun}{\mathrm{Bun}}

\renewcommand{\mod}{\operatorname{-mod}}

\newcommand{\LS}{\mathrm{LS}}

\newcommand{\Sym}{\mathrm{Sym}}

\newcommand{\Coh}{\mathrm{Coh}}

\newcommand{\DMod}{\mathrm{D}}
\newcommand{\ICoh}{\mathrm{IndCoh}}
\newcommand{\ICohc}{\overset{\circ}{\ICoh}}

\newcommand{\QCoh}{\mathrm{QCoh}}
\newcommand{\Sh}{\mathrm{Sh}}
\newcommand{\Sing}{\mathrm{Sing}}
\newcommand{\Singc}{\overset{\circ}{\Sing}}

\newcommand{\Rep}{\mathrm{Rep}}

\newcommand{\Vect}{\mathrm{Vect}}

\newcommand{\pt}{\mathrm{pt}}

\newcommand{\simto}{\overset{\sim}{\rightarrow}}

\newcommand{\tightoverset}[2]{%
  \mathop{#2}\limits^{\vbox to -.5ex{\kern-0.75ex\hbox{$#1$}\vss}}}

\newcommand\newkevintheorem[3][]{%
    \newtheorem{#2}[kevinthm]{#3}%
    \csletcs{old#2}{#2}%
    \expandafter\renewcommand\csname #2\endcsname{\stepcounter{subsubsection}\csname old#2\endcsname}%
}

\newcommand\newappendixtheorem[3][]{%
    \newtheorem{#2}[appendixthm]{#3}%
    \csletcs{old#2}{#2}%
    \expandafter\renewcommand\csname #2\endcsname{\stepcounter{subsubsection}\csname old#2\endcsname}%
}

\crefname{appsec}{Appendix}{Appendices}

\newtheoremstyle{kevintheoremstyle}
{}% space above
{}% space below
{}% body font
{0pt}% space to indent
{\bfseries}% head font
{.}% separator between head and body
{4pt}% space after head
{\thmnumber{#2}\thmname{ #1}\thmnote{ (#3)}}%head

\theoremstyle{kevintheoremstyle}

\newkevintheorem{theorem}{Theorem}
\crefname{theorem}{Theorem}{Theorems}
\newkevintheorem{lemma}{Lemma}
\crefname{lemma}{Lemma}{Lemmas}
\newkevintheorem{proposition}{Proposition}
\crefname{proposition}{Proposition}{Propositions}
\newkevintheorem{corollary}{Corollary}
\crefname{corollary}{Corollary}{Corollaries}
\newkevintheorem{conjecture}{Conjecture}
\crefname{conjecture}{Conjecture}{Conjectures}
\newkevintheorem{construction}{Construction}
\crefname{construction}{Construction}{Constructions}
\newkevintheorem{computation}{Computation}
\crefname{computation}{Computation}{Computation}

\newappendixtheorem{aproposition}{Proposition}
\crefname{aproposition}{Proposition}{Propositions}

\makeatletter
\def\l@subsection{\@tocline{2}{0pt}{2.5pc}{5pc}{}}
\makeatother

\numberwithin{equation}{subsection}

\usepackage[T1]{fontenc}

\begin{document}

\title{Coherent sheaves, sheared D-modules, and Hochschild cochains}

\author{Dario Beraldo}
\author{Kevin Lin}
\author{Wyatt Reeves}

\begin{abstract}
We show that the category of ind-coherent sheaves on a quasi-smooth scheme is naturally tensored over the category of sheared D-modules on its shifted cotangent bundle, commuting with its natural action of categorified Hoschschild cochains. We prove that it defines a Morita equivalence as such. We then extend these results to quasi-smooth Artin stacks. 

As a consequence of our formalism, we are able to articulate a precise sense in which the space of unramified automorphic functions over a function field localizes over the stack of arithmetic Arthur parameters.
\end{abstract}

\maketitle

\tableofcontents

\epigraph{The anatomist presents to the eye the most hideous and disagreeable objects; but his science is useful to the painter in delineating
even a Venus or an Helen.}{David Hume, \emph{An Enquiry Concerning
Human Understanding}}

\section{Introduction}

\subsection{Motivation}

\subsubsection{} The question of the existence of finite projective resolutions for modules was first raised by Hilbert, who showed in his famous ``syzygy theorem'' that the polynomial ring $k[x_1, \ldots, x_n]$ has global dimension $n$. Another major result in the subject, due to Auslander--Buchsbaum--Serre, says that a Noetherian local ring has bounded global dimension if and only if it is regular. This indicates that the failure of a coherent sheaf to have a finite projective resolution is entirely determined by its behavior around the singular points of the underlying variety.

\subsubsection{} If $Y$ is a local complete intersection scheme, one can even produce a quantitative measure of how badly a coherent sheaf on $Y$ fails to have a finite projective resolution. Namely, there is a scheme $\Sing(Y)$ lying over $Y$, which is isomorphic to $Y$ over the smooth locus, but which has positive fiber dimension over the singular locus. (The fiber dimension at a singular point measures ``how pathological'' the singularity is.) $\Sing(Y)$ carries a contracting action of $\Gm$ and the fixed locus for this action is isomorphic to $Y$. 

For any coherent sheaf $\FF \in \Coh(Y)$ one can produce a closed conical subset $\mathrm{supp}(\FF) \subseteq \Sing(Y)$, which is contained in the zero section if and only if $\FF$ admits a finite projective resolution. The assignment $\FF \mapsto \mathrm{supp}(\FF)$ gives a notion of ``microlocal support'' for coherent sheaves. 

\subsubsection{} Recently, the microlocal geometry of coherent sheaves has played a key role in the resolution of the geometric Langlands conjecture: 
\[D(\Bun_G) \simeq \ICoh_{\NN}(\LS_{\check{G}}).\]
Evidently, the theory of microsupport of coherent sheaves is required to even state the result. One should think that the theorem gives a ``spectral decomposition'' of the category $D(\Bun_G)$. To first order, $D(\Bun_G)$ is tensored over $\QCoh(\LS_{\check{G}})$, which gives the decomposition of $D(\Bun_G)$ by Langlands parameters. However, the full statement of the equivalence shows that D-modules on $\Bun_G$ have a more refined decomposition in terms of coherent microsupport on $\LS_{\check G}$. 

\newcommand{\bfF}{\mathbf{F}}
\newcommand{\bfQ}{\mathbf{Q}}
\newcommand{\ol}[1]{\overline{#1}}
\newcommand{\totext}[1]{\xrightarrow{#1}}
\newcommand{\ul}[1]{\underline{#1}}
\newcommand{\ShvCat}{\mathrm{ShvCat}}
\newcommand{\usot}[1]{\underset{#1}{\otimes}}
\newcommand{\bfL}{\mathbf{L}}
\newcommand{\bfG}{\mathbf{G}}
\newcommand{\vphi}{\varphi}
\newcommand{\sN}{\mathscr{N}}
\newcommand{\rlto}{\rightleftharpoons}

\subsubsection{} The microsupport filtration on $\ICoh_{\NN}(\LS_{\check{G}})$ is expected to have the following ``experimental signature'': if $X$ is a curve over $\bfF_q$, then the stack $\LS^{\res}(X)$ of $\ell$-adic local systems with restricted variation was introduced in \cite{arinkin2020stack}. In \textit{loc. cit.} they conjecture a restricted version of the geometric Langlands equivalence 
\[\Sh_{\NN}(\Bun_G) \overset{?}{\simeq} \ICoh_{\NN}(\LS^{\res}_{\check{G}})\]
which is expected to be closely related to the de Rham form of the geometric Langlands correspondence. This equivalence is expected to be equivariant for the natural Frobenius endomorphism on each side, and one should be able to recover the unramified part of the Langlands correspondence for function fields 
\[\Fun_c(\Bun_G(\bfF_q); \ol{\bfQ}_\ell) \simeq \Gamma(\LS^{\mathrm{arith}}_{\check{G}}; \omega) \]
by taking the categorical trace of both sides. Now the filtration of $\ICoh_{\NN}(\LS^{\res}_{\check{G}})$ by microsupport decategorifies to a filtration on unramified automorphic functions, which is expected to encapsulate the theory of \textit{Arthur parameters}. 

\subsubsection{} Motivated by these considerations, in this paper we undertake a study of the microlocal geometry of coherent sheaves from a categorical perspective; and especially with the goal of computing categorical traces. 

\subsection{Statement of Results}

\subsubsection{} For the rest of the introduction, $Y$ will denote a quasi-smooth (derived) scheme unless otherwise stated.

\subsubsection{Microlocalization} Our first goal is to make precise the sense in which coherent sheaves lie over $\Sing(Y)$. Recall the so-called category of singularities
    \[
    \ICohc(Y) = \ICoh(Y)/\QCoh(Y),
    \]
One of the main accomplishments of \cite{arinkin2018category} was using the theory of singular support from \cite{arinkin2015singular} to equip $\ICoh(Y)$ with the following additional structure:

\begin{construction}[{\cite[Theorem 0.4.2]{arinkin2018category}}]
\label{construction:AG2}
The category $\ICohc(Y)$ is naturally a crystal of categories over $\bfP\Sing(Y)$.
\end{construction}

\subsubsection{} It is natural to wonder whether the theory of singular support actually encodes some additional structure on $\ICoh(Y)$ itself. In fact, this question was already raised in \cite[Remark 1.4.3]{arinkin2018category}. We offer the following answer:

\begin{construction}[{\Cref{construction:AG2_epsilon_scheme}}]
\label{construction:AG2_epsilon}
The category $\ICoh(Y)$ is naturally tensored over $\DMod(\Sing(Y))^\Rightarrow$.
\end{construction}

\subsubsection{} Here the notation $(-)^\Rightarrow$ refers to the \textit{shearing} functor, which is an automorphism of the 2-category of categories acted on (weakly) by $\Gm$ (see \cite[\S A.2]{arinkin2015singular} or \cite{beraldo2020spectral}). The key example to bear in mind is that if $A$ is a graded algebra then the category $A \mod$ is acted on $\Gm$, and 
\[(A \mod)^\Rightarrow \simeq (A^\Rightarrow) \mod,\]
where $A^\Rightarrow$ is the graded algebra whose $i$-th graded piece is $A_i[-2i]$. Because we don't seriously consider questions of $t$-exactness in this paper, for us shearing plays a technically important but conceptually minor role. Let us note for future reference another important property of shearing, which is that it induces an equivalence on categories of equivariant objects:
\begin{equation}\label{eq:shearing-equivt-objects}
    \CC^{\Gm} \simeq (\CC^\Rightarrow)^\Gm
\end{equation}

\subsubsection{} When $Y$ is a global complete intersection scheme, an action of $D(\Sing(Y))^\Rightarrow$ on $\ICoh(Y)$ was previously constructed in \cite[\S 2.2.5]{beraldo2020spectral}. One should think of this as giving a construction of the action in local coordinates. We should emphasize that in the setting of derived algebraic geometry it is highly non-trivial to glue constructions along coordinate patches, due to the presence of infinite towers of higher homotopy coherences. In fact, our \Cref{construction:AG2_epsilon} proceeds by a method quite different from the one in \textit{loc. cit.}

\newcommand{\bfT}{\mathbf{T}}

\subsubsection{Morita Equivalence} To state our next result, we should first recall the category $\HH(Y)$, introduced in \cite{beraldo2021center}. For any eventually coconnective scheme with perfect cotangent complex, the HKR theorem (see \cite[Proposition F.1.5]{arinkin2015singular} and \cite[Corollary G.2.7(b)]{arinkin2015singular} for a more precise reference) gives an isomorphism of $\bfE_1$-algebras between the algebra $\mathrm{HC}(Y)$ of Hochschild cochains, and $\UU(\bfT_Y[-1])$, the universal enveloping algebra of the shifted tangent sheaf\footnote{Here we are regarding $\bfT_Y$ as an object of $\QCoh$ and forming the universal enveloping algebra inside of $\QCoh$. Likewise, $\mathrm{HC}(Y)$ refers to the quasi-coherent sheaf of Hochschild cochains.} of $Y$. As a result, the category $\HH(Y)$ of left modules for $\mathrm{HC}(Y)$ can be described explicitly as 
\begin{equation}\label{eq:HH-description}
\ICoh_0(Y \times Y^\wedge_Y) = \ICoh(Y \times Y^\wedge_Y) \underset{\ICoh(Y)}{\times} \QCoh(Y).    
\end{equation}
The category $\HH(Y)$ may be thought of as a categorification of either Hochschild cochains or the ring of differential operators. From the explicit description \ref{eq:HH-description} we see that $\HH(Y)$ itself has a monoidal structure, coming from convolution. In fact, $\HH(Y)$ acts on $\ICoh(Y)$ via convolution, and this action is known to encode important facts about singular support for $\ICoh$. 

\subsubsection{}In \S \ref{s:bimodule-structure}, we show that the action of $D(\Sing(Y))^\Rightarrow$ on $\ICoh(Y)$ from \Cref{construction:AG2_epsilon} commutes with the well-known convolution action of $\HH(Y)$ on $\ICoh(Y)$. The main theorem of this paper is that $\ICoh(Y)$, as a $(D(\Sing(Y))^\Rightarrow, \HH(Y))$-bimodule, defines a Morita equivalence:
\begin{theorem}[{\Cref{thm:morita_D_H}}]\label{thm:main-in-intro}
\[ D(\Sing(Y))^\Rightarrow \MOD \simeq \HH(Y) \MOD.\]  
\end{theorem}

\subsubsection{}\label{s:toy-model} Here is a simple result one category level down which captures the spirit of our \Cref{thm:main-in-intro}: if $Y$ is a smooth affine variety over a field of characteristic zero, then $\OO_Y$ is a bimodule for the action of locally constant functions $\Gamma(\ul{k}_Y)$ and the ring of differential operators $\DD_Y$. Famously, locally constant functions are precisely the ones that commute with the action of every differential operator: 
\[\End_{\DD_X}(\OO_Y) \simeq \Gamma(\ul{k}_Y)\]
so the bimodule $\OO_Y$ defines a fully-faithful functor 
\begin{equation}\label{eq:toy-model-embedding}
\Gamma(\ul{k}_Y) \mod \hookrightarrow \DD_Y \mod.    
\end{equation}
One should think that there are analogies 
\begin{center}
\begin{tabular}{c|c}
     \Cref{thm:main-in-intro} & \S \ref{s:toy-model} \\
     \hline
     $D(\Sing(Y))^\Rightarrow$ & $\Gamma(\ul{k}_Y)$ \\
     $\HH(Y)$ & $\DD_Y$ \\
     $\ICoh(Y)$ & $\OO_Y$
\end{tabular}
\end{center}
(See \S \ref{s:intro-categorification} for more discussion of this point.)

\subsubsection{} It is notably almost never the case that $\DD_X$ can be identified with endomorphisms of $\OO_X$ as a module for locally-constant functions, and for this reason the reader might be surprised by the comparison between our \Cref{thm:main-in-intro}, where $\ICoh(Y)$ defines a Morita equivalence, and the toy model of \S \ref{s:toy-model}, where $\OO_Y$ merely defines an embedding in one direction. 

The key mechanism that improves the situation one category level up is Gaitsgory's \textit{1-affineness theorem} for de Rham prestacks (\cite[Theorem 2.6.3]{gaitsgory2015sheaves}): if we identify $\ul{k}_Y$ with $\OO_{Y_{\dR}}$, then failure of the map (\ref{eq:toy-model-embedding}) to be an equivalence amounts to a failure of 
\[\Gamma(\OO_{Y_{\dR}}) \mod \hookrightarrow \QCoh(Y_{\dR})\]
to be an equivalence, which says that $Y_{\dR}$ is not 0-affine. What the 1-affineness result of \textit{loc. cit.} shows is that this obstruction goes away in our categorified setting: 
\[D(Y)\MOD \simeq \ShvCat(Y_{\dR}).\]

\subsubsection{}\label{s:intro-categorification} Let's try to briefly justify the table of analogies in \S \ref{s:toy-model}. Recall that for a smooth scheme $Y$ there are canonical equivalences 
\begin{align*}
D(\Sing(Y))^\Rightarrow &\simeq D(Y) \\
\HH(Y) &\simeq \QCoh(Y \times Y^\wedge_Y)\\
\ICoh(Y) &\simeq \QCoh(Y)
\end{align*}
The analogy between $\QCoh(Y)$ and $\OO_Y$ is well-known (see for instance \cite{gaitsgory2015sheaves}). The assertion that $D(Y)$ categorifies $\Gamma(\ul{k}_X)$ says that there is an analogy between a function being locally constant and a sheaf admitting a flat connection. The assertion that $\QCoh(Y \times Y^\wedge_Y)$ categorifies $\DD_X$ has the following interpretation: we should perceive the operation of acting by a differential operator as being like convolving with an integral kernel ``supported on the formal neighborhood of the diagonal''.

\subsubsection{} For a smooth scheme, it is well-known that $\QCoh(Y)$ defines a Morita equivalence between $D(Y) \MOD$ and $\HH(Y) \MOD$. In fact, for any scheme it was shown in \cite[Lemma 4.2.7]{beraldo2019sheaves} that $\ICoh(Y)$ defines a Morita equivalence 
\[D(Y) \MOD \simeq \ICoh(Y \times Y^\wedge_Y) \MOD.\]
Our \Cref{thm:main-in-intro} amounts to saying that for a quasi-smooth scheme $Y$, if you ``increase the size'' of $D(Y)$ to $D(\Sing(Y))^\Rightarrow$, then on the dual side $\ICoh(Y \times Y^\wedge_Y)$ ``decreases in size'' to precisely become $\HH(Y)$. 

\subsubsection{Extension to stacks} Let us emphasize now that so far in this introduction, $Y$ has always been a \textit{scheme}. In fact, our \Cref{thm:main-in-intro}, as stated, fails for even the simplest stack $BG$ (see \S\ref{s:BG-calculation}). However, the equivalence of \Cref{thm:main-in-intro} is natural enough that we are able to produce an equivalence 
\begin{construction}[{\Cref{construction:ShCat_morita}}]\label{con:extension-to-stacks-intro}
\[\ShvCat_{\SSh}(\YY) \simeq \ShvCat_{\HH}(\YY)\]    
\end{construction}
for any quasi-smooth Artin stack $\YY$\footnote{For our precise assumptions on $\YY$, see \S\ref{s:stack-assumptions}}. 

\subsubsection{} To make sense of the statement of \Cref{con:extension-to-stacks-intro}, we need to recall the notion of a sheaf of categories with coefficients, introduced in \cite{beraldo2019sheaves} (although let's say right away that our notion of a sheaf of categories with coefficients will be slightly different from that of \textit{loc. cit.}). 

First, in \textit{loc. cit.}, a coefficient system is defined to be a functor 
\[\bfA: \mathrm{AffSch} \to \mathrm{Mo}(\dgCat)\]
where $\mathrm{Mo}(\dgCat)$ is the Morita bimodule category whose objects are monoidal dg categories, morphisms are bimodule categories, etc. The composite functor
\[\mathrm{AffSch} \totext{\bfA} \mathrm{Mo}(\dgCat) \totext{(-)\MOD} \twoDGCAT\]
can then right Kan extended from $\mathrm{AffSch}$ to $\mathrm{Pstk}$, producing a functor called $\ShvCat_{\bfA}$\footnote{In \textit{loc. cit.} this object is actually called $\ShCat^{\bfA}$. We apologize for the minor inconsistency in notation.}. More explicitly, if $\YY$ is any prestack,
\[\ShvCat_{\bfA}(\YY) = \lim_{S \to \YY}\bfA(S) \MOD,\]
where the limit is formed over affine schemes $S$ mapping to $\YY$. 

\subsubsection{} One of the achievements of \cite{beraldo2019sheaves} was to produce the 2-category $\ShvCat_{\HH}(\YY)$ through the construction described above. Our construction of $\ShvCat_{\SSh}(\YY)$ doesn't quite fit into the pattern outlined above; instead, for a quasi-smooth Artin stack $\YY$ we start with the functor 
\[\Sm_{/\YY} \totext{D(\Sing(-))^\Rightarrow \MOD} \twoDGCAT\]
and then set 
\[\ShvCat_{\SSh}(\YY) = \lim_{S \to \YY} D(\Sing(S))^\Rightarrow \MOD,\]
where the limit is formed over the smooth site of $\YY$. 

\subsubsection{} Although this construction of $\ShvCat_{\SSh}(\YY)$ is not obviously functorial in $\YY$, our \Cref{con:extension-to-stacks-intro} makes $\ShvCat_{\SSh}(\YY)$ into a functor out of quasi-smooth Artin stacks by transport of structure from $\ShvCat_{\HH}$. We characterize this functoriality by showing that our Morita equivalence matches the $(\HH(\ZZ), \HH(\YY))$-bimodule $\HH_{\ZZ \to \YY}$ with the sheared crystal $\ul{D}(\Sing(\YY)_\ZZ)^\Rightarrow$ over $\Sing(\ZZ) \times \Sing(\YY)$.

\subsubsection{Application to Traces} Let us now recall the formalism of categorical traces, which has been considered in \cite{benzvi2019nonlineartraces, kondyrev2020categorical, gaitsgory2022toy}. If $V$ is a finite dimensional vector space over $k$, then the trace map is 
\[\tr: \End(V) \simeq V \otimes V^\vee \totext{\ev} k.\]
If $\phi$ is a particular endomorphism of $V$, then we can further rewrite the formula for the trace of $\phi$ as 
\[\tr(\phi; V) = \left\{k \totext{\unit} V \otimes V^\vee \totext{\phi \otimes \id} V \otimes V^\vee \totext{\ev} k\right\}\]
where we have implicitly used the canonical identification $k \simeq \End(k)$. This formula for the trace admits a substantial generalization: if $\CC$ is any monoidal category, and $c \in \CC$ is a \textit{dualizable object} of $\CC$, then for any endomorphism $\phi \in \End(c)$ we can form the trace $\tr(\phi; c)$, which is an element of $\End(1_\CC)$:
\[\tr(\phi; c) = \left\{1_\CC \totext{\unit_c} c \otimes c^\vee \totext{\phi \otimes \id} c \otimes c^\vee \totext{\ev_c} 1_\CC \right\}\]

\subsubsection{} If $\CC$ is itself a 2-category, then $\End(1_\CC)$ is a 1-category, and the categorical trace enjoys the following functoriality: if $c$ and $c'$ are dualizable objects of $\CC$, equipped with endomorphisms $\phi$ and $\psi$, respectively, then for every \textit{right adjointable} morphism $f$, and every natural transformation $\alpha$ as in the diagram below
\[
\begin{tikzcd}
c \arrow[d, "\phi"] \arrow[r, "f"]& c' \arrow[d, "\psi"]\\
c \arrow[r, "f"] & c' 
\arrow[Rightarrow, "\alpha", from=2-1, to=1-2]
\end{tikzcd}
\]
we obtain a \textit{Chern class map} 
\[\cl(f, \alpha): \tr(\phi; c) \to \tr(\psi; c').\]

\subsubsection{}\label{s:intro-trace-tautology} The starting point for the theory of enhanced traces is the following tautological observation: the category $\Vect$ is dualizable as an object of $\dgCat$, and
\[\tr(\id_{\Vect}; \Vect) \simeq k.\]
Moreover, for any finite-dimensional vector space $V$ with an endomorphism $\phi$, we obtain a lax-commutative diagram in $\dgCat$
\[
\begin{tikzcd}
\Vect \arrow[d, "\id_{\Vect}"'] \arrow[r, "V \otimes (-)"] & \Vect \arrow[d, "\id_{\Vect}"] \\
\Vect \arrow[r, "V \otimes (-)"'] & \Vect \arrow[Rightarrow, "\phi", from=2-1, to=1-2]
\end{tikzcd}
\]
and 
\[\cl(V, \phi) = \tr(\phi; V)\]
in $\End(k)$. 

\newcommand{\Fr}{\mathrm{Fr}}

\subsubsection{}\label{s:weil-sheaf-enh-trace} Now let $\FF$ be a $\ol{\bfQ}_\ell$ Weil sheaf on a scheme $\YY/\ol{\bfF}_q$. By \S \ref{s:intro-trace-tautology}, the trace of Frobenius on the compactly supported cohomology of $\FF$ may be written as $\cl(H^*_c(\YY; \FF), \Fr_q)$. We can factor $H^*_c(\YY; \FF): \Vect \to \Vect$ as the composite 
\[\Vect \totext{\FF} \Sh(\YY) \totext{\Gamma_c} \Vect\]
and using the Weil structure on $\FF$ we can form the lax-commutative diagram 
\[
\begin{tikzcd}
\Vect \arrow[r, "\FF"] \arrow[d, "\id_\Vect"'] & \Sh(\YY) \arrow[d, "\Fr_q^*"] \arrow[r, "\Gamma_c"]& \Vect \arrow[d, "\id_\Vect"]\\
\Vect \arrow[r, "\FF"] & \Sh(\YY) \arrow[r, "\Gamma_c"]& \Vect
\arrow[Rightarrow, "\alpha", from=2-1, to=1-2] 
\arrow[Rightarrow, from=2-2, to=1-3]
\end{tikzcd}
\]
The class $\cl(\FF, \alpha)$ in $\tr(\Fr_q; \Sh(\YY)) = \Fun(\YY(\bfF_q); \ol{\bfQ}_\ell)$ is the usual function attached to $\FF$ under the sheaf-function correspondence, and the class map associated to the second square in the diagram is given by summing over the $\bfF_q$ points of $\YY$. 

\subsubsection{}\label{s:weil-cat-enh-trace} One category level up, if a dg-category $\CC$ arises as the global sections of a ``Weil'' sheaf of categories $\ul{\CC}$ over a space $\YY$ equipped with an endomorphism $\phi$, then we can form the object
\[\tr^{\mathrm{enh}}(\phi; \ul{\CC}) \in \tr(\phi; \ShvCat(\YY)) = \QCoh(\YY^\phi)\]
which is by definition $\cl(\CC, \alpha)$, where $\alpha$ is the Weil structure for $\CC$. This enhanced trace has the property that 
\[\Gamma(\YY^\phi; \tr^{\mathrm{enh}}(\phi; \ul{\CC})) \simeq \tr(\phi; \CC).\]

\subsubsection{} Morally, the previous paragraphs \ref{s:weil-sheaf-enh-trace} and \ref{s:weil-cat-enh-trace} show that the enhanced trace of a Weil object is like a version of the trace that ``lives over $\YY$''. Following this pattern, one should expect that the enhanced trace of a category $\CC$ acted on by $\HH(\YY)$ should be an upgrade of the usual trace that also remembers how $\CC$ lives over the stack of singularities. The final results of our paper confirm this expectation: 

\begin{computation}[{\Cref{computation:bimod_trace}}]\label{comp:trace-of-H}
Let $\phi$ be a schematic endomorphism of a quasi-smooth stack. Then 
\[\tr(\HH_\phi(\YY), \HH(\YY) \MOD) \simeq D(\Sing_2(\YY^\phi))^\Rightarrow.\]
\end{computation}

\begin{computation}[{\Cref{computation:proj_Nilp_trace}}]\label{comp:enh-trace-icoh}
If $\NN$ is a subset of $\Sing(\YY)$ stable under $\phi$ and $\alpha$ is the natural equivariance data for $\ICoh_{\NN}(\YY)$ as an object of $\HH(\YY) \MOD$, then 
\[\cl(\ICoh_\NN(\YY), \alpha) \simeq (\omega_{\NN^\phi})^\Rightarrow\]
as an object of $D(\Sing_2(\YY^\phi))^\Rightarrow$.
\end{computation} 

In the statement of \Cref{comp:enh-trace-icoh}, the notation $(\omega_{\NN^\phi})^\Rightarrow$ denotes the object which corresponds to $\omega_{\NN^\phi}$ in $D(\Sing_2(\YY^\phi))$ under the equivalence (\ref{eq:shearing-equivt-objects}). Here $\omega_{\NN^\phi}$ is the (naturally $\Gm$-equivariant) object of $D(\Sing_2(\YY^\phi))$ which is the pushforward of the dualizing sheaf on the fixed-point stack
\[\NN^\phi = \NN \times_{\Sing(\YY \times \YY)} \Sing_\phi (\YY).\]

\subsubsection{} Now let's fix a smooth proper curve $X$ over $\ol{\bfF}_q$ and set $\YY = \LS^{\res}_{\check{G}}(X)$. Let $\phi$ be the endomorphism of $\LS^{\res}_{\check{G}}$ coming by transport of structure from the geometric Frobenius endomorphism of $X$. Then $\Sing_2(\YY^\phi)$ is the stack of arithmetic Arthur parameters (see \cite[\S 24.5.5]{arinkin2020stack}). Strictly speaking, our \Cref{comp:trace-of-H} and \Cref{comp:enh-trace-icoh} don't apply as stated in this situation because $\YY$ is a formal stack rather than an honest Artin stack. However, the singular support theory of formal stacks is expected to be very similar to the theory of Artin stacks (see for instance \cite[\S 24.7.5]{arinkin2020stack}). We expect that our \Cref{comp:trace-of-H} and \Cref{comp:enh-trace-icoh} extend ``as-is'' to formal stacks; this gives a precise expectation for how the space of automorphic functions localizes over the stack of Arthur parameters. 

\subsection{Examples}

\subsubsection{} In this section we will illustrate what our main results amount to in some simple examples. 

\subsubsection{The case of a shifted vector space} Let's first consider our Morita equivalence \Cref{thm:main-in-intro} for the simplest quasi-smooth scheme, $\Omega V$. Let $V$ be a vector space, viewed as an affine scheme. Then 
    \[
    \Omega V = \pt \underset{V}{\times} \pt
    \]
is a quasi-smooth affine scheme with
    \[
    \Sing(\Omega V)_\dR = (V^\vee)_\dR.
    \]
The main results of this paper are particularly transparent when one takes advantage of the Koszul duality equivalence
\begin{equation}
\label{eqn:KD}
    \KD : 
    \ICoh(\Omega V)
    \simeq
    \Sym(V[-2])\mod
    \simeq
    \QCoh(V^\vee)^\Rightarrow.
\end{equation}

\subsubsection{}
\label{subsubsection:KD_review}
Let us recall some basic features of the equivalence (\ref{eqn:KD}) for future reference:
\begin{itemize}
    \item[(i)] Koszul duality intertwines the action of $\Gm$ on $\Omega V$ by dilation with the action of $\Gm$ on $V^\vee$ by \emph{inverse} dilation.

    \item[(ii)] Koszul duality intertwines the Serre self-duality of $\ICoh(\Omega V)$ with the self-duality of $\QCoh(V^\vee)^\Rightarrow$, \emph{up to the negation involution of $V$}
        \[
        \tau : 
        V \to V
        \qquad
        v \mapsto -v.
        \]
\end{itemize}
Observe that $V$, and hence $\Omega V$ is a commutative group scheme. Therefore, $\ICoh(\Omega V)$ acquires the structure of a symmetric monoidal category under convolution.
\begin{itemize}
    \item[(iii)] Koszul duality intertwines the convolution symmetric monoidal structure on $\ICoh(\Omega V)$ with the pointwise tensor symmetric monoidal structure on $\QCoh(V^\vee)$.

    \item[(iv)] Koszul duality intertwines the pointwise tensor symmetric monoidal structure on $\ICoh(\Omega V)$ with the convolution symmetric monoidal structure on $\QCoh(V^\vee)$.
\end{itemize}

\subsubsection{} In terms of Koszul duality, the action of $\DMod(\Sing(\Omega V))^\Rightarrow$ on $\ICoh(\Omega V)$ is the action of $\DMod(V^\vee)^\Rightarrow$ on $\QCoh(V^\vee)^\Rightarrow$ induced by the symmetric monoidal functor
    \[
    \oblv^\Rightarrow :
    \DMod(V^\vee)^\Rightarrow
    \to
    \QCoh(V^\vee)^\Rightarrow.
    \]

\subsubsection{} By \S\ref{subsubsection:KD_review}(iii), (iv), the equivalence
\begin{equation}
\label{eqn:KD_KD_tau}
    \ICoh(\Omega V \times \Omega V) 
    \overset{\KD \otimes \KD}{\too}
    \QCoh(V^\vee \times V^\vee)^\Rightarrow
    \overset{\id \otimes \tau}{\too}
    \QCoh(V^\vee \times V^\vee)^\Rightarrow
\end{equation}
intertwines the convolution monoidal structures on its source and target, compatibly with their actions on the two sides of (\ref{eqn:KD}). 

\medskip

It follows that (\ref{eqn:KD_KD_tau}) restricts to an equivalence
    \[
    \HH(\Omega V) 
    \simeq
    \QCoh\big( 
    (V^\vee \times V^\vee)^\wedge_{V^\vee} 
    \big)^\Rightarrow
    \]
intertwining the action of $\HH(\YY)$ on $\ICoh(\YY)$ with the action of
    \[
    \QCoh\big( (V^\vee \times V^\vee)^\wedge_{V^\vee} \big)^\Rightarrow
    \]
on $\QCoh(V^\vee)^\Rightarrow$ obtained by recognizing
    \[
    (V^\vee \times V^\vee)^\wedge_{V^\vee} 
    \simeq 
    V^\vee \! \underset{(V^\vee)_\dR}{\times} \! V^\vee
    \]
as the infinitesimal groupoid of $V^\vee$.

\medskip

So the Morita equivalence of \Cref{thm:main-in-intro} in this case is a manifestation of the well-known fact that, for any smooth scheme $S$, the category $\QCoh(S)$ is a Morita bimodule for $\DMod(S)$ and $\QCoh((S \times S)^\wedge_S)$.

\subsubsection{The case of a classifying stack} Let us now examine the equivalence
    \[
    \ShCat_{\HH} \simeq \ShCat_{\mathbf{S}^\Rightarrow}
    \]
for the smooth stack $\YY = BG$. In this subsection we will only consider \emph{strong} action of groups on categories. 

\subsubsection{}\label{s:BG-calculation} In this case both sides of the equivalence are familiar to geometric representation theory:
\begin{itemize}
    \item On the one hand,
        \[
        \ShCat_\SSh(BG) \simeq G\MOD
        \]
    identifies with the $(\infty, 2)$-category of dg categories with a strong action of $G$.

    \item On the other hand,
        \[
        \HH(BG) 
        \simeq \ICoh((BG \times BG)^\wedge_{BG})
        \simeq (\fg \mod \otimes \fg \mod)^G
        \]
    is the category of Harish--Chandra bimodules for $G$.
\end{itemize}

Therefore, the equivalence (\ref{eqn:morita_ShCat}) in this case recovers the well-known equivalence (see \cite[Equation (2.6)]{beraldo2019sheaves}) between categories acted on strongly by $G$ and modules for the category of Harish--Chandra bimodules for $G$.

\subsubsection{} In the case at hand, the adjunction
    \[
    \LOC_\SSh :
    \DMod(\Sing(\YY))^\Rightarrow \MOD
    \rightleftharpoons
    \ShCat_\SSh(\YY)
    : \GAMMA_\SSh
    \]
identifies with the adjunction
    \[
    \DMod(BG) \MOD
    \rightleftharpoons
    G\MOD
    \]
defined by the trivial bimodule
    \[
    \Vect.
    \]

\medskip

This adjunction is a colocalization because
    \[
    \End_G(\Vect) \simeq \Vect^G \simeq \DMod(BG).
    \]
On the other hand,
    \[
    \End_{\DMod(BG)}(\Vect)
    \simeq
    \Vect \underset{\DMod(BG)}{\otimes} \Vect
    \simeq
    \DMod(G)_0
    \]
is the full subcategory of $\DMod(G)$ generated by the constant sheaf.

\subsubsection{} To summarize, the failure of the $(\LOC_\SSh, \GAMMA_\SSh)$ adjunction to be an equivalence is tied to the failure of D-affineness for stacks, which in turn manifests as the failure of the relative categorical Kunneth formula for the fiber product:
    \[
    \begin{tikzcd}
        G \arrow[d] \arrow[r] & \pt \arrow[d] \\
        \pt \arrow[r] & BG
    \end{tikzcd}
    \]

\subsection{Miscellania}

\subsubsection{} In this section we will review some of the important technical ideas in this paper that may be of interest to experts. 

\subsubsection{} At various points in this paper homotopy coherence issues arise in the context of crystals of categories. Our \Cref{prop:rel_crys_locally_ff} provides an approach for dealing with these issues. A simplified version of the statement is that if $\YY$ is a scheme and $\CC_1$ and $\CC_2$ are two module categories for $D(\YY)$, then it is a \textit{property} whether a functor $F: \CC_1 \to \CC_2$ is $D(\YY)$ linear. We use this proposition in a crucial way both when we construct an action of $D(\Sing (\YY))^\Rightarrow$ on $\ICoh(\YY)$ and when we extend our theory from schemes to stacks. 

We should note that this idea is implicit in the main construction of \cite{arinkin2018category}, which gives the category of singularities $\ICohc(Z)$ the structure of a crystal of categories over $\bfP(\Sing(Z))$. In this construction, the crystal structure is produced by specifying, for each $S \to \bfP(\Sing(Z))$ a \textit{full subcategory} of $\QCoh(S) \otimes \ICohc(Z)$. An even earlier antecedent for this sort of idea is \cite[\S 4.2.1]{gaitsgory2015sheaves}, which observes that for a closed embedding of schemes, the functor of coinduction of crystals of categories is actually fully-faithful. 

\subsubsection{} The main idea of our construction of the action of $D(\Sing (\YY))^\Rightarrow$ on $\ICoh(\YY)$ is that there is a close relationship between $\ICoh(\YY)$ and the category of singularities $\ICohc(\YY \times \Omega \bfA^1)$, which allows us to bootstrap from the results of \cite{arinkin2018category}. The idea of relating $\ICoh$ on one variety to the category of singularities on another has appeared already in the work of M. Umut Isik \cite{isik2013equivalence}. 

\subsubsection{} In the course of constructing the Morita equivalence between $D(\Sing(\YY))^\Rightarrow \MOD$ and $\HH(\YY) \MOD$ we give a new proof of an important result of Ben-Zvi, Gunningham, and Orem \cite[Theorem 1.3]{ben2020highest} concerning when a bimodule defines a Morita equivalence (see \S \ref{s:cat-framework}).

\subsection{Conventions}

\subsubsection{} The notation
    \[
    \mathsf{L} : \CC \rightleftharpoons \DD : \mathsf{R}
    \]
will always mean that $(\mathsf{L}, \mathsf{R})$ is a pair of adjoint functors. We say that an adjunction is a localization (resp. a colocalization) if the right (resp. left) adjoint is fully faithful.

\subsubsection{} In this paper we will only deal with prestacks that are locally almost of finite type \cite[\S 2.1.7]{gaitsgorystudyI}. We remind the reader that the class of such prestacks is stable under the formation of de Rham prestacks \cite[Proposition 4.1.1.4]{gaitsgorystudyII}. 

\subsubsection{}\label{s:stack-assumptions} We make the running assumption (see \cite[\S 1.3]{beraldo2019sheaves}) that our stacks are quasi-compact, perfect, bounded, and locally of finite presentation. A perfect stack has affine diagonal (see \cite[\S 1.1]{ben2010integral}) and is therefore QCA.

\subsubsection{} For a quasi-smooth algebraic stack $\YY$, we will denote
    \[
    \Sing(\YY) = \Spec_\YY \big( \Sym(\Tan_\YY[1]) \big),
    \]
the `total space' of $\Cot_\YY[-1]$. This notation differs slightly from the more standard notation of \cite{arinkin2015singular} and \cite{arinkin2018category}, which used this expression to denote $\Sing(\YY)^\cl$.

\subsection{Acknowledgements}

We are grateful to Dennis Gaitsgory for his support. WR was supported by the NSF Graduate Research Fellowship under Grant No. DGE1745303. 

\section{Construction of the action}

\label{section:constr}

The goal of this section is to perform the the following:

\begin{construction}
\label{construction:AG2_epsilon_scheme}
Let $Z$ be a quasi-smooth scheme. Then $\ICoh(Z)$ receives a canonical action of $\DMod(\Sing(Z))^\Rightarrow$.
\end{construction}

The outline of this section is as follows: 
\begin{itemize}
    \item In \S\ref{subsection:review_shcat} we briefly review material about sheaves of categories from \cite{gaitsgory2015sheaves}. 
    \item In \S\ref{subsection:sheared_act} we describe a general construction for equipping a category with the structure of left module over the shearing of a symmetric monoidal category acted on by $\Gm$. This construction is a sort of \textit{deperiodization} result.
    \item In \S\ref{subsection:2-geometric-incarnation} we apply the framework of \S\ref{subsection:sheared_act} to carry out \Cref{construction:AG2_epsilon_scheme}, modulo the technical result \Cref{prop:two_actions_commute}. The main geometric input is that the vector space $V$ is open in $\bfP(V \oplus \id)$, which allows us to reduce to the results of \cite{arinkin2018category}.
    \item In \S\ref{subsection:rel_crys_of_cat} we prove a ``rigidity'' theorem for relative crystals of categories, which gives us control over homotopy coherence problems.
    \item In \S\ref{subsection:two_actions_commute_proof} we apply the framework of \S\ref{subsection:rel_crys_of_cat} to prove \Cref{prop:two_actions_commute}.
\end{itemize}

\subsection{Review of sheaves of categories}

\label{subsection:review_shcat}

\subsubsection{} For a morphism $f : T \to S$ of affine schemes, we will write
    \[
    \IND : 
    \QCoh(S)\MOD \rightleftharpoons \QCoh(T) \MOD
    : \RES
    \]
for the adjunction defined by extension and restriction of scalars along the map of symmetric monoidal categories
    \[
    f^* : \QCoh(S) \to \QCoh(T).
    \]

\subsubsection{} Let us recall some features of the theory of sheaves of categories.

\medskip

To every prestack $\YY$, one associates the symmetric monoidal $(\infty, 2)$-category
    \[
    \ShCat(\YY) = \lim_{S \to \YY} \QCoh(S)\MOD
    \]
of sheaves of categories on $\YY$. This limit is taken over test affine schemes mapping to $\YY$, with the transition functors given by $\IND$. The value of a sheaf of categories $\CC$ at a test affine scheme $S$ is denoted
    \[
    \GAMMA(S; \CC) \in \QCoh(S)\MOD.
    \]
The monoidal unit of $\ShCat(\YY)$ is denoted $\QCoh_{/\YY}$. It is given by the assignment
    \[
    \GAMMA(S; \QCoh_{/\YY}) = \QCoh(S).
    \]

\medskip

To a morphism of prestacks $f : \ZZ \to \YY$, one associates the adjunction
\begin{equation}
\label{eqn:CORES_COIND}
    \CORES_f : 
    \ShCat(\YY) \rightleftharpoons \ShCat(\ZZ)
    : \COIND_f.
\end{equation}
For any cartesian diagram of prestacks
\[
\begin{tikzcd}
    \ZZ' \arrow[r, "g' "] \arrow[d, "f' "] & \ZZ \arrow[d, "f"] \\
    \YY' \arrow[r, "g"] & \YY
\end{tikzcd}
\]
the natural transformation
    \[
    \CORES_{g} \circ \COIND_{f}
    \to
    \COIND_{f'} \circ \CORES_{g'}
    \]
obtained by adjunction from the identification
    \[
    \CORES_{f'} \circ \CORES_g \simto \CORES_{g'} \circ \CORES_f 
    \]
is an isomorphism.

\subsubsection{} Note that
    \[
    \COIND_f(\QCoh_{/\ZZ})
    \]
is naturally a sheaf of symmetric monoidal categories over $\ShCat(\YY)$, and the adjunction (\ref{eqn:CORES_COIND}) enhances to an adjunction
    \[
    \LOC_{\ZZ/\YY} : 
    \COIND_f(\QCoh_{/\ZZ})\MOD (\ShCat(\YY)) \rightleftharpoons \ShCat(\ZZ) 
    : \GAMMA_{\ZZ/\YY}
    \]
The morphism $f$ is called relatively 1-affine when this adjunction is an equivalence.

\subsubsection{} A prestack $\YY$ is said to be 1-affine if the morphism $\YY \to \pt$ is relatively 1-affine. We remind the reader that:
\begin{itemize}
    \item Every (eventually coconnective, laft, with affine diagonal) algebraic stack is 1-affine. 
    \item The de Rham prestack $Z_\dR$ of any finite type scheme $Z$ is 1-affine.
\end{itemize}

\subsection{Actions of sheared categories}
\label{subsection:sheared_act}

\subsubsection{} Let $\Gm\MOD$ denote the symmetric monoidal $(\infty, 2)$-category of dg categories with a weak action of $\Gm$. Recall from \cite[\S A.2]{arinkin2015singular} that this category is equipped with a symmetric monoidal auto-equivalence, known as shearing:
    \[
    (-)^\Rightarrow : \Gm\MOD \to \Gm\MOD.
    \]
We write $(-)^\Leftarrow$ for its inverse. 

\subsubsection{}
\label{subsubsection:sheared_action_setup}
Let $\AA$ be a symmetric monoidal category acted on by $\Gm$ (that is, $\AA$ is a commutative algebra object of $\Gm\MOD$). Since shearing is symmetric monoidal, $\AA^\Rightarrow$ is another symmetric monoidal category acted on by $\Gm$.

\medskip

The goal of this subsection is to give an answer to the question of what it takes equip a dg category $\MM$ with an action of $\AA^\Rightarrow$ (without regard for its $\Gm$-action).

\medskip

That is, we will give a description of 
    \[
    \AA^\Rightarrow\MOD
    =
    \AA^\Rightarrow\MOD(\dgCat).
    \]

\subsubsection{A toy model} The construction we are about to give has an analogue one category level down. Namely, for a graded algebra $A$, and a vector space $M$, the following data are equivalent:
\begin{itemize}
    \item[(a)] a non-graded $A$-module structure on $M$;
    \item[(b)] a graded $A[t^{\pm 1}]$-module structure on $M[t^{\pm 1}]$.
\end{itemize}
Indeed, if $M$ is given a non-graded $A$-module structure, then we endow $M[t^{\pm 1}]$ with the structure of a graded $A$-module, where a homogeneous element $a$ of degree $n$ acts by
    \[
    a \cdot (m t^k) = (am) t^{n + k}.
    \]
This action commutes with the natural graded action of $k[t^{\pm 1}]$ on $M[t^{\pm 1}]$ and hence makes $M[t^{\pm 1}]$ into a graded $A[t^{\pm 1}]$-module. To recover the original non-graded action of $A$ of $M$, observe that $A$ is the degree zero part of $A[t^{\pm 1}]$ and $M$ is the degree zero part of $M[t^{\pm 1}]$.

\subsubsection{} Since the morphism $f: \pt \to B\Gm$ is relatively 1-affine, we obtain an identification
    \[
    \dgCat \simeq \QCoh(\Gm)\MOD(\Gm\MOD)
    \qquad
    \MM \mapsto \MM \otimes \QCoh(\Gm)
    \]
where $\QCoh(\Gm)$ is viewed as a $\Gm$-equivariant symmetric monoidal category via its pointwise tensor product.

\medskip

Now let $\BB$ be any symmetric monoidal category. From the monadic adjunction
    \[
    \IND : 
    \dgCat
    \rightleftharpoons
    \BB \MOD(\dgCat)
    : \OBLV,
    \]
we obtain a commutative diagram
\begin{equation}\label{eq:shearing-action}
\begin{tikzcd}
    \BB \MOD \arrow[d, "\oblv"] \arrow[r, "\sim"] & \big(\BB \otimes \QCoh(\Gm) \big)\MOD(\Gm\MOD) \arrow[d, "\oblv"] \\
    \dgCat \arrow[r, "\sim"] & \QCoh(\Gm) \MOD (\Gm \MOD)
\end{tikzcd}    
\end{equation}
in which both horizontal arrows are equivalences.

\subsubsection{} Now let $\BB$ be the underlying symmetric monoidal category of $\AA^\Rightarrow$. By the projection formula for sheaves of categories we obtain the equivalences of commutative algebras
\begin{align*}
    \BB \otimes \QCoh(\Gm) 
    & \simeq \COIND_f \circ \CORES_f\left(\AA^\Rightarrow\right) \\
    & \qquad \simeq \AA^\Rightarrow \otimes \COIND_f(\Vect) 
    \simeq \AA^\Rightarrow \otimes \QCoh(\Gm)
\end{align*}
in $\Gm \MOD$.

As a result, for our choice of $\BB$, the diagram (\ref{eq:shearing-action}) rewrites as
\begin{equation}
\label{eqn:shearing_action_rewrite}
\begin{tikzcd}
    \AA^\Rightarrow \MOD \arrow[d, "\oblv"] \arrow[r, "\sim"] & \big( \AA^\Rightarrow \otimes \QCoh(\Gm) \big)\MOD(\Gm\MOD) \arrow[d, "\oblv"] \\
    \dgCat \arrow[r, "\sim"] & \QCoh(\Gm) \MOD (\Gm \MOD)
\end{tikzcd}
\end{equation}
The difference between this diagram and (\ref{eq:shearing-action}) is that here the $\Gm$ action on $\AA^\Rightarrow \otimes \QCoh(\Gm)$ is the natural one coming from the preexisting action on $\AA^\Rightarrow$, while in (\ref{eq:shearing-action}) the $\Gm$ action on $\BB$ is trivial.

\subsubsection{} Note that shearing induces a commutative diagram
\[
\begin{tikzcd}
    \big( \AA^\Rightarrow \otimes \QCoh(\Gm) \big)\MOD(\Gm\MOD) \arrow[d, "\oblv"] \arrow[r, "(-)^\Leftarrow"] & \big( \AA \otimes \QCoh(\Gm)^\Leftarrow \big)\MOD(\Gm\MOD) \arrow[d, "\oblv"] \\
    \QCoh(\Gm) \MOD (\Gm \MOD) \arrow[r, "(-)^\Leftarrow"] & \QCoh(\Gm)^\Leftarrow \MOD (\Gm \MOD)
\end{tikzcd}
\]
in which the horizontal arrows are equivalences. Concatenating this with (\ref{eqn:shearing_action_rewrite}), we obtain the commutative diagram
\[
\begin{tikzcd}
    \AA^\Rightarrow \MOD \arrow[d, "\oblv"] \arrow[r, "\sim"] & \big( \AA \otimes \QCoh(\Gm)^\Leftarrow \big)\MOD(\Gm\MOD) \arrow[d, "\oblv"] \\
    \dgCat \arrow[r, "\sim"] & \QCoh(\Gm)^\Leftarrow \MOD (\Gm \MOD)
\end{tikzcd}
\]
where the lower horizontal equivalence is given by
    \[
    \MM \mapsto \MM \otimes \QCoh(\Gm)^\Leftarrow.
    \]

\subsubsection{} This is our answer to the question posed in \S\ref{subsubsection:sheared_action_setup}: giving an action of $\AA^\Rightarrow$ on $\MM$ is equivalent to giving a $\Gm$-equivariant action of $\AA$ on
    \[
    \MM \otimes \QCoh(\Gm)^\Leftarrow,
    \]
commuting with the natural $\Gm$-equivariant action of $\QCoh(\Gm)^\Leftarrow$ on the right factor.

\subsubsection{Remark} Since 
\[\QCoh(\Gm)^\Leftarrow \simeq k[\beta^\pm] \mod\]
where $|\beta| = 2$, the category $\MM \otimes \QCoh(\Gm)^\Leftarrow$ is a 2-periodization of $\MM$. One may interpret the result of this subsection as a \textit{deperiodization} statement, showing what additional data is needed to reconstruct an action of $\AA^\Rightarrow$ on $\MM$ from an action of $\AA$ on the 2-periodization of $\MM$.

\subsection{Geometric incarnation of shearing}
\label{subsection:2-geometric-incarnation}
\subsubsection{} For the rest this section, we will let $Z$ denote a fixed choice of quasi-smooth scheme. In light of \S\ref{subsection:sheared_act}, Construction \ref{construction:AG2_epsilon_scheme} admits the following reformulation:

\begin{construction}
\label{construction:untwisted_action}
There is a natural lift of
    \[
    \ICoh(Z) \otimes \QCoh(\Gm)^\Leftarrow
    \in
    \QCoh(\Gm)^\Leftarrow\MOD(\Gm\MOD)
    \]
to an object of
    \[
    \big( \DMod(\Sing(Z)) \otimes \QCoh(\Gm)^\Leftarrow \big) \MOD
    (\Gm \MOD).
    \]
\end{construction}

\subsubsection{} Our next task is to reformulate Construction \ref{construction:untwisted_action} by giving a geometric incarnation of $\QCoh(\Gm)^\Leftarrow$, viewed as a symmetric monoidal category with its dilation action of $\Gm$.

\medskip

Regard $\ICoh(\Omega \bfA^1)$ as a symmetric monoidal category under convolution. Under the Koszul duality equivalence (\ref{eqn:KD}) for $V = \bfA^1$, the monoidal localization
    \[
    (j^*)^\Rightarrow : \QCoh(\bfA^1)^\Rightarrow \to \QCoh(\Gm)^\Rightarrow
    \]
goes over the monoidal localization
    \[
    \ICoh(\Omega \bfA^1) \to \ICohc(\Omega \bfA^1).
    \]
Therefore, we obtain an $\Gm$-equivariant equivalence of symmetric monoidal categories
    \[
    \ICohc(\Omega \bfA^1) \simeq \QCoh(\Gm)^\Rightarrow,
    \]
intertwining the \emph{inverse} dilation action on the LHS with the dilation action on the RHS. After inverting the $\Gm$ action on the RHS, we obtain a symmetric monoidal equivalence
    \[
    \ICohc(\Omega \bfA^1) \simeq \QCoh(\Gm)^\Leftarrow
    \]
which is equivariant for the dilation actions on both sides. 

\medskip

It follows that Construction \ref{construction:untwisted_action} admits the following reformulation:

\begin{construction}
\label{construction:untwisted_action_Koszul}
There is a natural lift of
\begin{equation}
\label{eqn:ICohc_OmegaA1_lift}
    \ICoh(Z) \otimes \ICohc(\Omega \bfA^1)
    \in \ICohc(\Omega\bfA^1)\MOD(\Gm\MOD).
\end{equation}
to an object of
    \[
    \left(
    \DMod(\Sing(Z)) \otimes \ICohc(\Omega \bfA^1)
    \right) \MOD
    (\Gm \MOD).
    \]
\end{construction}

\subsubsection{} We will carry out Construction \ref{construction:untwisted_action_Koszul} by studying the product
    \[
    Z \times \Omega \bfA^1,
    \]
equipped with the action of $\Gm$ by dilation on the second factor.

\medskip

To start off, $\bfP\Sing(Z \times \Omega \bfA^1)$ inherits a $\Gm$ action by the naturality of $\Sing(-)$, and we claim that
    \[
    \ICohc(Z \times \Omega \bfA^1)
    \in
    \Gm \MOD
    \]
naturally lifts to an object
\begin{equation}
\label{eqn:ICohc_SingZ_lift_prelim}
    \ICohc(Z \times \Omega \bfA^1)
    \in
    \DMod(\bfP\Sing(Z \times \Omega \bfA^1))\MOD(\Gm\MOD).
\end{equation}

\medskip

In other words, we are saying that
\begin{equation}
\label{eqn:ICohc_Z_OmegaA1_Gm}
    \ICohc((Z \times \Omega\bfA^1)/\Gm) 
    \simeq
    \ICohc(Z \times \Omega\bfA^1)^\Gm 
    \in
    \Rep(\Gm)\MOD
\end{equation}
admits a natural lift to
\begin{equation*}
\begin{split}
    & \ShCat(\bfP\Sing(Z \times \Omega \bfA^1)_\dR/\Gm) \\
    & \qquad \qquad \simeq
    \DMod(\bfP\Sing(Z \times \Omega \bfA^1))^\Gm\MOD
    (\Rep(\Gm)\MOD).
\end{split}
\end{equation*}
Indeed, applying \cite[Proposition-Construction 1.5.2]{arinkin2018category} to the quotient stack
    \[
    (Z \times \Omega \bfA^1)/\Gm
    \]
produces a lift of (\ref{eqn:ICohc_Z_OmegaA1_Gm}) to a sheaf of categories over
    \[
    \bfP\Sing((Z \times \Omega \bfA^1)/\Gm)_\dR
    \underset{((Z \times \Omega \bfA^1)/\Gm)_\dR}{\times}
    (Z \times \Omega \bfA^1)/\Gm,
    \]
and this last stack maps to
\begin{equation*}
\begin{split}
    & \bfP\Sing((Z \times \Omega \bfA^1)/\Gm)_\dR
    \underset{((Z \times \Omega \bfA^1)/\Gm)_\dR}{\times}
    (Z \times \Omega \bfA^1)_\dR/\Gm \\
    & \qquad \qquad \simeq
    \bfP\Sing(Z \times \Omega \bfA^1)_\dR/\Gm.
\end{split}
\end{equation*}

\subsubsection{} 
\label{subsubsection:PSing_Z_OmegaA1_action}
One has
    \[
    \Sing(Z \times \Omega \bfA^1)_\dR
    \simeq
    (\Sing(Z) \oplus \id)_\dR,
    \]
where the latter expression denotes the de Rham prestack of $\Sing(Z) \underset{Z}{\times} \bfA^1_Z$. 

\medskip

In terms of this identification, the action of $\Gm$ on $\Sing(Z \times \Omega \bfA^1)$ obtained from the dilation action of $\Gm$ on $\Omega \bfA^1$ by the naturality of the $\Sing(-)$ construction matches with \emph{inverse} dilation on $\bfA^1_Z$.

\medskip

Therefore, the open embedding
\begin{equation}
\label{eqn:FB}
    \Sing(Z) \mon \bfP(\Sing(Z) \oplus \id) \simeq \bfP\Sing(Z \times \Omega \bfA^1)
    \qquad
    \xi \mapsto [\xi : 1].
\end{equation}
intertwines the action of $\Gm$ on its source by dilation with the preceding action of $\Gm$ on its target, and restricting (\ref{eqn:ICohc_SingZ_lift_prelim}) along this map defines a lift of
\begin{equation*}
\begin{split}
    \ICoh(Z) 
    & \otimes \ICohc(\Omega \bfA^1) \\
    & \qquad \simeq
    \DMod(\Sing(Z))
    \underset{\DMod(\bfP(\Sing(Z) \oplus \id))}{\otimes}
    \ICohc(Z \times \Omega \bfA^1)
    \in
    \Gm \MOD
\end{split}
\end{equation*}
to an object
\begin{equation}
\label{eqn:ICohc_SingZ_lift}
    \ICoh(Z) \otimes \ICohc(\Omega \bfA^1)
    \in \DMod(\Sing(Z))\MOD(\Gm \MOD).
\end{equation}
We claim:

\begin{proposition}
\label{prop:two_actions_commute}
The $\Gm$-equivariant action $\DMod(\bfP \Sing(Z))$ on
    \[
    \ICoh(Z) \otimes \ICohc(\Omega \bfA^1)
    \in \Gm \MOD
    \]
encoded by (\ref{eqn:ICohc_SingZ_lift}) commutes with the $\Gm$-equivariant action of $\ICohc(\Omega \bfA^1)$ encoded by (\ref{eqn:ICohc_OmegaA1_lift}).
\end{proposition}

\subsubsection{} Observe that Construction \ref{construction:untwisted_action_Koszul} follows immediately from Proposition \ref{prop:two_actions_commute}. However, it is not clear \emph{a priori} that Proposition \ref{prop:two_actions_commute} is actually a proposition (as opposed to a construction). This will be shown in \S\ref{subsection:rel_crys_of_cat}, and the proof of Proposition \ref{prop:two_actions_commute} will appear in \S\ref{subsection:two_actions_commute_proof}.

\subsection{Relative crystals of categories}
\label{subsection:rel_crys_of_cat}

\subsubsection{} In this subsection, we will prove a general assertion about relative crystals of categories. A consequence of this result will be that the commutativity of the two actions appearing in Proposition \ref{prop:two_actions_commute} is merely a property and not additional data.

\subsubsection{} Let $\ZZ \to \YY$ be a schematic morphism of prestacks, and write
    \[
    p : \ZZ_{\dR/\YY} 
    = \YY_\dR \underset{\ZZ_\dR}{\times} \ZZ
    \to \YY
    \]
for its formal completion.

\begin{proposition}
\label{prop:rel_crys_locally_ff}
The functor
    \[
    \COIND_p : \ShCat(\ZZ_{\dR/\YY})
    \to
    \ShCat(\YY)
    \]
is locally fully faithful\footnote{A functor of $(\infty, 2)$-categories is said to be \emph{locally fully faithful} if it induces fully faithful functors on mapping categories.}.
\end{proposition}

\subsubsection{Proof} Note that it suffices to consider the case where $\YY$ is an affine scheme. Indeed, under the identifications
    \[
    \ShCat(\ZZ_{\dR/\YY}) \simeq \lim_S \ShCat((\ZZ_S)_{\dR/S})
    \qquad
    \ShCat(\YY) = \lim_S \ShCat(S),
    \]
where the limits are taken over all test affine schemes mapping to $\YY$, the functor $\COIND_p$ is the one induced by the functors
    \[
    \COIND_{p_S} : 
    \ShCat((\ZZ_S)_{\dR/S})
    \to
    \ShCat(S).
    \]

\subsubsection{} Therefore, we may replace $\ZZ \to \YY$ with a map of schemes $Z \to Y$, with $Y$ affine. Then $Y$ and $Z_{\dR/Y}$ are both 1-affine, and the $(\CORES, \COIND)$ adjunction identifies with
    \[
    \IND :
    \QCoh(Y)\MOD
    \rightleftharpoons
    \QCoh(Z_{\dR/Y})\MOD
    : \RES.
    \]
In order to show that $\RES$ is locally fully faithful, it suffices to check that the counit
    \[
    \IND \circ \RES(\CC) \to \CC
    \]
is a colocalization; i.e. that it admits a fully faithful left adjoint.

\medskip

The counit under consideration is given by the morphism of $\QCoh(Z_{\dR/Y})$-bimodules
\begin{equation}
\label{eqn:rel_DMod_diag}
\begin{split}
    & \QCoh(Z_{\dR/Y}) \underset{\QCoh(Y)}{\otimes} \QCoh(Z_{\dR/Y}) \\
    & \hspace{8em} \simeq
    \QCoh((Z \times_Y Z)_{\dR/Y})
    \overset{\Delta_{\dR/Y}^*}{\longrightarrow}
    \QCoh(Z_{\dR/Y}),
\end{split}
\end{equation}
so it suffices to show that $\Delta_{\dR/Y}^*$ admits a fully faithful $\QCoh(Z_{\dR/Y})$-bilinear left adjoint.

\medskip

By Kashiwara's lemma, the functor
    \[
    \Delta_\dR^* :
    \QCoh((Z \times_Y Z)_\dR)
    \to
    \QCoh(Z_\dR)
    \]
admits a fully faithful $\QCoh((Z \times_Y Z)_\dR)$-linear left adjoint. Tensoring up along
    \[
    - \underset{\QCoh(Y_\dR)}{\otimes}
    \QCoh(Y)
    \]
shows that (\ref{eqn:rel_DMod_diag}) admits a left adjoint with the desired properties.

\qed[Proposition \ref{prop:rel_crys_locally_ff}]

\subsubsection{Porism} Let us record some consequences of Proposition \ref{prop:rel_crys_locally_ff} (or, more precisely, its proof).

\medskip

Let $\CC$ be a sheaf of categories over $\ZZ_{\dR/\YY}$. Then the data of $\CC$ is encoded by the data of $\COIND_p(\CC)$ along with compatible choices of colocalizing subcategories
\begin{equation}
\label{eqn:rel_crys_subcat}
    \CC_S
    \mon
    \QCoh(S) \otimes \GAMMA(\YY; \COIND_p(\CC))
\end{equation}
for every test affine scheme mapping to $\YY$. 

\medskip

We note for future reference that the functor (\ref{eqn:rel_crys_subcat}) and its right adjoint are $\QCoh(S)$-linear.

\medskip

Let $\DD$ be another sheaf of categories over $\ZZ_{\dR/\YY}$. Then a morphism $\CC \to \DD$ in $\ShCat(\ZZ)$ amounts to a morphism
    \[
    F : \COIND_p(\CC) \to \COIND_p(\DD)
    \]
in $\ShCat(\YY)$ with the property: for every test affine scheme mapping to $\YY$, the functor
\begin{equation*}
\begin{split}
    \id_{\QCoh(S)} \otimes \GAMMA(S; F) :
    \QCoh(S) 
    & \otimes \GAMMA(S; \COIND_p(\CC)) \\
    & \qquad \to
    \QCoh(S) \otimes \GAMMA(S; \COIND_p(\DD))
\end{split}
\end{equation*}
takes $\CC_S$ into $\DD_S$.

\medskip

Let $\AA$ be a sheaf of monoidal categories over $\YY$. Then an action of $\CORES_p(\AA)$ on $\CC$ amounts to an action of $\AA$ on $\COIND_p(\CC)$ with the property: for every test affine scheme mapping to $\YY$, the action map
\begin{equation*}
\begin{split}
    & \QCoh(S) \otimes \GAMMA(S; \COIND_p(\CC)) \underset{\QCoh(S)}{\otimes} \GAMMA(S; \AA) \\
    & \qquad\qquad \simeq
    \QCoh(S) \otimes \GAMMA\left(S; \COIND_p(\CC) \underset{\YY}{\otimes} \AA \right) \\
    & \hspace{8em} \overset{\id \otimes \GAMMA(S; \act)}{\longrightarrow} 
    \QCoh(S) \otimes \GAMMA(S; \COIND_p(\CC))
\end{split}
\end{equation*}
restricts to a functor
    \[
    \CC_S \underset{\QCoh(S)}{\otimes} \GAMMA(S; \AA)
    \to
    \CC_S.
    \]

\subsection{Proof of Proposition \ref{prop:two_actions_commute}}
\label{subsection:two_actions_commute_proof}

\subsubsection{} We apply the paradigm of \S\ref{subsection:rel_crys_of_cat} in the case
    \[
    \ZZ 
    = \Sing(Z)/\Gm \subset
    \bfP(\Sing(Z) \oplus \id)/\Gm
    \qquad
    \YY = B \Gm.
    \]
We interpret (\ref{eqn:ICohc_SingZ_lift}) as making
    \[
    \CC^\Gm 
    = \left( \ICoh(Z) \otimes \ICohc(\Omega \bfA^1) \right)^\Gm
    \]
into a sheaf of categories over $\ZZ_{\dR/\YY}$, and (\ref{eqn:ICohc_OmegaA1_lift}) as encoding an action of
    \[
    \AA^\Gm = \ICohc(\Omega \bfA^1)^\Gm,
    \]
viewed as a constant sheaf of symmetric monoidal categories on $\YY$, on $\COIND_p(\CC^\Gm)$.

\medskip

In these terms, the assertion of Proposition \ref{prop:two_actions_commute} says that the action of $\AA^\Gm$ on $\COIND_p(\CC^\Gm)$ lifts to an action of $\CORES_p(\AA^\Gm)$ on $\CC^\Gm$. Since $\ZZ$ and $\YY$ are 1-affine, this amounts to the assertion: for every test affine scheme $S$ mapping to $\ZZ_{\dR/\YY}$, the full subcategory
\begin{equation}
\label{eqn:commuting_action_subcat}
    \QCoh(S) \underset{\QCoh(\ZZ_{\dR/\YY})}{\otimes} \CC^\Gm
    \subset
    \QCoh(S) \underset{\QCoh(\YY)}{\otimes} \CC^\Gm
\end{equation}
is stable under the action of $\AA^\Gm$.

\medskip

Observe that the inclusion (\ref{eqn:commuting_action_subcat}) rewrites as
    \[
    \left( 
    \QCoh(S') \underset{\DMod(\ZZ')}{\otimes} \CC
    \right)^\Gm
    \subset
    \big( 
    \QCoh(S') \otimes \CC
    \big)^\Gm,
    \]
where
    \[
    \ZZ' = \ZZ \underset{B\Gm}{\times} \pt
    \qquad
    S' = S \underset{B\Gm}{\times} \pt.
    \]
Therefore, it suffices to show: for every test affine scheme $T$ mapping to $\ZZ'_\dR$, the full subcategory
\begin{equation}
\label{eqn:commuting_action_subcat_deeq}
    \QCoh(T) \underset{\DMod(\ZZ')}{\otimes} \CC
    \subset
    \QCoh(T) \otimes \CC
\end{equation}
is stable under the action of $\AA$.

\subsubsection{} To proceed further, we will need to recall the description of (\ref{eqn:commuting_action_subcat_deeq}) from \cite[\S2.4]{arinkin2018category}. Write
    \[
    \Gamma_f 
    \mon T_\red \times \Sing(Z)
    \subset T_\red \times \bfP\Sing(Z \times \Omega \bfA^1)
    \]
for the graph of the associated morphism
    \[
    f : T_\red \to \Sing(Z) \subset \bfP(\Sing(Z) \oplus \id).
    \]
It corresponds to a closed conical subscheme
    \[
    \mathrm{Cone}(\Gamma_f) \mon T_\red \times \Sing(Z \times \Omega \bfA^1).
    \]
By definition, (\ref{eqn:commuting_action_subcat_deeq}) is the quotient of the full subcategory
\begin{equation}
\label{eqn:cone_subcat}
\begin{split}
    & \left(
    \QCoh(T) \otimes \ICoh(Z \times \Omega \bfA^1) 
    \right)_{\mathrm{Cone}(\Gamma_f)} \\
    & \hspace{8em} \subset
    \QCoh(T) \otimes \ICoh(Z \times \Omega \bfA^1) 
\end{split}
\end{equation}
by the full subcategory
    \[
    \QCoh(T) \otimes \QCoh(Z \times \Omega \bfA^1).
    \]
It remains for us to check that (\ref{eqn:cone_subcat}) is stable under the action of $\ICoh(\Omega \bfA^1)$.

\subsubsection{} Let us make the identification
\begin{equation}
\label{eqn:cone_subcat_big}
\begin{split}
    & \left(
    \QCoh(T) \otimes \ICoh(Z \times \Omega \bfA^1) 
    \right)_{\mathrm{Cone}(\Gamma_f)}
    \otimes \ICoh(\Omega \bfA^1) \\
    & \qquad \simeq
    \left(
    \QCoh(T) \otimes \ICoh(Z \times \Omega \bfA^1 \times \Omega \bfA^1) 
    \right)_{\mathrm{Cone}(\Gamma_f) \oplus \id}
\end{split}
\end{equation}
as full subcategories of
\begin{equation}
\label{eqn:T_ICoh_prod}
    \QCoh(T) \otimes \ICoh(Z \times \Omega \bfA^1 \times \Omega \bfA^1)
\end{equation}
In these terms, what remains to be seen is that the action map
\begin{equation*} 
\begin{split}
    \id_{\QCoh(T)} \otimes (\id \times \add)_* :
    \QCoh(T)
    & \otimes 
    \ICoh(Z \times \Omega \bfA^1 \times \Omega \bfA^1) \\
    & \qquad \qquad \to
    \QCoh(T) 
    \otimes 
    \ICoh(Z \times \Omega \bfA^1)
\end{split}
\end{equation*}
takes (\ref{eqn:cone_subcat_big}) into (\ref{eqn:cone_subcat}).

\subsubsection{} Consider the factorization of $\add$ into the composite 
\begin{equation}
\label{eqn:add_sigma_pr1}
    \Omega \bfA^1 \times \Omega \bfA^1
    \overset{\sigma}{\longrightarrow}
    \Omega \bfA^1 \times \Omega \bfA^1
    \overset{\pr_1}{\longrightarrow}
    \Omega \bfA^1
\end{equation} 
of the automorphism $\sigma$ defined by the formula
    \[
    (s, t) \mapsto (s + t, t)
    \]
and the co\"ordinate projection $\pr_1$. 

\medskip

The automorphism of
    \[
    S_\red \times \Sing(Z \times \Omega \bfA^1 \times \Omega \bfA^1)
    \]
induced by $\sigma$ restricts to an automorphism of\footnote{Recall that for a vector scheme $V$, $\Sing(\Omega V \times \Omega V)^{\cl} \simeq V^\vee \times V^\vee$, and duality exchanges shearing on the first factor with shearing on the second factor.} $\mathrm{Cone}(\Gamma_f) \oplus \id$, so 
    \[
    \id_{\QCoh(T)} \otimes (\id \times \sigma)_*
    \]
defines an automorphism of (\ref{eqn:T_ICoh_prod}) restricting to an automorphism of (\ref{eqn:cone_subcat_big}). 

\medskip

The functor
    \[
    \id_{\QCoh(T)} \otimes (\id \times \pr_1)_*
    \simeq
    \id_{\QCoh(T) \otimes \ICoh(Z \times \Omega\bfA^1)} \otimes \Gamma_{\ICoh}
    \]
takes (\ref{eqn:cone_subcat_big}) into (\ref{eqn:cone_subcat}).

\qed[Proposition \ref{prop:two_actions_commute}]

\section{Compatibilities}

\subsubsection{} We will show in Proposition \ref{prop:backwards_compat} that the action of $\DMod(\Sing(Z))^\Rightarrow$ on $\ICoh(Z)$ constructed in \S\ref{section:constr} has the expected backwards compatibility with the action of $\DMod(\bfP\Sing(Z))$ on $\ICohc(Z)$ constructed in \cite{arinkin2018category}. When there is possibility for confusion, we will refer to the former as the \emph{inhomogeneous action} and the latter as the \emph{homogeneous action}.

\medskip

The outline of this section is as follows:
\begin{itemize}
    \item In \S\ref{subsection:3-the-basics} we prove a few straightforward results about the action of $D(\Sing(Z))^\Rightarrow$ on $\ICoh(Z)$.
    \item In \S\ref{subsection:3-rel-d-mod} we record some results about the relation between $\ICoh$ on formal schemes and the action of $D(\Sing(Z))^\Rightarrow$. In particular we give a simple description of the category $\ICoh_0(Y^\wedge_Z)$ in terms of this action. The proofs of the results in this section are deferred to Subsections \ref{subsection:compat_iota_adjunction_thm} and \ref{subsection:compat_iota_adjunction_prop}.
    \item In \S\ref{subsection:compat_iota_adjunction_thm} we show how to reduce the statements in \S\ref{subsection:3-rel-d-mod} to \Cref{prop:iota_adjunction_linear}, a statement about the Arinkin--Gaitsgory homogeneous action of $D(\bfP \Sing(Z))$ on the category of singularities, whose proof is deferred to \S\ref{subsection:compat_iota_adjunction_prop}. 
    \item In \S\ref{subsection:compat_iota_adjunction_prop} we prove \Cref{prop:iota_adjunction_linear}. 
    \item In \S\ref{subsection:3-rel-AG} we prove a natural compatibility statement relating the inhomogeneous action and the homogeneous action. (Our proof of this compatibility takes \Cref{prop:iota_adjunction_linear} as input, which is why this subsection doesn't appear earlier in this section.)
    \item In \S\ref{subsection:global_complete_intersection} we give an alternative description of the inhomogeneous action when $Z$ is presented as a global complete intersection scheme. This version of the action was considered previously in \cite[Section 2.2]{beraldo2020spectral}, so this section gives a backwards compatibility between our construction and the construction of \textit{loc. cit.}
\end{itemize}

\subsection{The basics} \label{subsection:3-the-basics}

\subsubsection{} In this subsection we will collect some features of the inhomogenous action which are more or less immediate from its construction.

\begin{proposition}
\label{prop:D_action_compat}
Let $\pi : \Sing(Z) \to Z$ denote the canonical projection. The restriction of the inhomogeneous action along the symmetric monoidal functor
    \[
    (\pi^!)^\Rightarrow : \DMod(Z) \to \DMod(\Sing(Z))^\Rightarrow
    \]
recovers the usual action of $\DMod(Z)$ on $\ICoh(Z)$.
\end{proposition}

\subsubsection{Proof} This follows from \cite[Corollary 1.4.5(b)]{arinkin2018category}.

\subsubsection{} Now consider the zero section $Z \to \Sing(Z)$ and view $\DMod(Z)$ as tensored over $\DMod(\Sing(Z))$ using this map. We claim:

\begin{proposition}
\label{prop:ICoh_zero_compat}
The full subcategory
    \[
    \DMod(Z)
    \underset{\DMod(\Sing(Z))^\Rightarrow}{\otimes}
    \ICoh(Z)
    \subset
    \ICoh(Z)
    \]
is equal to
    \[
    \QCoh(Z) \subset \ICoh(Z).
    \]
\end{proposition}

\subsubsection{Proof} Recall the action of $\DMod(\Sing(Z))$ on 
    \[
    \ICoh(Z) \otimes \ICohc(\Omega \bfA^1)
    \]
from (\ref{eqn:ICohc_SingZ_lift}). We need to show that the full subcategory
\begin{equation*}
\begin{split}
    \DMod(Z) 
    \underset{\DMod(\Sing(Z))}{\otimes} 
    & \left(
    \ICoh(Z) \otimes \ICohc(\Omega \bfA^1)
    \right) \\
    & \qquad \subset
    \ICoh(Z) \otimes \ICohc(\Omega \bfA^1)
\end{split}
\end{equation*}
is equal to
    \[
    \QCoh(Z) \otimes \ICohc(\Omega \bfA^1)
    \subset 
    \ICoh(Z) \otimes \ICohc(\Omega \bfA^1).
    \]
This is a consequence of the following assertion, Proposition \ref{prop:prod_ss}.

\qed[Proposition \ref{prop:ICoh_zero_compat}]

\begin{proposition}
\label{prop:prod_ss}
Let $Z_1$ and $Z_2$ be schemes, and consider the closed embedding
    \[
    \bfP\Sing(Z_1) \times Z_2 \mon \bfP\Sing(Z_1 \times Z_2)
    \qquad
    \xi \mapsto [\xi : 0].
    \]
Then
\begin{equation*}
\begin{split}
    & \DMod(\bfP\Sing(Z_1) \times Z_2) \underset{\DMod(\bfP\Sing(Z_1 \times Z_2))}{\otimes} \ICohc(Z_1 \times Z_2) \\
    & \hspace{16em} \simeq
    \ICohc(Z_1) \otimes \QCoh(Z_2)
\end{split}
\end{equation*}
as full subcategories of
    \[
    \ICohc(Z_1 \times Z_2).
    \]
\end{proposition}

\subsubsection{Proof} This follows from \cite[Theorem-Construction 1.4.2(a)]{arinkin2018category} and \cite[Lemma 4.6.4]{arinkin2015singular}.

\qed[Proposition \ref{prop:prod_ss}]

\subsection{Relative D-modules} \label{subsection:3-rel-d-mod}

\subsubsection{} Let $f : Z \to Y$ be a morphism of quasi-smooth schemes. We will be interested in the category
    \[
    \ICoh(Y^\wedge_Z).
    \]
By \cite[Proposition 3.1.2]{arinkin2018category}, the functor of IndCoh $!$-pullback
    \[
    f^! : 
    \ICoh(Y) 
    \to 
    \ICoh(Y^\wedge_Z)
    \]
induces an equivalence of categories
\begin{equation}
\label{eqn:rel_DMod_tensor}
    \DMod(Z) \underset{\DMod(Y)}{\otimes} \ICoh(Y)
    \simto
    \ICoh(Y^\wedge_Z).
\end{equation}

\medskip

In view of Proposition \ref{prop:D_action_compat}, we may regard (\ref{eqn:rel_DMod_tensor}) as defining an equivalence of categories
\begin{equation}
\label{eqn:rel_DMod_tensor_Sing}
    \DMod(\Sing(Y)_Z)^\Rightarrow \underset{\DMod(\Sing(Y))^\Rightarrow}{\otimes} \ICoh(Y)
    \simto
    \ICoh(Y^\wedge_Z),
\end{equation}
where
    \[
    \Sing(Y)_Z = Z \underset{Y}{\times} \Sing(Y).
    \]
In particular, this equips $\ICoh(Y^\wedge_Z)$ an action of $\DMod(\Sing(Y)_Z)^\Rightarrow$. 

\subsubsection{} Now consider the natural map
    \[
    \iota : Z \to Y^\wedge_Z
    \]
and the resulting adjunction
\begin{equation}
\label{eqn:rel_DMod_iota_adjunction}
    \iota_* :
    \ICoh(Z) \rightleftharpoons \ICoh(Y^\wedge_Z)
    : \iota^!.
\end{equation}
Note that $\ICoh(Y^\wedge_Z)$ is naturally tensored over $\DMod(\Sing(Z))^\Rightarrow$ via the singular codifferential
    \[
    \Sing(f) :
    \Sing(Y)_Z \to \Sing(Z).
    \]
We will prove the following assertion over the course of \S\ref{subsection:compat_iota_adjunction_thm} and \S\ref{subsection:compat_iota_adjunction_prop}:

\begin{theorem}
\label{thm:iota_adjunction_linear}
The $(\iota_*, \iota^!)$ adjunction is naturally $\DMod(\Sing(Z))^\Rightarrow$-linear.
\end{theorem}

\subsubsection{} 
\label{subsubsection:ICoh0}
Recall that
    \[
    \ICoh_0(Y^\wedge_Z) \subset \ICoh(Y^\wedge_Z)
    \]
denotes the full subcategory of objects whose images under the functor $\iota^!$ of (\ref{eqn:rel_DMod_iota_adjunction}) belong to
    \[
    \QCoh(Z) \subset \ICoh(Z).
    \]
Recall also the `relative' scheme of singularities defined by the fiber product:
    \[
    \begin{tikzcd}
        \Sing(Z/Y) \arrow[r] \arrow[d] & Z \arrow[d] \\
        \Sing(Y)_Z \arrow[r, "\Sing(f)"] & \Sing(Z)
    \end{tikzcd}
    \]
One consequence of Theorem \ref{thm:iota_adjunction_linear} is:

\begin{corollary}
\label{cor:icoh-0-tensor-formula}
The equivalence (\ref{eqn:rel_DMod_tensor_Sing}) restricts to an equivalence
    \[
    \DMod(\Sing(Z/Y))^\Rightarrow 
    \underset{\DMod(\Sing(Y))^\Rightarrow}{\otimes}
    \ICoh(Y)
    \simto
    \ICoh_0(Y^\wedge_Z)
    \]
\end{corollary}

\subsection{Proof of Theorem \ref{thm:iota_adjunction_linear}}
\label{subsection:compat_iota_adjunction_thm}

\subsubsection{} In this subsection we will formulate Proposition \ref{prop:iota_adjunction_linear}, an analogue of Theorem \ref{thm:iota_adjunction_linear} for the homogenous action, and use it to deduce Theorem \ref{thm:iota_adjunction_linear} itself. The proof of Proposition \ref{prop:iota_adjunction_linear} will appear in \S\ref{subsection:compat_iota_adjunction_prop}.

\subsubsection{} Recall from \cite[Proposition 3.1.8]{arinkin2018category} that (\ref{eqn:rel_DMod_tensor}) restricts to an equivalence
\begin{equation}
\label{eqn:rel_DMod_tensor_circ}
    \DMod(Z) \underset{\DMod(Y)}{\otimes} \ICohc(Y)
    \simto
    \ICohc(Y^\wedge_Z)
\end{equation}
onto the full subcategory
    \[
    \ICohc(Y^\wedge_Z)
    \subset
    \ICoh(Y^\wedge_Z)
    \]
which is right-orthogonal to $\QCoh(Y^\wedge_Z)$.

\medskip

Let us view (\ref{eqn:rel_DMod_tensor_circ}) as tensored over $\DMod(\bfP\Sing(Y)_Z)$ using the homogeneous action of $\DMod(\bfP\Sing(Y))$ on $\ICohc(Y)$.

\medskip

Now the singular codifferential does not induce a morphism from the entirety of $\bfP\Sing(Y)_Z$ to $\bfP\Sing(Z)$. Instead, we only have a morphism
\begin{equation}
\label{eqn:P_singular_codifferential}
    \bfP\Sing(f) :
    \bfP \Sing(Y)_\Zc \to \bfP \Sing(\ZZ)
\end{equation}
out of the open subscheme
    \[
    \bfP \Sing(Y)_\Zc \mon \bfP \Sing(Y)
    \]
defined via the fiber product
    \[
    \begin{tikzcd}
        \Sing(Y)_{\Zc} \arrow[d, hook] \arrow[r] & \Singc(Z) \arrow[d, hook] \\
        \Sing(Y)_Z \arrow[r, "\Sing(f)"] & \Sing(Z)
    \end{tikzcd}
    \]
where $\Singc(Z)$ is the complement of the zero section in $\Sing(Z)$. According to \cite[Theorem 3.2.9]{arinkin2018category}, there is a natural identification 
    \[
    \ICohc_0(Y^\wedge_Z) 
    = \DMod(\bfP(\Sing(f)^{-1}(0))) 
    \underset{\DMod(\bfP \Sing(Y)_Z)}{\otimes} 
    \ICohc(Y^\wedge_Z)
    \]
as full subcategories of $\ICohc(Y^\wedge_Z)$. Combined with \cite[Proposition 3.2.3]{arinkin2018category} this shows that (\ref{eqn:rel_DMod_iota_adjunction}) restricts to an adjunction
\begin{equation}
\label{eqn:rel_DMod_iota_adjunction_circ}
    \ICohc(Z)
    \rightleftharpoons    
    \DMod(\bfP\Sing(Y)_\Zc)
    \underset{\DMod(\bfP\Sing(Y)_Z)}{\otimes}
    \ICohc(Y^\wedge_Z).
\end{equation}
Both categories appearing in this adjunction are naturally tensored over $\DMod(\bfP\Sing(Z))$; the LHS via the homogenous action and the RHS via the morphism (\ref{eqn:P_singular_codifferential}). The analogue of Theorem \ref{thm:iota_adjunction_linear} for the homogenous action is then:

\begin{proposition}
\label{prop:iota_adjunction_linear}
The adjunction (\ref{eqn:rel_DMod_iota_adjunction_circ}) is $\DMod(\bfP\Sing(Z))$-linear.
\end{proposition}

\subsubsection{Proof of Theorem \ref{thm:iota_adjunction_linear}, assuming Proposition \ref{prop:iota_adjunction_linear}} Consider the adjunction
\begin{equation}
\label{eqn:rel_DMod_Gm_adjunction}
\begin{split}
    & \iota_* \otimes \id : 
    \ICoh(Z) \otimes \QCoh(\Gm)^\Leftarrow \\
    & \hspace{8em} \rightleftharpoons
    \ICoh(Y^\wedge_Z) \otimes \QCoh(\Gm)^\Leftarrow
    : \iota^! \otimes \id,
\end{split}
\end{equation}    
taking place in
    \[
    \QCoh(\Gm)^\Leftarrow\MOD(\Gm\MOD).
    \]
Both categories appearing in this adjunction are equipped with lifts to
    \[
    \big(
    \DMod(\Sing(Z)) 
    \otimes \QCoh(\Gm)^\Leftarrow 
    \big)\MOD
    (\Gm \MOD),
    \]
and we are claiming that the adjunction itself lifts to this category. Note that Proposition \ref{prop:rel_crys_locally_ff} implies that this is a property and not additional structure, and that it suffices to check that (\ref{eqn:rel_DMod_Gm_adjunction}) is $\DMod(\Sing(Z))$-linear as an adjunction of plain dg categories. \\

Consider the adjunction (\ref{eqn:rel_DMod_iota_adjunction_circ}) for the map
    \[
    Z \times \Omega \bfA^1
    \to Y \times \Omega \bfA^1.
    \]
By Proposition \ref{prop:iota_adjunction_linear}, it respects the homogeneous action of
    \[
    \DMod(\bfP(\Sing(Z) \oplus \id)),
    \]
so after restriction along the open embedding (\ref{eqn:FB}), we obtain that the adjunction
    \[
    \ICoh(Z) \otimes \ICohc(\Omega \bfA^1)
    \rightleftharpoons
    \ICoh(Y^\wedge_Z) \otimes \ICohc(\Omega \bfA^1)
    \]
is linear over\footnote{Let us point out that the open subset $\Sing(Y) \subset \bfP(\Sing(Y) \oplus \id)$ is always contained in $\bfP \Sing(Y \times \Omega \bfA^1)_{\overset{\circ}{Z \times \Omega \bfA^1}}$.} $\DMod(\Sing(Z))$.

\qed[Theorem \ref{thm:iota_adjunction_linear}]

\subsection{Proof of Proposition \ref{prop:iota_adjunction_linear}}
\label{subsection:compat_iota_adjunction_prop}

\subsubsection{} Let $Y$ be a scheme and let $\CC$ be a category tensored over $\DMod(Y)$. So we are in the framework of \S\ref{subsection:rel_crys_of_cat}, but in the absolute case. In particular, for every test affine scheme $S$ mapping to $Y$, the colocalization (\ref{eqn:rel_crys_subcat}) becomes a $\QCoh(S)$-linear colocalization
\begin{equation}
\label{eqn:abs_crys_subcat}
    \CC_S = \QCoh(S) \underset{\DMod(Y)}{\otimes} \CC
    \rightleftharpoons
    \QCoh(S) \otimes \CC.
\end{equation}

\begin{proposition}
\label{prop:pointwise_subcat}
Suppose that $\CC$ is compactly generated and that $S$ is eventually coconnective. Let $c$ be an object of $\QCoh(S) \otimes \CC$. Then the following are equivalent:
\begin{itemize}
    \item[(a)] $c$ belongs to $\CC_S$;
    \item[(b)] for every geometric point $i : s \to S$, the functor
        \[
        i^* \otimes \id :
        \QCoh(S) \otimes \CC
        \to
        \QCoh(s) \otimes \CC
        \]
    takes $c$ into $\CC_s$.
\end{itemize}
\end{proposition}

\subsubsection{Proof} Note that $c$ belongs to $\CC_S$ if and only if the counit of (\ref{eqn:abs_crys_subcat}) is an isomorphism when evaluated on $c$. Therefore, the result follows from the fact that the collection of functors $\iota^* \otimes \id$ appearing in (b) are jointly conservative (for a proof of this assertion, see \cite[Lemma 21.4.6]{arinkin2020stack}).

\qed[Proposition \ref{prop:pointwise_subcat}]

\begin{corollary}
\label{cor:pointwise_linear}
Let $\CC$ and $\DD$ be compactly generated categories tensored over $\DMod(Y)$, and let
    \[
    F : \CC \to \DD
    \]
be a functor of plain dg categories. The following are equivalent:
\begin{itemize}
    \item[(a)] $F$ is $\DMod(Y)$-linear;
    \item[(b)] for every geometric point $\eta \to Y$ the functor
        \[
        \id \otimes F : 
        \QCoh(\eta) \otimes \CC
        \to
        \QCoh(\eta) \otimes \DD
        \]
    takes $\CC_\eta$ into $\DD_\eta$.
\end{itemize}
\end{corollary}

\subsubsection{Proof} This follows from Proposition \ref{prop:pointwise_subcat} and the fact that $Y_\dR$ is right Kan extended from a functor on classical schemes.

\qed[Corollary \ref{cor:pointwise_linear}]

\subsubsection{Proof of Proposition \ref{prop:iota_adjunction_linear}} Let $\eta$ be an arbitrary geometric point of $\bfP\Sing(Z)$. By Corollary \ref{cor:pointwise_linear}, it suffices to check that the adjunction
    \[
    \id_{\QCoh(\eta)} \otimes (\ref{eqn:rel_DMod_iota_adjunction_circ})
    \]
restricts to an adjunction
\begin{equation*}
\begin{split}
    & \left( 
    \QCoh(\eta) \otimes \ICohc(Z) 
    \right)_\eta \\
    & \qquad \qquad \rightleftharpoons
    \left(
    \QCoh(\eta) \otimes 
    \DMod(\bfP\Sing(Y)_\Zc)
    \underset{\DMod(\bfP\Sing(Y)_Z)}{\otimes}
    \ICohc(Y^\wedge_Z)
    \right)_\eta.
\end{split}
\end{equation*}
After extending the base field, we may assume that $\eta$ is a closed point. Then it remains to check that  (\ref{eqn:rel_DMod_iota_adjunction_circ}) restricts to an adjunction
\begin{equation*}
\begin{split}
    & \Vect 
    \underset{\DMod(\bfP\Sing(Z))}{\otimes}
    \ICohc(Z) \\
    & \qquad \qquad \rightleftharpoons
    \Vect \underset{\DMod(\bfP\Sing(Z))}{\otimes}
    \DMod(\bfP\Sing(Y)_\Zc)
    \underset{\DMod(\bfP\Sing(Y)_Z)}{\otimes}
    \ICohc(Y^\wedge_Z).
\end{split}
\end{equation*}
This follows from \cite[Theorem 3.2.9]{arinkin2018category} and \cite[Proposition 3.2.3]{arinkin2018category}.

\qed[Proposition \ref{prop:iota_adjunction_linear}]

\subsection{Relation to the action of Arinkin--Gaitsgory} \label{subsection:3-rel-AG}

\subsubsection{} To formulate the backwards compatability between the inhomogeneous and homogeneous actions, consider the full subcategory
    \[
    \DMod(\Singc(Z))^\Rightarrow \underset{\DMod(\Sing(Z))^\Rightarrow}{\otimes} \ICoh(Z)
    \subset
    \ICoh(Z).
    \]
By Proposition \ref{prop:ICoh_zero_compat}, it is equal to
    \[
    \ICohc(Z)
    \subset
    \ICoh(Z).
    \]
Therefore, we obtain an action of $\DMod(\Singc(Z))^\Rightarrow$ on $\ICohc(Z)$, and we may consider its restriction along the symmetric monoidal functor\footnote{We remind the reader that in this paper the notation $(-)^\Gm$ will always refer to \emph{weak} $\Gm$-invariants.}
    \[
    \DMod(\bfP\Sing(Z))
    \to
    \DMod(\Singc(Z))^\Gm
    \simeq
    \big(
    \DMod(\Singc(Z))^\Rightarrow
    \big)^\Gm
    \to 
    \DMod(\Singc(Z))^\Rightarrow.
    \]
We claim:

\begin{proposition}
\label{prop:backwards_compat}
This action agrees with the homogeneous action.
\end{proposition}

\subsubsection{Proof of Proposition \ref{prop:backwards_compat}} Consider the action of $\DMod(\Singc(Z))$ on the full subcategory 
\begin{equation*}
\begin{split}
    \DMod(\Singc(Z)) 
    \underset{\DMod(\Sing(Z))}{\otimes} 
    & \left(
    \ICoh(Z) \otimes \ICohc(\Omega \bfA^1)
    \right) \\
    & \qquad \subset
    \ICoh(Z) \otimes \ICohc(\Omega \bfA^1).
\end{split}
\end{equation*}
After taking $\Gm$-invariants, we obtain an action of $\DMod(\Singc(Z))^\Gm$ on
\begin{equation}
\label{eqn:backwards_compat_eq1}
\begin{split}
    & \left(
    \ICohc(Z) \otimes \ICohc(\Omega \bfA^1)
    \right)^\Gm \\
    & \qquad \qquad \simeq
    \left(
    \ICohc(Z) \otimes \QCoh(\Gm)^\Leftarrow
    \right)^\Gm
    \simeq
    \ICohc(Z).
\end{split}
\end{equation}
This is the action obtained by restricting the action of $\DMod(\Singc(Z))^\Rightarrow$ on $\ICohc(Z)$ along the functor
    \[
    \DMod(\Singc(Z))^\Gm
    \simeq
    \big(
    \DMod(\Singc(Z))^\Rightarrow
    \big)^\Gm
    \to
    \DMod(\Singc(Z))^\Rightarrow.
    \]

\medskip
    
Therefore, it suffices to show that the functor
\begin{equation}
\label{eqn:backwards_compat_eq2}
\begin{split}
    & \ICohc(Z) \otimes \ICohc(\Omega \bfA^1) \\
    & \qquad \to
    \left(
    \ICohc(Z) \otimes \ICohc(\Omega \bfA^1)
    \right)_\Gm \\
    & \qquad \qquad \simeq \left(
    \ICohc(Z) \otimes \ICohc(\Omega \bfA^1)
    \right)^\Gm 
    \overset{(\ref{eqn:backwards_compat_eq1})}{\longrightarrow}
    \ICohc(Z)
\end{split}
\end{equation}
intertwines the action of $\DMod(\Singc(Z))$ on its source with the homogeneous action of $\DMod(\bfP\Sing(Z))$ on its target, via the quotient map
\begin{equation}
\label{eqn:backwards_compat_quotient}
    \Singc(Z) \to \bfP \Sing(Z).
\end{equation}

\subsubsection{} Note that the functor
    \[
    \ICohc(\Omega \bfA^1)
    \to
    \ICohc(\Omega \bfA^1)_\Gm
    \simto 
    \Vect
    \]
identifies with the composite
    \[
    \ICohc(\Omega \bfA^1)
    \mon
    \ICoh(\Omega \bfA^1)
    \overset{i^!}{\longrightarrow}
    \Vect,
    \] 
where $i: \pt \to \Omega \bfA^1$ is the inclusion of zero into the vector scheme $\Omega \bfA^1$. Therefore, we may reinterpret (\ref{eqn:backwards_compat_eq2}) as the composite
    \[
    \ICohc(Z) \otimes \ICohc(\Omega \bfA^1)
    \mon
    \ICohc(Z \times \Omega \bfA^1)
    \overset{(\id \times i)^!}{\longrightarrow}
    \ICohc(Z).
    \]

\medskip

Now observe that $i$ is a nil-isomorphism, so the natural map
    \[
    (Z \times \Omega \bfA^1) ^\wedge_Z 
    \to
    Z \times \Omega \bfA^1
    \]
is an isomorphism. Therefore, Proposition \ref{prop:iota_adjunction_linear} implies that $(\id \times i)^!$ intertwines the action of $\DMod(\bfP\Sing(Z \times \Omega \bfA^1)_{\overset{\circ}{Z}})$ on $\ICohc(Z \times \Omega \bfA^1)$ with the action of $\DMod(\bfP\Sing(Z))$ on $\ICohc(Z)$, via the map
    \[
    \bfP\Sing(\id \times i) : 
    \bfP\Sing(Z \times \Omega \bfA^1)_{\overset{\circ}{Z}}
    \to
    \bfP\Sing(Z).
    \]
The result follows because the restriction of $\bfP\Sing(\id \times i)$ to the open subscheme
    \[
    \Singc(Z)
    \subset
    \bfP\Sing(Z \times \Omega \bfA^1)_{\overset{\circ}{Z}}
    \]
is the quotient map (\ref{eqn:backwards_compat_quotient}). 

\qed[Proposition \ref{prop:backwards_compat}]

\begin{proposition}
\label{prop:ICoh_N_compat}
For every closed conical subset $\Nilp \subset \Sing(Z)$, the full subcategory
    \[
    \DMod(\Nilp)^\Rightarrow \underset{\DMod(\Sing(Z))^\Rightarrow}{\otimes}
    \ICoh(Z)
    \subset
    \ICoh(Z)
    \]
is equal to
    \[
    \ICoh_\Nilp(Z) \subset \ICoh(Z).
    \]
\end{proposition}

\subsubsection{Proof} By \cref{prop:ICoh_zero_compat}, it suffices to show that the full subcategory
\begin{equation}
\label{eqn:ICoh_Nilp_tensor}
    \DMod(\bfP\Nilp)^\Rightarrow \underset{\DMod(\bfP\Sing(Z))^\Rightarrow}{\otimes}
    \ICohc(Z)
    \subset
    \ICohc(Z)
\end{equation}
is equal to
\begin{equation}
\label{eqn:ICohc_Nilp}
    \ICohc_\Nilp(Z) \subset \ICohc(Z).
\end{equation}
Although the tensor product in (\ref{eqn:ICoh_Nilp_tensor}) is formed using the inhomogenous action, \cref{prop:backwards_compat} implies that we may replace it with the homogenous action. Now the equality of (\ref{eqn:ICoh_Nilp_tensor}) and (\ref{eqn:ICohc_Nilp}) is \cite[Theorem-Construction 1.4.2(a)]{arinkin2018category}. 

\qed[Proposition \ref{prop:ICoh_N_compat}]

\subsection{Case of a global complete intersection}
\label{subsection:global_complete_intersection}

\subsubsection{} Suppose that $Z$ has been presented as a global complete intersection; that is, it has been written as a fiber product
\begin{equation}
\label{eqn:global_complete_intersection}
    \begin{tikzcd}
        Z \arrow[r] \arrow[d] & U \arrow[d] \\
        \pt \arrow[r] & V
    \end{tikzcd}
\end{equation}
where $U$ is smooth and $\pt \to V$ is the inclusion of the zero vector into a vector space. In this situation, it is possible to give a description of the action of $\DMod(\Sing(Z))^\Rightarrow$ on $\ICoh(Z)$ in the spirit of \cite[\S5.4]{arinkin2015singular}.

\begin{construction}
\label{construction:global_complete_intersection_action}
For every choice of presentation (\ref{eqn:global_complete_intersection}), one obtains an action of $\DMod(Z \times V^\vee)^\Rightarrow$ on $\ICoh(Z)$.
\end{construction}

\subsubsection{} Note that
\begin{equation}
\label{eqn:groupoid_parallel}
    Z \times_U Z 
    \simeq
    Z \times \Omega V
\end{equation}
as groupoids over $Z$, so we obtain an action of
    \[
    \ICoh(Z \times \Omega V)
    \simeq
    \ICoh(Z) \otimes \ICoh(\Omega V)
    \]
on $\ICoh(Z)$ by convolution. To obtain Construction \ref{construction:global_complete_intersection_action}, we precompose this action with the symmetric monoidal functor
\begin{equation}
\label{eqn:complete_intersection_oblv}
\begin{split}
    \DMod(Z \times V^\vee)^\Rightarrow 
    & \overset{\oblv^\Rightarrow}{\longrightarrow}
    \ICoh(Z \times V^\vee)^\Rightarrow \\
    & \qquad \qquad \simeq
    \ICoh(Z) \otimes \QCoh(V^\vee)^\Rightarrow \\
    & \hspace{8em} \simeq
    \ICoh(Z) \otimes \ICoh(\Omega V).
\end{split}
\end{equation}

\subsubsection{} Since $U$ is smooth, the morphism
    \[
    \Cot_Z
    \to
    \Cot_{Z/U} 
    \simeq \OO_Z \otimes V^\vee [1]
    \]
induces a closed embedding
\begin{equation}
\label{eqn:complete_intersection_sing_embed}
    \delta : \Sing(Z) \mon Z \times V^\vee.
\end{equation}
We claim:

\begin{proposition}
\label{prop:compat_complete_intersection}
The action of Construction \ref{construction:global_complete_intersection_action} agrees with composition of the inhomogeneous action with the monoidal colocalization
    \[
    (\delta^!)^\Rightarrow : 
    \DMod(Z \times V^\vee)^\Rightarrow 
    \to 
    \DMod(\Sing(Z))^\Rightarrow.
    \]
\end{proposition}

\subsubsection{A tautology}
\label{subsubsection:action_matching_tautology}
Let $\AA$ be a monoidal category for which the multiplication map
    \[
    \mult : \AA \otimes \AA \to \AA
    \]
admits a fully faithful $(\AA, \AA)$-bilinear left adjoint. As in \S\ref{subsection:rel_crys_of_cat}, this implies that the forgetful functor
    \[
    \AA \MOD \to \DGCAT
    \]
is locally fully faithful. 

\medskip

Suppose now that a category $\MM$ has been given two actions of $\AA$, to be denoted 
    \[
    \act_1 : \AA \otimes \MM \to \MM
    \qquad
    \act_2 : \AA \otimes \MM \to \MM.
    \]
Tensoring up the $(\mult^L, \mult)$ adjunction shows that $\act_1$ and $\act_2$ admit fully faithful $\AA$-linear left adjoints realizing $\MM$ as a full subcategory of $\AA \otimes \MM$ stable under the action of $\AA$. Therefore, the two actions of $\AA$ on $\MM$ coincide if and only if $\act_1$ and $\act_2$ are isomorphic as functors of plain dg categories.

\medskip

We claim that $\act_1$ and $\act_2$ are isomorphic as functors once $\act_1$ is $\AA$-linear for the tautological action of $\AA$ on $\AA \otimes \MM$ and the second action of $\AA$ on $\MM$. Indeed, this follows from precomposing the commutative diagram
    \[
    \begin{tikzcd}
        \AA \otimes \AA \otimes \MM \arrow[r, "\id \otimes \act_1"] \arrow[d, "\mult \otimes \id"] & \AA \otimes \MM \arrow[d, "\act_2"] \\
        \AA \otimes \MM \arrow[r, "\act_1"] & \MM
    \end{tikzcd}
    \]
with the map
    \[
    \id \otimes \unit \otimes \id : 
    \AA \otimes \MM \to \AA \otimes \AA \otimes \MM
    \]

\subsubsection{Proof of Proposition \ref{prop:compat_complete_intersection}} We apply the setup of \S\ref{subsubsection:action_matching_tautology} to the data
    \[
    \AA = \DMod(Z \times V^\vee)^\Rightarrow
    \qquad
    \MM = \ICoh(Z),
    \]
with $\act_1$ being the action of Construction \ref{construction:global_complete_intersection_action} and $\act_2$ being the composite of the inhomogenous action with $(\delta^!)^\Rightarrow$.

\medskip

To check that the action map 
    \[
    \act_1: \AA \otimes \MM \to \MM
    \]
is $D(Z \times V^\vee)^\Rightarrow$-linear, it suffices to check separately that it is linear for $D(Z)$ and for $D(V^\vee)^\Rightarrow$. 

\medskip

Since the action of Construction \ref{construction:global_complete_intersection_action} is $\DMod(Z)$-linear, the linearity of $\act_1$ for $D(Z)$ follows from \cref{prop:D_action_compat}. Therefore, it remains to show that $\act_1$ is $D(V^\vee)^\Rightarrow$-linear. 

\medskip

We claim that the functor
\begin{equation}
\label{eqn:complete_intersection_V_act}
    \ICoh(Z) \otimes \ICoh(\Omega V) 
    \to \ICoh(Z)
\end{equation}
of IndCoh $*$-pushforward along the action map
    \[
    \act : Z \times \Omega V \to Z
    \]
is already $D(V^\vee)^\Rightarrow$-linear in the appropriate sense. It's enough to show that (\ref{eqn:complete_intersection_V_act}) intertwines the action of $\DMod(V^\vee)^\Rightarrow \otimes \DMod(V^\vee)^\Rightarrow$ on its source with the action of $\DMod(V^\vee)^\Rightarrow$ on its target, via the map
    \[
    \Delta : V^\vee \to V^\vee \times V^\vee. 
    \]

\medskip

Indeed, this follows from Theorem \ref{thm:iota_adjunction_linear} for the map $\act$, together with the commutative diagram
    \[
    \begin{tikzcd}
        \Sing(Z) \arrow[r, "\Sing(\act)"] \arrow[d, "\delta"] & \Sing(Z) \times V^\vee \arrow[d, "\delta \times \id"] \\
        V^\vee \arrow[r, "\Delta"] & V^\vee \times V^\vee. 
    \end{tikzcd}
    \]

\qed[Proposition \ref{prop:compat_complete_intersection}]

\section{Morita equivalence}

In this section we shall formulate and prove the second main result of our paper, Theorem \ref{thm:morita_D_H}, which says that $\ICoh(Z)$ defines a Morita equivalence between $\DMod(\Sing(Z))^\Rightarrow$ and the monoidal category $\HH(Z)$.

\medskip

We then study the interaction of this equivalence with a morphism of quasi-smooth schemes $Z \to Y$.

\medskip 

The outline of this section is as follows:
\begin{itemize}
    \item In \S\ref{s:bimodule-structure} we describe how the action of $D(\Sing(Z))^\Rightarrow$ on $\ICoh(Z)$ behaves under duality, and as a consequence we show that $\ICoh(Z)$ upgrades to a $(D(\Sing(Z))^\Rightarrow, \HH(Z))$-bimodule.
    \item In \S\ref{s:cat-framework} we give a new proof of a theorem of D. Ben-Zvi, S. Gunningham, and H. Orem on Morita bimodules. Their result \cite[Theorem 1.3]{ben2020highest} appears as our \Cref{prop:fully_faithful_morita} and \Cref{cor:conservative_morita_equivalence}. 
    \item In \S\ref{subsection:4-our-sit-fits-in-fw} we show that in our situation the hypotheses of \Cref{prop:fully_faithful_morita} are satisfied. 
    \item In \S\ref{s:4-pf-of-morita-thm} we show that the additional hypothesis of \Cref{cor:conservative_morita_equivalence} is satisfied, proving \Cref{thm:morita_D_H}. 
    \item In \S\ref{s:4-transfer-bimodules} we calculate the image of the transfer bimodule $\HH_{Z \to Y}$ under our Morita equivalence. (Recall from \cite{beraldo2019sheaves} that $\HH_{Z \to Y}$ characterizes the functoriality of the category $\HH(Z)$ as $Z$ varies.)
    \item In \S\ref{s:4-comp-transf-bimod} we calculate the image of the composition map 
    \[\HH_{Z \leftarrow Y} \underset{\HH(Y)}{\otimes} \HH_{Y \leftarrow X} \to \HH_{Z \leftarrow X}\]
    under our Morita equivalence.
\end{itemize}

\subsection{Bimodule structure on the category of ind-coherent sheaves}\label{s:bimodule-structure}

\subsubsection{} 
\label{subsubsection:D_mod_duality}
Let $\MM$ be a module over
    \[
    \AA = \DMod(\Sing(Z))^\Rightarrow,
    \]
and suppose that $\MM$ admits a dual $\MM^\vee$ as a plain dg category. By \cite[\S A.3.9]{chen2023nearby}\footnote{Strictly speaking, \textit{loc. cit.} only considers categories of the form $D(Z)$, without shearing. However, the proof in \textit{loc. cit.} only relies on the category-theoretic properties of $D(Z)$ as an object of $\dgCat$. Since shearing is a symmetric-monoidal automorphism of $\Gm\MOD$ the assumptions listed in \textit{loc. cit.} are also satisfied in our situation.}, the composite
    \[
    \Vect 
    \to \MM \otimes \MM^\vee 
    \to \MM \otimes_\AA \MM^\vee
    \]
forms the unit of an $\AA$-module duality. 

\medskip

The counit of the duality
    \[
    \ev : \MM^\vee \otimes \MM \to \AA
    \]
is given by the composite
    \[
    \MM^\vee \otimes \MM
    \overset{\id \otimes \co\act}{\longrightarrow}
    \MM^\vee \otimes \MM \otimes \AA
    \to
    \AA,
    \]
where $\co\act$ is the morphism
    \[
    \co\act : 
    \MM
    \to
    \MM \otimes \AA^\vee
    \simto
    \MM \otimes \AA
    \]
obtained by transposing the action map and applying Verdier duality.

\medskip

We observe for future reference that the functor
    \[
    \AA \simeq
    \AA^\vee
    \overset{\ev^\vee}{\longrightarrow}
    (\MM^\vee \otimes \MM)^\vee
    \simeq
    \MM^\vee \otimes \MM
    \simeq
    \End(\MM)
    \]
is the map of $(\AA, \AA)$-bimodules classifying the action of $\AA$ on $\MM$.

\subsubsection{} Now let us consider the module $\MM = \ICoh(Z)$. In this case, Serre duality provides a canonical equivalence of dg categories
\begin{equation}
\label{eqn:Serre_duality}
    \ICoh(Z)^\vee \simeq \ICoh(Z).
\end{equation}
An irksome feature of the theory is that this equivalence \emph{does not} intertwine the $\AA$-module structure on the LHS from (\ref{subsubsection:D_mod_duality}) with the pre-existing $\AA$-module structure on the RHS. 

\medskip

To correct this, let
    \[
    \ICoh(Z)^\tau
    \]
denote the category $\ICoh(Z)$ equipped with the action of $\AA$ obtained by precomposing the standard action with the involution of $\AA$ induced by the negation map: 
    \[
    \tau : \Sing(Z) \to \Sing(Z)
    \qquad
    \xi \mapsto -\xi.
    \]

\begin{proposition}
\label{prop:Serre_duality_tau}
The equivalence of plain dg categories (\ref{eqn:Serre_duality}) lifts to an an equivalence of $\DMod(\Sing(Z))^\Rightarrow$-modules
    \[
    \ICoh(Z)^\vee \simeq \ICoh(Z)^\tau.
    \]
\end{proposition}

\subsubsection{Proof} We need to lift the map of left $\DMod(\Sing(Z))^\Rightarrow$-modules
\begin{equation}
\label{eqn:coev_Serre}
\begin{split}
    \DMod(\Sing(Z))^\Rightarrow
    & \overset{\coev^\vee}{\too}
    \ICoh(Z) \otimes \ICoh(Z)^\vee \\
    & \qquad \qquad \simeq
    \ICoh(Z) \otimes \ICoh(Z)^\tau
\end{split}
\end{equation}
to a map of bimodules. It suffices to show that the functor (\ref{eqn:coev_Serre}) factors through the full subcategory
\begin{equation}
\label{eqn:tau_tensor_subcat}
\begin{split}
    \DMod(\Sing(Z))^\Rightarrow 
    \underset{\DMod(\Sing(Z))^\Rightarrow \otimes \DMod(\Sing(Z))^\Rightarrow}{\otimes}
    & \big( \ICoh(Z) \otimes \ICoh(Z)^\tau \big) \\
    & \qquad \qquad \subset
    \ICoh(Z) \otimes \ICoh(Z)^\tau,
\end{split}
\end{equation}
where $\DMod(\Sing(Z))^\Rightarrow \otimes \DMod(\Sing(Z))^\Rightarrow$ acts on $\DMod(\Sing(Z))^\Rightarrow$ via the diagonal embedding
    \[
    \Sing(Z) \mon \Sing(Z) \times \Sing(Z).
    \]

\medskip

Observe that (\ref{eqn:coev_Serre}) is given by acting on the object
    \[
    \Delta_*(\omega_Z) \in \ICoh(Z \times Z),
    \]
and (\ref{eqn:tau_tensor_subcat}) is a submodule inclusion. Therefore, it suffices to show that $\Delta_*(\omega_Z)$ belongs to the submodule (\ref{eqn:tau_tensor_subcat}).

\subsubsection{} Let us rewrite the full subcategory (\ref{eqn:tau_tensor_subcat}) as
\begin{equation}
\label{eqn:tau_tensor_subcat_rewrite}
\begin{split}
    \DMod(\Sing(Z))^\Rightarrow 
    \underset{\DMod(\Sing(Z \times Z))^\Rightarrow}{\otimes}
    & \ICoh(Z \times Z) \\
    & \qquad \qquad \subset
    \ICoh(Z \times Z),
\end{split}
\end{equation}
where $\DMod(\Sing(Z \times Z))^\Rightarrow$ now acts on $\ICoh(Z \times Z)$ in the usual way and on $\DMod(\Sing(Z))^\Rightarrow$ via the closed embedding
    \[
    \id \times \tau : 
    \Sing(Z) \to \Sing(Z) \times \Sing(Z) \simeq \Sing(Z \times Z).
    \]
This embedding identifies with the morphism
    \[
    \Sing(Z/Z \times Z) \to \Sing(Z \times Z),
    \]
so Corollary \ref{cor:icoh-0-tensor-formula} identifies (\ref{eqn:tau_tensor_subcat_rewrite}) with the full subcategory
    \[
    \ICoh_0((Z \times Z)^\wedge_Z)
    \subset
    \ICoh(Z \times Z)
    \simeq
    \ICoh(Z) \otimes \ICoh(Z)^\tau.
    \]

\qed[Proposition \ref{prop:Serre_duality_tau}]

\subsubsection{Porism}
\label{subsubsection:bimod_construction} 
The proof of Proposition \ref{prop:Serre_duality_tau} shows that
    \[
    \End_{\DMod(\Sing(Z))^\Rightarrow}(\ICoh(Z))
    \simeq
    \HH(Z)
    \]
as full monoidal subcategories of
    \[
    \End(\ICoh(Z)) \simeq \ICoh(Z \times Z),
    \]
acting on the \emph{right} on $\ICoh(Z)$.

\medskip

In particular, we have upgraded $\ICoh(Z)$ to a bimodule
    \[
    \ICoh(Z)
    \in
    (\DMod(\Sing(Z))^\Rightarrow, \HH(Z)) \BIMOD.
    \]

\begin{theorem}
\label{thm:morita_D_H}
Let $Z$ be a quasi-smooth scheme. The bimodule $\ICoh(Z)$ defines a Morita equivalence
    \[
    \HH(Z)\MOD
    \simto
    \DMod(\Sing(Z))^\Rightarrow\MOD
    \]
\end{theorem}

\subsubsection{} The symmetric monoidal structure on $\DMod(\Sing(Z))^\Rightarrow$ gives a canonical identification of its categories of left and right modules. We will use this identification freely.

\medskip

The monoidal involution
    \[
    \sigma : 
    \HH(Z)^\rev
    \simto
    \HH(Z)
    \]
induced by the swap map on $Z \times Z$ gives a canonical identification between the categories of left and right modules over $\HH(Z)$. This identification exchanges $\ICoh(Z)$, viewed as a right $\HH(Z)$-module under convolution, with $\ICoh(Z)$, viewed as a left $\HH(Z)$-module under convolution. We will use this identification freely.

\medskip 

In particular, we may regard $\ICoh(Z)^\tau$ as a bimodule
\begin{equation}
\label{eqn:ICoh_tau_bimod}
    \ICoh(Z)^\tau
    \in
    (\HH(Z), \DMod(\Sing(Z))^\Rightarrow)\BIMOD.
\end{equation}

\subsubsection{}
\label{subsubsection:morita_inverse_tau}
By definition, the equivalence of Theorem \ref{thm:morita_D_H} fits into the commutative diagrams:
\begin{equation*}
\begin{tikzcd}[column sep={4.5em, between origins}]
    \HH(Z) \MOD\arrow[dr, "\ICoh(Z)" below left] \arrow[rr, "\sim"] & & \DMod(\Sing(Z))^\Rightarrow \MOD  \arrow[dl, "\OBLV"] \\
    & \DGCAT &
\end{tikzcd}
\qquad
\begin{tikzcd}[column sep={4.5em, between origins}]
    \HH(Z) \MOD \arrow[rr, "\sim"] & & \DMod(\Sing(Z))^\Rightarrow \MOD \\
    & \DGCAT \arrow[ur, "\ICoh(Z)" below right] \arrow[ul, "\IND" below left] &
\end{tikzcd}
\end{equation*}
After passing to left (resp. right) adjoints, we obtain the diagrams:
\begin{equation*}
\begin{tikzcd}[column sep={4.5em, between origins}]
    \HH(Z) \MOD \arrow[rr, "\sim"] & & \DMod(\Sing(Z))^\Rightarrow \MOD \\
    & \DGCAT \arrow[ur, "\IND" below right] \arrow[ul, "\ICoh(Z)" below left] &
\end{tikzcd}
\qquad
\begin{tikzcd}[column sep={4.5em, between origins}]
    \HH(Z) \MOD\arrow[dr, "\OBLV" below left] \arrow[rr, "\sim"] & & \DMod(\Sing(Z))^\Rightarrow \MOD  \arrow[dl, "\ICoh(Z)^\tau"] \\
    & \DGCAT &
\end{tikzcd}
\end{equation*}
This shows in particular that the inverse to the Morita equivalence of Theorem \ref{thm:morita_D_H} is given by the bimodule (\ref{eqn:ICoh_tau_bimod}).

%\subsubsection{Remark} The symmetric monoidal structure on $\DMod(\Sing(Z))^\Rightarrow$ induces a symmetric monoidal structure on $\DMod(\Sing(Z))^\Rightarrow \MOD$. On the other hand, one can give $\HH(Z)$ the structure of a cocommutative Hopf algebra, so $\HH(Z)\MOD$ is also equipped with a canonical symmetric monoidal structure. Then the equivalence of Theorem \ref{thm:morita_D_H} is symmetric monoidal---we will not pursue this point further in this work.

\subsection{Categorical framework for the argument}\label{s:cat-framework}

\subsubsection{} The assertion of Theorem \ref{thm:morita_D_H} takes place in context of the Morita $(\infty, 2)$-category of algebras and bimodules. For the convenience of the reader, we will recall the basic features of this theory which are relevant to us.

\subsubsection{} Let $\AA$ and $\BB$ be monoidal categories, and let $\MM$ be a $(\BB, \AA)$-bimodule. Suppose that $\MM$ admits a left dual $\MM^\vee$, i.e. an $(\AA, \BB)$-bimodule equipped with a pair of bimodule morphisms
\begin{equation}
\label{eqn:bimodule_duality}
    \coev : \BB \to \MM \otimes_\AA \MM^\vee
    \qquad
    \ev : \MM^\vee \otimes_\BB \MM \to \AA
\end{equation}
satisfying the triangle identities.

\subsubsection{} 
\label{subsubsection:right_adjointable_duality}
Suppose that the morphisms appearing in the duality data (\ref{eqn:bimodule_duality}) admit continuous right adjoints, which are moreover linear for the apparent bimodule structures.
    \[
    \BB \leftarrow \MM \otimes_\AA \MM^\vee : \coev^R
    \qquad
    \MM^\vee \otimes_\BB \MM \leftarrow \AA : \ev^R.
    \]
Since passing to adjoints is functorial, these right adjoints identify $\MM^\vee$ as a right dual of $\MM$.

\begin{proposition}
\label{prop:fully_faithful_morita}
Suppose in the situation of \S\ref{subsubsection:right_adjointable_duality} that $\coev$ is an equivalence. Then $\ev$ is fully faithful.
\end{proposition}

\subsubsection{Proof} It suffices\footnote{The counit of any coalgebra structure on the monoidal unit of a monoidal category $\AA$ is an isomorphism. Indeed, counitality says that $\comult : \id \to \id \otimes \id \simeq \id$ and $\counit : \id \to \id$ satisfy $\counit \circ \comult = \mathrm{id}$. Since $\End_\AA(\id)$ is an $\mathbf{E}_2$-algebra, this implies $\counit$ and $\comult$ are mutually inverse.

Viewing the comonad $\ev^R \circ \ev$ as a coalgebra in the endofunctor category of $\MM^\vee \otimes_\BB \MM$ shows that the counit of the $(\ev, \ev^R)$ is an equivalence once we know that $\ev^R \circ \ev$ happens to be isomorphic to $\id$.} to show that the morphism
    \[
    \MM^\vee \otimes_\BB \MM 
    \overset{\ev}{\to} 
    \AA
    \overset{\ev^R}{\to}
    \MM^\vee \otimes_\BB \MM
    \]
is isomorphic to the identity
    \[
    \MM^\vee \otimes_\BB \MM 
    \overset{\id}{\to}
    \MM^\vee \otimes_\BB \MM.
    \]
After transposition, we see that it is equivalent to compare the morphisms
    \[
    \id \otimes \ev^R : 
    \MM 
    \to 
    \MM \otimes_\AA \MM^\vee \otimes_\BB \MM
    \]
and
    \[
    \coev \otimes \id : 
    \MM 
    \to 
    \MM \otimes_\AA \MM^\vee \otimes_\BB \MM.
    \]
Since $\coev$ is an equivalence, so is $\coev^R$. Therefore it suffices to see that both morphisms are the identity after composition with $\coev^R \otimes \id$. For the first morphism this is a triangular identity. For the second morphism we use $\coev^R \simeq \coev^{-1}$.

\qed[Proposition \ref{prop:fully_faithful_morita}]

\medskip 

Note that Proposition \ref{prop:fully_faithful_morita} provides a more direct proof of \cite[Theorem 1.3]{ben2020highest}:

\begin{corollary}
\label{cor:conservative_morita_equivalence}
Suppose in the situation of \S\ref{subsubsection:right_adjointable_duality} that $\coev$ is an equivalence. Then $\MM$ defines a Morita equivalence if and only if $\ev^R$ is conservative.
\end{corollary}

\subsection{Our situation fits into the framework} \label{subsection:4-our-sit-fits-in-fw}

\subsubsection{} The goal of this subsection is to show that the hypothesis of Corollary \ref{cor:conservative_morita_equivalence} is applicable to the data
\begin{equation}
\label{eqn:A_B_M}
    \AA = \DMod(\Sing(Z))^\Rightarrow
    \qquad
    \BB = \HH(Z)
    \qquad
    \MM = \ICoh(Z).
\end{equation}
That is:
\begin{itemize}
    \item $\MM$ admits a left dual as a $(\BB, \AA)$-bimodule (Proposition \ref{prop:M_left_dual}),
    \item the coevaluation map for this duality is an equivalence (see \S\ref{subsubsection:bimod_construction}), and
    \item the evaluation map admits an $(\AA, \AA)$-bilinear right adjoint (Proposition \ref{prop:ev_continuous_right_adjoint}).
\end{itemize}

\begin{proposition}
\label{prop:M_left_dual}
In the situation of (\ref{eqn:A_B_M}), $\MM$ admits a left dual as a $(\BB, \AA)$-bimodule.
\end{proposition}

\subsubsection{Proof} We need to show that the functor
    \[
    \AA \MOD \to \BB \MOD
    \]
admits a left adjoint. Indeed, its composite with the forgetful functor
    \[
    \oblv : \BB \MOD \to \dgCat
    \]
admits a left adjoint because $\MM$ is dualizable as an $\AA$-module, and $\BB \MOD$ is generated under colimits by the image of $\ind$.

\qed[Proposition \ref{prop:M_left_dual}]

\subsubsection{} Note that the proof of Proposition \ref{prop:M_left_dual} identifies the dg category underlying the left dual of $\MM$ with its dual $\MM^\vee$ as a plain dg category. In these terms, the resulting duality data
    \[
    \coev : \BB \to \MM \otimes_\AA \MM^\vee
    \qquad
    \ev : \MM^\vee \otimes_\BB \MM \to \AA
    \]
can be described as follows:
\begin{itemize}
    \item[(i)] The composite
        \[
        \BB
        \overset{\coev}{\longrightarrow} \MM \otimes_\AA \MM^\vee
        \simeq \End_\AA(\MM)
        \]
    is the action map. In the case at hand, we saw in \S\ref{subsubsection:bimod_construction} that this map in equivalence.

    \item[(ii)] The composite
    \begin{equation}
    \label{eqn:mu_ev}
        \MM^\vee \otimes \MM
        \to \MM^\vee \otimes_\BB \MM
        \overset{\ev}{\longrightarrow} \AA
    \end{equation}
    is the counit of the $\AA$-module duality between $\MM^\vee$ and $\MM$.
\end{itemize}

\begin{proposition}
\label{prop:ev_continuous_right_adjoint}
In the situation of (\ref{eqn:A_B_M}), $\ev$ admits a continuous right adjoint, which is moreover $(\AA, \AA)$-bilinear.
\end{proposition}

\subsubsection{Proof} Since $\BB$ is rigid, the functor
    \[
    \mu : 
    \MM^\vee \otimes \MM
    \to \MM^\vee \otimes_\BB \MM
    \]
admits a continuous right adjoint given by the tensor product
    \[
    \mult^R \otimes \id : 
    \BB \underset{\BB \otimes \BB^\rev}{\otimes} (\MM^\vee \otimes \MM)
    \to (\BB \otimes \BB^\rev) \underset{\BB \otimes \BB^\rev}{\otimes} (\MM^\vee \otimes \MM).
    \]
Now $\mu^R$ is visibly $(\AA, \AA)$-bilinear, and it is conservative because the image of $\mu$ generates its target. Therefore, it suffices to check that the composite (\ref{eqn:mu_ev}) admits a continuous $(\AA, \AA)$-bilinear right adjoint. 

\medskip

Recall from \S\ref{subsubsection:D_mod_duality} that the composite (\ref{eqn:mu_ev}) identifies with the functor dual to the action map
\begin{equation}
\label{eqn:A_M_act}
    \AA \to \MM \otimes \MM^\vee.
\end{equation}
Therefore, it suffices show that (\ref{eqn:A_M_act}) admits an $(\AA, \AA)$-bilinear \emph{left} adjoint.

\subsubsection{} After Zariski localization, we may assume that $Z$ is presented as a global complete intersection (\ref{eqn:global_complete_intersection}). Recall the closed immersion (\ref{eqn:complete_intersection_sing_embed}) and the identification (\ref{eqn:groupoid_parallel})
    \[
    \delta : \Sing(Z) \mon Z \times V^\vee
    \qquad
    Z \times \Omega V 
    \simeq Z \times_U Z
    \]
from \S\ref{subsection:global_complete_intersection}.

\medskip

Since $(\delta_!, \delta^!)$ adjunction is a colocalization, it suffices to check that the composite of monoidal functors
\begin{equation}
\label{eqn:delta_act}
    \DMod(Z \times V^\vee)^\Rightarrow
    \overset{(\delta^!)^\Rightarrow}{\longrightarrow} 
    \DMod(\Sing(Z))^\Rightarrow
    \overset{(\ref{eqn:A_M_act})}{\longrightarrow}
    \ICoh(Z \times Z)
\end{equation}
admits a $\DMod(Z \times V^\vee)^\Rightarrow$-bilinear left adjoint. 

\medskip

By Proposition \ref{prop:compat_complete_intersection}, we may replace (\ref{eqn:delta_act}) with the sequence of monoidal functors:
\begin{align*}
    & \DMod(Z) \otimes \DMod(V^\vee)^\Rightarrow
    \overset{\oblv \otimes \oblv^\Rightarrow}{\longrightarrow}
    \ICoh(Z) \otimes \QCoh(V^\vee)^\Rightarrow \\
    & \qquad \simto
    \ICoh(Z) \otimes \ICoh(\Omega V)
    \simto
    \ICoh(Z \times_U Z)
    \overset{(\Delta_{Z/U})_*}{\longrightarrow}
    \ICoh(Z \times Z),
\end{align*}
where the last functor is that of IndCoh $*$-pushforward along the relative diagonal
    \[
    \Delta_{Z/U} : Z \times_U Z \to Z \times Z.
    \]

\subsubsection{} The left adjoint
    \[
    \ind \otimes \ind^\Rightarrow : 
    \ICoh(Z) \otimes \QCoh(V^\vee)^\Rightarrow
    \to 
    \DMod(Z) \otimes \DMod(V^\vee)^\Rightarrow
    \]
of $\oblv \otimes \oblv^\Rightarrow$ is $\DMod(Z) \otimes \DMod(V^\vee)^\Rightarrow$-bilinear. Therefore, it remains to show that $(\Delta_{Z/U})_*$ admits an $\ICoh(Z \times_U Z)$-bilinear left adjoint.

\medskip

Now $\Delta_{Z/U}$ is quasi-smooth because it is a base change of the diagonal
    \[
    U \to U \times U.
    \]
In particular, it has finite tor amplitude and admits an IndCoh $*$-pullback. The requisite bilinearity follows from a diagram chase using the base change compatibilities of IndCoh $*$-pullback against IndCoh $*$-pushforward and IndCoh $!$-pullback.

\qed[Proposition \ref{prop:ev_continuous_right_adjoint}]

\subsection{Proof of Theorem \ref{thm:morita_D_H}}\label{s:4-pf-of-morita-thm}

\subsubsection{} By Corollary \ref{cor:conservative_morita_equivalence}, it remains to check that the functor
    \[
    \ev^R : \AA \to \MM^\vee \otimes_\BB \MM
    \]
is conservative. It suffices to prove that the composite
    \[
    \AA \overset{\ev^R}{\to} \MM^\vee \otimes_\BB \MM \overset{\mu^R}{\to} \MM^\vee \otimes \MM
    \]
is conservative. As before, we may Zariski localize and assume that $Z$ is presented as a global complete intersection. 

\subsubsection{} The description of
    \[
    (\mu \circ \ev)^{R, \vee} \simeq (\mu \circ \ev)^{\vee, L}
    \]
from the proof of Proposition \ref{prop:ev_continuous_right_adjoint} identifies $\mu^R \circ \ev^R$ with
\begin{align*}
    & \DMod(\Sing(Z))^\Rightarrow \\
    & \quad \overset{(\delta^!)^{\Rightarrow, \vee}}{\longrightarrow}
    \DMod(Z) \otimes \DMod(V^\vee)^\Rightarrow
    \overset{\ind^\vee \otimes (\ind^\Rightarrow)^\vee}{\longrightarrow}
    \ICoh(Z) \otimes \QCoh(V^\vee)^\Rightarrow \\
    & \qquad \simto
    \ICoh(Z) \otimes \ICoh(\Omega V)
    \simto
    \ICoh(Z \times_U Z)
    \overset{(\Delta_{Z/U}^*)^\vee}{\longrightarrow}
    \ICoh(Z \times Z).
\end{align*}

\medskip

Now, $\Delta_{Z/U}^*$ differs from $\Delta_{Z/U}^!$ by a shifted line bundle because $\Delta_{Z/U}$ is quasi-smooth and hence Gorenstein. Therefore, it is equivalent to show that the functor
\begin{equation}
\label{eqn:complete_intersection_act_oblv_push}
\begin{split}
    & \DMod(\Sing(Z))^\Rightarrow
    \overset{(\delta_*)^\Rightarrow}{\longrightarrow}
    \DMod(Z) \otimes \DMod(V^\vee)^\Rightarrow
    \overset{\oblv \otimes \oblv}{\longrightarrow}
    \ICoh(Z) \otimes \QCoh(V^\vee)^\Rightarrow \\
    & \qquad \simto
    \ICoh(Z) \otimes \ICoh(\Omega V)
    \simto
    \ICoh(Z \times_U Z) \\
    & \qquad \qquad \overset{(\Delta_{Z/U})_*}{\longrightarrow}
    \ICoh(Z \times Z)
\end{split}
\end{equation}
is conservative.

\subsubsection{} We shall analyse (\ref{eqn:complete_intersection_act_oblv_push}) by factoring into the composite of
\begin{equation}
\label{eqn:complete_intersection_act_oblv}
\begin{split}
    & \DMod(\Sing(Z))^\Rightarrow 
    \overset{(\delta_*)^\Rightarrow}{\longrightarrow}
    \DMod(Z) \otimes \DMod(V^\vee)^\Rightarrow
    \overset{\oblv \otimes \oblv^\Rightarrow}{\longrightarrow}
    \ICoh(Z) \otimes \QCoh(V^\vee)^\Rightarrow \\
    & \qquad \simto
    \ICoh(Z \times \Omega V),
\end{split}
\end{equation}
followed by the functor of IndCoh $*$-pushforward along the map
    \[
    \alpha : 
    Z \times \Omega V 
    \simeq
    Z \times_U Z
    \to
    Z \times Z.
    \]

\subsubsection{} Observe that the functor (\ref{eqn:complete_intersection_act_oblv}) is conservative. Furthermore, by writing it as the composite
\begin{equation*}
\begin{split}
    \DMod(\Sing(Z))^\Rightarrow 
    \overset{(\delta_*)^\Rightarrow}{\longrightarrow}
    & \DMod(Z) \otimes \DMod(V^\vee)^\Rightarrow 
    \overset{\oblv \otimes \oblv^\Rightarrow}{\longrightarrow}
    \QCoh(Z) \otimes \QCoh(V^\vee)^\Rightarrow \\
    & \qquad \mon
    \ICoh(Z) \otimes \QCoh(V^\vee)^\Rightarrow
    \simto
    \ICoh(Z \times \Omega V),
\end{split}
\end{equation*}
we see that (\ref{eqn:complete_intersection_act_oblv}) intertwines the action of $\DMod(\Sing(Z))^\Rightarrow$ on its source with the action of $\DMod(\Sing(Z \times \Omega V))^\Rightarrow$ on its target, via the map
\begin{equation}
\label{eqn:0_delta}
    \Sing(Z) 
    \overset{\delta}{\mon} 
    Z \times V^\vee
    \mon
    \Sing(Z) \times V^\vee
    \simeq
    \Sing(Z \times \Omega V).
\end{equation}    
It follows that the image of (\ref{eqn:complete_intersection_act_oblv}) is contained in the full subcategory of $\ICoh(Z \times \Omega V)$ consisting of objects whose singular support is contained in the image of (\ref{eqn:0_delta}).

\medskip

By \cite[Proposition 7.3.10]{arinkin2015singular}, the functor of IndCoh $*$-pushforward along $\alpha$ is conservative when restricted to those objects whose singular support is contained in the image of the singular codifferential
    \[
    \Sing(\alpha) : \Sing(Z \times Z)_Z \mon \Sing(Z \times \Omega V).
    \]
We observe that this locus contains the image of (\ref{eqn:0_delta}).

\qed[Theorem \ref{thm:morita_D_H}]

\subsection{Transfer bimodules}\label{s:4-transfer-bimodules}

\subsubsection{}
\label{subsubsection:H_transfer_dual}
Let $f : Z \to Y$ be a morphism of quasi-smooth schemes. We may form the category
    \[
    \HH_{Z \to Y} = \ICoh_0((Z \times Y)^\wedge_Z)
    \]
associated to the graph of $f$ by the procedure of \S\ref{subsubsection:ICoh0}, regarded as an $(\HH(Z), \HH(Y))$-bimodule under convolution. According to\footnote{Strictly speaking, the material in this subsection is logically independent of this citation, and the results of this section can be used to give a new proof of \cite[Proposition 5.2.8]{beraldo2019sheaves} for quasi-smooth schemes. Namely, one can recover the ambidextrous duality between $\HH_{Z \to Y}$ and $\HH_{Y \leftarrow Z}$ from \cref{construction:morita_transfer} and the ambidextrous duality between (\ref{eqn:corresp_SSh_bimod}) and (\ref{eqn:corresp_SSh_bimod_reversed}).} \cite[Proposition 5.2.8]{beraldo2019sheaves}, the bimodule $\HH_{Z \to Y}$ admits a right dual (and also a left dual) canonically identified with
    \[
    \HH_{Y \leftarrow Z} = \ICoh_0((Y \times Z)^\wedge_Z),
    \]
viewed as an $(\HH(Y), \HH(Z))$-bimodule.

\medskip

In particular, we obtain a pair of adjoint functors
    \[
    \HH_{Z \to Y} : 
    \HH(Y) \MOD 
    \rightleftharpoons
    \HH(Z) \MOD
    : \HH_{Y \leftarrow Z}.
    \]
The purpose of this subsection is to describe the adjunction
    \[
    \DMod(\Sing(Y))^\Rightarrow
    \rightleftharpoons
    \DMod(\Sing(Z))^\Rightarrow
    \]
corresponding to this one under the Morita equivalence of Theorem \ref{thm:morita_D_H}.

\subsubsection{} Consider the correspondence
\begin{equation}
\label{eqn:singular_codiff_correspondence}
    \begin{tikzcd}
        \Sing(Y)_Z \arrow[r, "\Sing(f)"] \arrow[d] & \Sing(Z) \\
        \Sing(Y)
    \end{tikzcd}
\end{equation}
induced by the singular codifferential of $f$. It defines a bimodule
\begin{equation}
\label{eqn:corresp_SSh_bimod}
    \DMod(\Sing(Y)_Z)^\Rightarrow
    \in
    (\DMod(\Sing(Z))^\Rightarrow, \DMod(\Sing(Y))^\Rightarrow) \BIMOD.
\end{equation}
We claim:

\begin{construction}
\label{construction:morita_transfer}
The following diagram commutes:
\begin{equation}
\label{eqn:morita_pull}
\begin{tikzcd}[column sep=4em]
    \HH(Y) \MOD \arrow[d, "\HH_{Z \to Y}" left] \arrow[r, "\ICoh(Y)"] & \DMod(\Sing(Y))^\Rightarrow \MOD \arrow[d, "\DMod(\Sing(Y)_Z)^\Rightarrow"] \\
    \HH(Z) \MOD \arrow[r, "\ICoh(Z)" below] & \DMod(\Sing(Z))^\Rightarrow \MOD
\end{tikzcd}
\end{equation}
\end{construction}

\subsubsection{} Note that we may also view the correspondence (\ref{eqn:singular_codiff_correspondence}) as defining a bimodule
\begin{equation}
\label{eqn:corresp_SSh_bimod_reversed}
    \DMod(\Sing(Y)_Z)^\Rightarrow
    \in
    (\DMod(\Sing(Y))^\Rightarrow, \DMod(\Sing(Z))^\Rightarrow) \BIMOD.
\end{equation}
This bimodule is canonically right dual (and also left dual) to (\ref{eqn:corresp_SSh_bimod}), so after taking right adjoints of the vertical arrows, we obtain from (\ref{eqn:morita_pull}) the commutative diagram:
\begin{equation}
\label{eqn:morita_push}
\begin{tikzcd}[column sep=4em]
    \HH(Z) \MOD \arrow[d, "\HH_{Y \leftarrow Z}" left] \arrow[r, "\ICoh(Z)"] & \DMod(\Sing(Z))^\Rightarrow \MOD \arrow[d, "\DMod(\Sing(Y)_Z)^\Rightarrow"] \\
    \HH(Y) \MOD \arrow[r, "\ICoh(Y)" below] & \DMod(\Sing(Y))^\Rightarrow \MOD
\end{tikzcd}
\end{equation}

\subsubsection{Execution of Construction \ref{construction:morita_transfer}} We need to construct an isomorphism
    \[
    \ICoh(Z) 
    \underset{\HH(Z)}{\otimes} 
    \HH_{Z \to Y}
    \simeq
    \DMod(\Sing(Y)_Z)^\Rightarrow
    \underset{\DMod(\Sing(Y))^\Rightarrow}{\otimes}
    \ICoh(Y)
    \]
of $(\HH(Z), \DMod(\Sing(Y))^\Rightarrow)$-bimodules.

\medskip

In view of \S\ref{subsubsection:morita_inverse_tau}, it is equivalent to produce an isomorphism
\begin{equation}
\label{eqn:transfer_push_bimod_iso}
    \HH_{Z \to Y}
    \simeq
    \ICoh(Z)^\tau 
    \underset{\DMod(\Sing(Z))^\Rightarrow}{\otimes}
    \DMod(\Sing(Y)_Z)^\Rightarrow
    \underset{\DMod(\Sing(Y))^\Rightarrow}{\otimes}
    \ICoh(Y)
\end{equation}
of $(\HH(Z), \HH(Y))$-bimodules.

\medskip

This amounts to an isomorphism of $\HH(Z \times Y)$-modules
\begin{equation}
\label{eqn:ICoh_0_prod}
    \ICoh_0((Z \times Y)^\wedge_Z)
    \simeq
    \ICoh(Z \times Y) 
    \underset{\DMod(\Sing(Z \times Y))^\Rightarrow}{\otimes}
    \DMod(\Sing(Y)_Z)^\Rightarrow
\end{equation}
where $\DMod(\Sing(Y)_Z)^\Rightarrow$ is regarded as a module over $\DMod(\Sing(Z \times Y))^\Rightarrow$ via the map
\begin{equation}
\label{eqn:singular_codifferential_prod_sign}
    \Sing(Y)_Z
    \overset{\text{(\ref{eqn:singular_codiff_correspondence})}}{\longrightarrow}
    \Sing(Z) \times \Sing(Y)
    \overset{\tau \times \id}{\longrightarrow}
    \Sing(Z) \times \Sing(Y).
\end{equation}

\medskip

Both sides of (\ref{eqn:ICoh_0_prod}) are naturally full subcategories of
    \[
    \ICoh(Z \times Y)
    \]
stable under the action of $\HH(Z \times Y)$, and we claim that they are equal as such. Indeed, (\ref{eqn:singular_codifferential_prod_sign}) identifies with the closed embedding
    \[
    \Sing(Z/Z \times Y)
    \mon
    \Sing(Z \times Y)
    \]
so the claim follows from \cref{cor:icoh-0-tensor-formula}. 

\qed[\cref{construction:morita_transfer}]

\subsection{Composition of transfer bimodules}\label{s:4-comp-transf-bimod}

\subsubsection{} Let us record for future use the sense in which the formation of the diagrams (\ref{eqn:morita_push}) is compatible with composition. 

\medskip

Suppose that
    \[
    Z \to Y \to X
    \]
are morphisms of quasi-smooth schemes. According to\footnote{Strictly speaking, the material in this subsection is logically independent of this citation. Namely, the diagram (\ref{eqn:transfer_composite_push}) with its left face removed implies that there exists \emph{some} isomorphism of the form (\ref{eqn:transfer_composite_H}), which \cref{construction:transfer_bimod_composite} identifies with the convolution map.} \cite[Theorem 4.3.4]{beraldo2019sheaves}, the convolution functor
\begin{equation}
\label{eqn:transfer_composite_H}
    \HH_{Z \leftarrow Y}
    \underset{\HH(Y)}{\otimes}
    \HH_{Y \leftarrow X}
    \to
    \HH_{Z \leftarrow X}
\end{equation}
is an isomorphism of $(\HH(Z), \HH(X))$-bimodules. 

\medskip

On the other hand, the cartesian diagram
    \[
    \begin{tikzcd}
        \Sing(X)_Z \arrow[r] \arrow[d] & \Sing(Y)_Z \arrow[r] \arrow[d] & \Sing(Z) \\
        \Sing(X)_Y \arrow[r] \arrow[d] & \Sing(Y) & \\
        \Sing(X)
    \end{tikzcd}
    \]
defines an isomorphism
\begin{equation}
\label{eqn:transfer_composite_SSh}
    \DMod(\Sing(Y)_Z)^\Rightarrow
    \underset{\DMod(\Sing(Y))^\Rightarrow}{\otimes}
    \DMod(\Sing(X)_Y)^\Rightarrow
    \simeq
    \DMod(\Sing(X)_Y)^\Rightarrow
\end{equation}
of $(\DMod(\Sing(Z))^\Rightarrow, \DMod(\Sing(X))^\Rightarrow)$-bimodules.

\medskip

The situation at this point is that we have obtained commutativity data for every face of the diagram below.
\begin{equation}
\label{eqn:transfer_composite_push}
\begin{tikzcd}[column sep={5em, between origins}]
    \HH(Z) \MOD \arrow[dd] \arrow[rr] \arrow[dddr] & & \DMod(\Sing(Z))^\Rightarrow \MOD \arrow[dd] \arrow[dddr] & \\
    & &   & \\
    \HH(Y) \MOD \arrow[dr] \arrow[rr] & & \DMod(\Sing(Y))^\Rightarrow \MOD \arrow[dr] & \\
    & \HH(X) \MOD \arrow[rr] & & \DMod(\Sing(X))^\Rightarrow \MOD
\end{tikzcd}
\end{equation}
Namely, the left and right faces are given by the isomorphisms (\ref{eqn:transfer_composite_H}) and (\ref{eqn:transfer_composite_SSh}), respectively, while the remaining faces are given by (\ref{eqn:morita_push}).

\medskip

The goal of this subsection is to construct a 3-morphism representing the interior of (\ref{eqn:transfer_composite_push}):

\begin{construction}
\label{construction:transfer_bimod_composite}
The diagram (\ref{eqn:transfer_composite_push}) commutes.
\end{construction}

\subsubsection{Execution of \ref{construction:transfer_bimod_composite}} We need to show that the isomorphisms (\ref{eqn:transfer_push_bimod_iso}) used to perform \cref{construction:morita_transfer} identify (\ref{eqn:transfer_composite_H}) with the composite
\begin{equation}
\label{eqn:transfer_composite_SSh_2}
\begin{split}
    & \ICoh(X)^\tau 
    \underset{\DMod(\Sing(X))^\Rightarrow}{\otimes} \DMod(\Sing(X)_Y)^\Rightarrow 
    \underset{\DMod(\Sing(Y))^\Rightarrow}{\otimes}
    \ICoh(Y) \\
    & \qquad \qquad \underset{\HH(Y)}{\otimes}
    \ICoh(Y)^\tau
    \underset{\DMod(\Sing(Y))^\Rightarrow}{\otimes}
    \DMod(\Sing(Y)_Z)^\Rightarrow
    \underset{\DMod(\Sing(Z))^\Rightarrow}{\otimes}
    \ICoh(Z)  \\
    & \quad \simto
    \ICoh(X)^\tau 
    \underset{\DMod(\Sing(X))^\Rightarrow}{\otimes} \DMod(\Sing(X)_Y)^\Rightarrow 
    \underset{\DMod(\Sing(Y))^\Rightarrow}{\otimes} 
    \DMod(\Sing(Y))^\Rightarrow \\
    & \hspace{8em} \underset{\DMod(\Sing(Y))^\Rightarrow}{\otimes}
    \DMod(\Sing(Y)_Z)^\Rightarrow
    \underset{\DMod(\Sing(Z))^\Rightarrow}{\otimes}
    \ICoh(Z)  \\
    & \qquad \simto
    \ICoh(X)^\tau
    \underset{\DMod(\Sing(X))^\Rightarrow}{\otimes} \DMod(\Sing(X)_Y)^\Rightarrow \\
    & \hspace{8em} \underset{\DMod(\Sing(Y))^\Rightarrow}{\otimes}
    \DMod(\Sing(Y)_Z)^\Rightarrow
    \underset{\DMod(\Sing(Z))^\Rightarrow}{\otimes}
    \ICoh(Z)  \\
    & \qquad \quad \overset{\text{(\ref{eqn:transfer_composite_SSh})}}{\longrightarrow} 
    \ICoh(X)^\tau 
    \underset{\DMod(\Sing(X))^\Rightarrow}{\otimes} \DMod(\Sing(X)_Z)^\Rightarrow \underset{\DMod(\Sing(Z))^\Rightarrow}{\otimes}
    \ICoh(Z)
\end{split}
\end{equation}

\subsubsection{} Recall that the identification of \cref{construction:morita_transfer} matches the targets of (\ref{eqn:transfer_composite_H}) and (\ref{eqn:transfer_composite_SSh_2}) as full subcategories of
\begin{equation}
\label{eqn:ambient_target}
    \ICoh(X) \otimes \ICoh(Z)
\end{equation}
The sources of (\ref{eqn:transfer_composite_H}) and (\ref{eqn:transfer_composite_SSh_2}) are likewise identified as full subcategories of
\begin{equation}
\label{eqn:ambient_source}
    \ICoh(X) 
    \otimes
    \ICoh(Y) \underset{\HH(Y)}{\otimes} \ICoh(Y)
    \otimes
    \ICoh(Z).
\end{equation}

\medskip

Furthermore, under these embeddings, both (\ref{eqn:transfer_composite_H}) and (\ref{eqn:transfer_composite_SSh_2}) are induced by functors
    \[
    \text{(\ref{eqn:ambient_source})}
    \to
    \text{(\ref{eqn:ambient_target})}
    \]
between the ambient categories. Namely, the ambient functor for (\ref{eqn:transfer_composite_H}) is the one obtained from the pairing
\begin{equation}
\label{eqn:ICoh_pairing_H}
    \ICoh(Y) \underset{\HH(Y)}{\otimes} \ICoh(Y)
    \to
    \Vect
\end{equation}
and the ambient functor for (\ref{eqn:transfer_composite_SSh_2}) is the one obtained from the composite
\begin{equation}
\label{eqn:ICoh_pairing_Gamma}
    \ICoh(Y) \underset{\HH(Y)}{\otimes} \ICoh(Y)^\tau
    \simto
    \DMod(\Sing(Y))^\Rightarrow
    \overset{(\Gamma_\dR)^\Rightarrow}{\longrightarrow}
    \Vect.
\end{equation}
So it remains to see that the two functors (\ref{eqn:ICoh_pairing_H}) and (\ref{eqn:ICoh_pairing_Gamma}) are canonically isomorphic.

\subsubsection{} 
\label{subsubsection:ev_comparison}
Let us use the notation of (\ref{eqn:A_B_M}) for convenience. Then we are trying to identify the composites
    \[
    \MM^\vee \otimes_\BB \MM
    \overset{\ev}{\longrightarrow}
    \AA 
    \overset{(\Gamma_\dR)^\Rightarrow}{\to} 
    \Vect
    \qquad \text{and} \qquad
    \MM^\vee \otimes_\BB \MM \to \Vect.
    \] 
where the second map is the pairing on $\MM$ as a left-dualizable $\BB$-module.
\medskip

According to \S\ref{subsubsection:D_mod_duality}, the left composite is given by
    \[
    \MM^\vee \otimes_\BB \MM
    \overset{\id \otimes \co\act}{\longrightarrow}
    \MM^\vee \otimes_\BB \MM \otimes \AA
    \to
    \AA
    \overset{(\Gamma_\dR)^\Rightarrow}{\longrightarrow}
    \Vect,
    \]
which rewrites as
    \[
    \MM^\vee \otimes_\BB \MM
    \overset{\id \otimes \co\act}{\longrightarrow}
    \MM^\vee \otimes_\BB \MM \otimes \AA
    \overset{\id \otimes (\Gamma_\dR)^\Rightarrow}{\longrightarrow}
    \MM^\vee \otimes_\BB \MM
    \to
    \Vect.
    \]
Therefore, it is enough to show that the composite
\begin{equation}
\label{eqn:coact_Gamma}
    \MM 
    \overset{\co\act}{\longrightarrow}
    \MM \otimes \AA
    \overset{(\Gamma_\dR)^\Rightarrow}{\longrightarrow}
    \MM
\end{equation}
is the identity.

\medskip

For this, it suffices to consider the universal case $\MM = \AA$, where the assertion is clear because (\ref{eqn:coact_Gamma}) becomes the composite
    \[
    \AA
    \overset{(\Delta_*)^\Rightarrow}{\longrightarrow}
    \AA \otimes \AA
    \overset{(\pr_{2, *})^\Rightarrow}{\longrightarrow}
    \AA.
    \]

\qed[\cref{construction:transfer_bimod_composite}]

\section{Extension to stacks}

In this section, we will establish an version of the Morita equivalence of Theorem \ref{thm:morita_D_H} when the quasi-smooth scheme $Z$ is replaced with a quasi-smooth stack $\YY$.

\medskip

Let us point out that although we do have a natural bimodule structure
    \[
    \ICoh(\YY)
    \in
    (\DMod(\Sing(\YY))^\Rightarrow, \HH(\YY))
    \BIMOD,
    \]
the resulting functor
\begin{equation}
\label{eqn:H_to_SSh_stack}
    \HH(\YY) \MOD
    \to
    \DMod(\Sing(\YY))^\Rightarrow \MOD
\end{equation}
will usually not be an equivalence. This is essentially due to the failure of 1-affineness for the de Rham space of an algebraic stack. Because of this, the correct generalization of Theorem \ref{thm:morita_D_H} to the setting of stacks will involve sheaves of $(\infty, 2)$-categories.

\medskip

We will then return to the category $\DMod(\Sing(\YY))^\Rightarrow \MOD$ and its relation to $\HH(\YY)$ and show that (\ref{eqn:H_to_SSh_stack}) is a colocalization in the $(\infty, 2)$-category $\twoDGCAT$ of stable $(\infty, 2)$-categories tensored over $\DGCAT$.

\medskip 

The outline of this section is as follows:
\begin{itemize}
    \item In \S\ref{s:5-sheaf-morita} we describe a category $\ShCat_{\SSh}(\YY)$, associated to any quasi-smooth stack $\YY$, which is related to $D(\Sing(\YY))^\Rightarrow \MOD$ in the same way that crystals of categories are related to $D(\YY) \MOD$.
    \item In \S\ref{subsubsection:ShCat_morita_construction} we prove that $\ShCat_{\SSh}(\YY)$ is equivalent to the category $\ShCat_\HH(\YY)$ defined in \cite{beraldo2019sheaves}. This amounts to showing that our Morita equivalence \Cref{thm:morita_D_H} is sufficiently functorial as the scheme $Z$ varies.
    \item In \S\ref{s:5-rel-d-mod-stacks} we show that the result of \S\ref{s:4-transfer-bimodules} on the image of transfer bimodules under our Morita equivalence extends to stacks, conditional on a construction involving relative $D$-modules \Cref{construction:rel_D_stack}.
    \item In \S\ref{subsubsection:rel_D_stack_construction} we carry out \Cref{construction:rel_D_stack}. 
    \item In \S\ref{s:5-sheared-crystals} we apply the results of the previous subsections to study the relationship between $D(\Sing(\YY))^\Rightarrow \MOD$ and $\ShCat_{\SSh}(\YY)$. 
\end{itemize}

\subsection{A sheaf version of the Morita equivalence}\label{s:5-sheaf-morita}

\subsubsection{} Let $\YY$ be a quasi-smooth algebraic stack, and let $\Sm_{/\YY}$ denote the category whose objects are schemes $S$ equipped with a smooth morphism $S \to \YY$ and whose morphisms are composites $T \to S \to \YY$ with $T \to S$ smooth.

\subsubsection{} Recall that \cite{beraldo2019sheaves} assembles the assignment
\begin{equation}
\label{eqn:H_MOD_obj_mor}
    S \rightsquigarrow \HH(S)\MOD
    \qquad
    (T \to S) \rightsquigarrow \HH_{T \to S}  \otimes_{\HH(S)} (-)
\end{equation}
into a functor
    \[
    \HH\MOD : 
    (\Sm_{/\YY})^\op \to \twoDGCAT.
    \]

\medskip

Now let\footnote{Smooth descent for $\HH$ \cite[Theorem 6.1.2]{beraldo2019sheaves} implies that this limit agrees with the notion of sheaves of categories with $\HH$-coefficients considered in \emph{loc. cit.}}
    \[
    \ShCat_\HH(\YY) = \lim_{S \to \YY} \HH(S) \MOD
    \]
denote the limit of the functor $\HH\MOD$. The $\HH$-affineness result of \emph{loc. cit.} exhibits an adjoint equivalence
    \[
    \LOC_\HH : 
    \HH(\YY)\MOD
    \rightleftharpoons
    \ShCat_\HH(\YY)
    : \GAMMA_\HH.
    \]
In particular, viewing $\ICoh(\YY)$ a module over $\HH(\YY)$, we obtain a canonical object
\begin{equation}
\label{eqn:LOC_H_ICoh}
    \LOC_\HH(\ICoh(\YY))
    \in
    \ShCat_\HH(\YY).
\end{equation}

\subsubsection{} 
\label{subsubsection:Loc_ICoh_obj_mor}
According to \cite[Lemma 6.6.1]{beraldo2019sheaves}, the value of (\ref{eqn:LOC_H_ICoh}) at an object $S \in \Sm_{/\YY}$ is
    \[
    \LOC_\HH(\ICoh(\YY))_S 
    \simeq
    \ICoh(S)
    \]
and its value at a morphism $T \to S$ is given by the convolution map
    \[
    \HH_{T \to S} \underset{\HH(S)}{\otimes} \ICoh(S)
    \simto
    \ICoh(T).
    \]

\subsubsection{} Now let us describe the counterparts of $\HH\MOD$ and $\Sh\Cat_\HH(\YY)$ that will appear on the other side of the Morita equivalence.

\medskip

Let $f : T \to S$ be any morphism in $\Sm_{/\YY}$. Since $f$ is required to be smooth, its singular codifferential is an isomorphism
    \[
    \Sing(S)_T \simto \Sing(T)
    \]
and the correspondence (\ref{eqn:singular_codiff_correspondence}) induces a map
\begin{equation}
\label{eqn:Sing_T_S}
    \Sing(T) \to \Sing(S).
\end{equation}
In particular, we may consider the adjunction
    \[
    \IND : 
    \DMod(\Sing(S))^\Rightarrow 
    \rightleftharpoons
    \DMod(\Sing(T))^\Rightarrow
    : \OBLV.
    \]

\medskip

The left adjoints assemble into a functor
    \[
    \mathbf{S}^\Rightarrow\MOD : (\Sm_{/\YY})^\op \to \twoDGCAT.
    \]
given at the level of objects and morphisms by the assignment
    \[
    (T \to S) 
    \rightsquigarrow 
    \Big(
    \DMod(\Sing(S))^\Rightarrow \MOD
    \overset{\IND}{\too}
    \DMod(\Sing(T))^\Rightarrow \MOD
    \Big).
    \]

\medskip

Let
\begin{equation}
\label{eqn:ShCat_SSh_limit}
    \ShCat_\SSh(\YY) = \lim_{S \to \YY} \DMod(\Sing(S))^\Rightarrow \MOD
\end{equation}
denote the limit of $\SSh\MOD$, viewed as a functor. It is equipped with tautological evaluation functors
    \[
    \GAMMA(S; -) : 
    \ShCat_\SSh(\YY) 
    \to 
    \DMod(\Sing(S))^\Rightarrow \MOD.
    \]

\medskip

For any stack $\mathscr{W}$ equipped with a $\Gm$-equivariant morphism
    \[
    \mathscr{W} \to \Sing(\YY),
    \]
we obtain an object
    \[
    \ul\DMod(\WW)^\Rightarrow
    \in
    \ShCat_\SSh(\YY)
    \]
with
    \[
    \GAMMA(S; \ul\DMod(\WW))
    \simeq
    \DMod(\WW_S)^\Rightarrow
    \qquad
    \WW_S = S \times_\YY \WW.
    \]

\begin{construction}
\label{construction:ShCat_morita}
There is an equivalence of categories
\begin{equation}
\label{eqn:morita_ShCat}
    \ShCat_\HH(\YY)
    \simeq
    \ShCat_\SSh(\YY)
\end{equation}
exchanging the objects
    \[
    \LOC_\HH(\ICoh(\YY))
    \qquad \text{and} \qquad
    \ul\DMod(\Sing(\YY))^\Rightarrow.
    \]
\end{construction}

\subsubsection{} We will obtain \cref{construction:ShCat_morita} from a slightly stronger result, \cref{construction:functor_morita}, whose statement we will now explain.

\medskip

Let $\DGCAT$ denote the constant functor
\begin{equation}
\label{eqn:Sm_Y_op_2dgcat}
    (\Sm_{/\YY})^\op
    \to
    \twoDGCAT
\end{equation}
with value $\DGCAT$. Then we may regard $\LOC_H(\ICoh(\YY))$ as a natural transformation
\begin{equation}
\label{eqn:ICoh_dgcat_H}
    \DGCAT \to \HH \MOD
\end{equation}
of functors of the shape (\ref{eqn:Sm_Y_op_2dgcat}), whose value at an object $S \in \Sm_{/\YY}$ is
\begin{equation}
\label{eqn:ICoh_dgcat_H_S}
    \ICoh(S) : \DGCAT \to \HH(S) \MOD.
\end{equation}
Similarly, we may regard $\ul\DMod(\Sing(\YY))^\Rightarrow$ as a natural transformation
\begin{equation}
\label{eqn:ind_dgcat_SSh}
    \DGCAT \to \SSh \MOD,
\end{equation}
whose value at an object $S \in \Sm_{/\YY}$ is
\begin{equation}
\label{eqn:ind_dgcat_SSh_S}
    \IND : \DGCAT \to \DMod(\Sing(S))^\Rightarrow \MOD.
\end{equation}

\medskip

To perform \cref{construction:ShCat_morita}, it suffices to perform the following and pass to the limit:

\begin{construction}
\label{construction:functor_morita}
There is a commutative diagram
    \[
    \begin{tikzcd}[column sep={3em, between origins}]
        \HH \MOD \arrow[rr, "\sim"] & & \SSh \MOD \\
        & \DGCAT \arrow[ul, "\text{(\ref{eqn:ICoh_dgcat_H})}"] \arrow[ur, "\text{(\ref{eqn:ind_dgcat_SSh})}" below right] &
    \end{tikzcd}
    \]
of functors of the shape (\ref{eqn:Sm_Y_op_2dgcat}). 
\end{construction}

\subsubsection{} 

\medskip

The value of the putative natural isomorphism
    \[
    \HH \MOD \simeq \SSh \MOD
    \]
at an individual morphism $T \to S$ is the diagram
\begin{equation}
\label{eqn:morita_pull_sm}
    \begin{tikzcd}[column sep=4em]
        \HH(S) \MOD \arrow[d, "\HH_{T \to S}" left] \arrow[r, "\ICoh(S)"] & \DMod(\Sing(S))^\Rightarrow \MOD \arrow[d, "\IND"] \\
        \HH(T) \MOD \arrow[r, "\ICoh(T)"] & \DMod(\Sing(T))^\Rightarrow \MOD
    \end{tikzcd}
\end{equation}
provided by Construction \ref{construction:morita_transfer}.

\medskip

Therefore, the content of Construction \ref{construction:functor_morita} is that the diagrams (\ref{eqn:morita_pull_sm}) can be assembled into an natural transformation. This will be done in \S\ref{subsubsection:ShCat_morita_construction}.

\subsection{Execution of Construction \ref{construction:functor_morita}}
\label{subsubsection:ShCat_morita_construction}

\subsubsection{} It turns out that \cref{construction:functor_morita} as stated is difficult to approach directly, so we will have provide a more tractable reformulation. This will appear below as \cref{construction:functor_R_morita}.

\subsubsection{} Recall from \S\ref{subsubsection:H_transfer_dual} that all of the bimodules $\HH_{T \to S}$ admit right duals. Therefore, by passing to right adjoints in $\twoDGCAT$, we obtain from $\HH\MOD$ a functor
    \[
    \HH \MOD^R : \Sm_{/\YY} \to \twoDGCAT
    \]
given at the level of objects and morphisms by the assignment
    \[
    S \rightsquigarrow \HH(S)\MOD
    \qquad
    (T \to S) \rightsquigarrow \HH_{S \leftarrow T} \underset{\HH(T)}{\otimes} (-).
    \]

\medskip

Applying the same procedure to $\SSh\MOD$, we obtain a functor
    \[
    \SSh \MOD^R : \Sm_{/\YY} \to \twoDGCAT
    \]
given at the level of objects and morphisms by the assignment    
    \[
    (T \to S) 
    \rightsquigarrow
    \Big(
    \DMod(\Sing(T))^\Rightarrow \MOD
    \overset{\OBLV}{\too} 
    \DMod(\Sing(S))^\Rightarrow \MOD
    \Big).
    \]

\medskip

Now let $\DGCAT^R$ denote the constant functor
\begin{equation}
\label{eqn:Sm_Y_2dgcat}
    \Sm_{/\YY} \to \twoDGCAT
\end{equation}
with value $\DGCAT$ (alternatively, it is obtained from the functor $\DGCAT$ by passing to right adjoints). Since $\ICoh(S)$ is self-dual as a module over $\HH(S)$, all of the functors (\ref{eqn:ICoh_dgcat_H_S}) admit right adjoints in $\twoDGCAT$. Therefore, by passing to right adjoints, we obtain from (\ref{eqn:ICoh_dgcat_H}) a natural transformation
\begin{equation}
\label{eqn:ICoh_dgcat_H_R}
    \HH\MOD^R
    \to
    \DGCAT^R
\end{equation}
of functors of the shape (\ref{eqn:Sm_Y_2dgcat}). Similarly, by virtue of the $(\IND, \OBLV)$ adjunction, we obtain from (\ref{eqn:ind_dgcat_SSh}) a natural transformation
\begin{equation}
\label{eqn:ind_dgcat_SSh_R}
    \SSh\MOD^R
    \to
    \DGCAT^R
\end{equation}
whose value at $S \in \Sm_{/\YY}$ is the forgetful functor
\begin{equation}
\label{eqn:oblv_S}
    \OBLV : \DMod(\Sing(S))^\Rightarrow \MOD \to \DGCAT.
\end{equation}

\medskip

Construction \ref{construction:functor_morita} is equivalent to:

\begin{construction}
\label{construction:functor_R_morita}
There is a commutative diagram
    \[
    \begin{tikzcd}[column sep={3.5em, between origins}]
        \HH \MOD^R \arrow[dr, "\text{(\ref{eqn:ICoh_dgcat_H_R})}" below left] \arrow[rr, "\sim"] & & \SSh \MOD^R \arrow[dl, "\text{(\ref{eqn:ind_dgcat_SSh_R})}"] \\
        & \DGCAT^R &
    \end{tikzcd}
    \]
of functors of the shape (\ref{eqn:Sm_Y_2dgcat}). 
\end{construction}

\subsubsection{} The feature of \cref{construction:functor_R_morita} which makes it manageable is the fact that the forgetful functors (\ref{eqn:oblv_S}) are locally fully faithful. This means (see \cref{aproposition:lff_lifting}) that in order to perform \cref{construction:functor_R_morita} it is equivalent to do the following:
\begin{itemize}
    \item[(i)] for every object $S \in \Sm_{/\YY}$, specify a commutative diagram
        \[
        \begin{tikzcd}[column sep={4.5em, between origins}]
            \HH(S) \MOD \arrow[rr, "\sim"] \arrow[dr, "\ICoh(S)" below left] & & \DMod(\Sing(S))^\Rightarrow \MOD \arrow[dl, "\OBLV"] \\
            & \DGCAT &
        \end{tikzcd}
        \]
    \item[(ii)] for every morphism $T \to S$ in $\Sm_{/\YY}$, verify that the resulting diagram:
        \begin{equation}
        \label{eqn:oblv_prism}
        \begin{tikzcd}
            \HH(T) \MOD \arrow[rr, "\sim"] \arrow[dd, "\HH_{S \leftarrow T}"] \arrow[rrrd, "\ICoh(T)" below left] & & \DMod(\Sing(T))^\Rightarrow \MOD \arrow[dd, "\OBLV"] \arrow[dr, "\OBLV"] & \\
            &  & & \DGCAT \arrow[dd] \\
            \HH(S) \MOD \arrow[rr, "\sim"] \arrow[rrrd, "\ICoh(S)" below left] & & \DMod(\Sing(S))^\Rightarrow \MOD \arrow[dr, "\OBLV"] & \\
            &  & & \DGCAT
        \end{tikzcd}
        \end{equation}
    (in which the back face is empty) can be completed to a commutative prism.
\end{itemize}
Note that the data required by (i) is provided by \cref{thm:morita_D_H}, so it remains to verify the assertion of (ii). This will occupy the remainder of the subsection.

\subsubsection{} The front face of (\ref{eqn:oblv_prism}) encodes the isomorphism of right $\HH(T)$-modules
\begin{equation}
\label{eqn:ICoh_transfer_convolution}
    \ICoh(S) \underset{\HH(S)}{\otimes} \HH_{S \leftarrow T}
    \to
    \ICoh(T)
\end{equation}
given by convolution. We need to check that (\ref{eqn:ICoh_transfer_convolution}) is linear for the pre-existing left actions of $\DMod(\Sing(S))^\Rightarrow$ on its source and target. 

\medskip

Indeed, let us write (\ref{eqn:ICoh_transfer_convolution}) as the composite
\begin{equation}
\label{eqn:ICoh_transfer_convolution_2}
    \begin{split}
        & \ICoh(S) \underset{\HH(S)}{\otimes} \HH_{S \leftarrow T}
        \simto
        \ICoh(S) 
        \underset{\HH(S)}{\otimes} 
        \ICoh(S)^\tau
        \underset{\DMod(\Sing(S))^\Rightarrow}{\otimes}
        \ICoh(T) \\
        & \qquad \mon
        \ICoh(S) 
        \underset{\HH(S)}{\otimes} 
        \ICoh(S)^\tau
        \otimes \ICoh(T) 
        \overset{\ev \otimes \id}{\too}
        \ICoh(T).
    \end{split}
\end{equation}
By \S\ref{subsubsection:ev_comparison}, we may further rewrite (\ref{eqn:ICoh_transfer_convolution_2}) as
\begin{equation*}
    \begin{split}
        & \ICoh(S) \underset{\HH(S)}{\otimes} \HH_{S \leftarrow T} 
        \simto
        \ICoh(S) 
        \underset{\HH(S)}{\otimes} 
        \ICoh(S)^\tau
        \underset{\DMod(\Sing(S))^\Rightarrow}{\otimes}
        \ICoh(T) \\
        & \qquad \simto
        \DMod(\Sing(S))^\Rightarrow
        \underset{\DMod(\Sing(S))^\Rightarrow}{\otimes}
        \ICoh(T) 
        \mon
        \DMod(\Sing(S))^\Rightarrow \otimes \ICoh(T) \\
        & \qquad \overset{(\Gamma_\dR)^\Rightarrow \otimes \id}{\too}
        \ICoh(T),
    \end{split}
\end{equation*}
This last expression identifies with
\begin{equation*}
    \begin{split}
        & \ICoh(S) \underset{\HH(S)}{\otimes} \HH_{S \leftarrow T} 
        \simto
        \ICoh(S) 
        \underset{\HH(S)}{\otimes} 
        \ICoh(S)^\tau
        \underset{\DMod(\Sing(S))^\Rightarrow}{\otimes}
        \ICoh(T) \\
        & \qquad \simto
        \DMod(\Sing(S))^\Rightarrow
        \underset{\DMod(\Sing(S))^\Rightarrow}{\otimes}
        \ICoh(T) 
        \simto
        \ICoh(T).
    \end{split}
\end{equation*}

\qed[\Cref{construction:functor_R_morita}]

\subsection{Relative D-modules over stacks}\label{s:5-rel-d-mod-stacks}

\subsubsection{} Let $\ZZ \to \YY$ be a schematic morphism of quasi-smooth stacks. Then $\ICoh(\YY^\wedge_\ZZ)$ is naturally acted on by $\HH(\YY)$, and we may form the object
    \[
    \LOC(\ICoh(\YY^\wedge_\ZZ)) \in \ShCat_\HH(\YY).
    \]
On the other hand, from the morphism
    \[
    \Sing(\YY)_\ZZ \to \Sing(\YY),
    \]
we obtain an object
    \[
    \ul\DMod(\Sing(\YY)_\ZZ)^\Rightarrow \in \ShCat_\SSh(\YY).
    \]
We note for future reference that for $S \in \Sm_{/\YY}$, we have
\begin{equation}
\label{eqn:sing_Sm_Y}
    \Sing(S)_{\ZZ_S} 
    \simeq 
    \left(\Sing(\YY)_\ZZ\right)_S.
\end{equation}

\begin{construction}
\label{construction:rel_D_stack}
The Morita equivalence (\ref{eqn:morita_ShCat}) exchanges the objects
    \[
    \ul\DMod(\Sing(\YY)_\ZZ)^\Rightarrow
    \qquad
    \text{and}
    \qquad
    \LOC(\ICoh(\YY^\wedge_\ZZ)).
    \]
\end{construction}

\subsubsection{} 
\label{subsubsection:rel_D_stack_explain}
Let us explain why \cref{construction:rel_D_stack} is nontrivial. Consider the object
    \[
    \DD \in \ShCat_\SSh(\YY)
    \]
corresponding to $\LOC(\ICoh(\YY^\wedge_\ZZ))$ under the Morita equivalence (\ref{eqn:morita_ShCat}). 

\medskip

The value of $\LOC(\ICoh(\YY^\wedge_\ZZ))$ at an object $S \in \Sm_{/\YY}$ identifies with
    \[
    \LOC(\ICoh(\YY^\wedge_\ZZ))_S 
    \simeq 
    \HH_{S \to \YY} \underset{\HH(\YY)}{\otimes} \ICoh(\YY^\wedge_\ZZ)
    \simeq
    \ICoh(S^\wedge_{\ZZ_S})
    \]
with its natural action of $\HH(S)$. Therefore, the equivalence (\ref{eqn:rel_DMod_tensor_Sing}), combined with the Morita equivalence for schemes, produces an equivalence of categories
\begin{equation}
\label{eqn:rel_D_stack_S}
    \DMod(\Sing(S)_{\ZZ_S})^\Rightarrow
    \simeq
    \DD_S
\end{equation}
linear over $\DMod(\Sing(S))^\Rightarrow$. However, for the purposes of \cref{construction:rel_D_stack}, this observation by itself is insufficient---it is not obvious how to supply the functoriality required to promote the data of (\ref{eqn:rel_D_stack_S}) to an equivalence in $\ShCat_\SSh(\YY)$.

\medskip

So the challenge is to produce a morphism
\begin{equation}
\label{eqn:rel_D_stack_map}
    \ul\DMod(\Sing(\YY)_\ZZ)^\Rightarrow
    \to
    \DD
\end{equation}
in $\ShCat_\SSh(\YY)$. Once this has been done, it will be straightforward to verify that the value of (\ref{eqn:rel_D_stack_map}) at every object of $\Sm_{/\YY}$ agrees with the equivalence (\ref{eqn:rel_D_stack_S}), and hence that (\ref{eqn:rel_D_stack_map}) itself is an equivalence.

\medskip

We shall delay the execution of \cref{construction:rel_D_stack} to \S\ref{subsubsection:rel_D_stack_construction}, and use the remainder of this subsection to explain some of its consequences.

\subsubsection{} The colocalization
    \[
    \ICoh_0(\YY^\wedge_\ZZ)
    \rightleftharpoons
    \ICoh(\YY^\wedge_\ZZ)
    \]
of $\HH(\YY)$-modules induces a colocalization
\begin{equation}
    \LOC(\ICoh_0(\YY^\wedge_\ZZ))
    \rightleftharpoons
    \LOC(\ICoh(\YY^\wedge_\ZZ)).
\end{equation}

\medskip

On the other hand, observe that we may identify
    \[
    \Sing(\ZZ_S/S) 
    \simeq 
    S \underset{\YY}{\times} \Sing(\ZZ/\YY)
    \]
as closed subschemes of (\ref{eqn:sing_Sm_Y}). In particular, we obtain a colocalization
    \[
    \ul\DMod(\Sing(\ZZ/\YY))^\Rightarrow
    \rightleftharpoons
    \ul\DMod(\Sing(\YY)_\ZZ)^\Rightarrow
    \]
in $\ShCat_\SSh(\YY)$ whose value an object $S \in \Sm_{/\YY}$ is given by the colocalization
    \[
    \DMod(\Sing(\ZZ_S/S))^\Rightarrow
    \rightleftharpoons
    \DMod(\Sing(S)_{\ZZ_S})^\Rightarrow.
    \]

\begin{proposition}
\label{prop:rel_D_stack_0}
\cref{construction:rel_D_stack} exchanges the colocalizations
    \[
    \ul\DMod(\Sing(\ZZ/\YY))^\Rightarrow
    \qquad \text{and} \qquad
    \LOC(\ICoh_0(\YY^\wedge_\ZZ)).
    \]
\end{proposition}

\subsubsection{} 
\label{subsubsection:rel_H_stack}
Under our running assumptions on stacks (see \S\ref{s:stack-assumptions}) for every schematic morphism $f : \ZZ \to \YY$ of quasi-smooth stacks the graph $\Gamma_f: \ZZ \to \ZZ \times \YY$ is also schematic. Recall that one may consider the $(\HH(\ZZ), \HH(\YY))$-bimodule
    \[
    \HH_{\ZZ \to \YY} = \ICoh_0((\ZZ \times \YY)^\wedge_\ZZ).
    \]

\medskip

Recall also the correspondence of stacks
\begin{equation}
\label{eqn:singular_codifferential_stack}
    \begin{tikzcd}
        \Sing(\YY)_\ZZ \arrow[r] \arrow[d] & \Sing(\ZZ) \\
        \Sing(\YY)
    \end{tikzcd}
\end{equation}
induced by the singular codifferential of $f$. Let us form the object
    \[
    \ul\DMod(\Sing(\YY)_\ZZ)^\Rightarrow
    \in
    \ShCat_\SSh(\ZZ \times \YY)
    \]
obtained from the morphism
\begin{equation}
\label{eqn:sing_codifferential_prod_stack}
    \Sing(\YY)_\ZZ 
    \overset{\text{(\ref{eqn:singular_codifferential_stack})}}{\too} 
    \Sing(\ZZ) \times \Sing(\YY)
    \overset{\tau \times \id}{\too}
    \Sing(\ZZ) \times \Sing(\YY).
\end{equation}

\medskip

Since (\ref{eqn:sing_codifferential_prod_stack}) identifies with the map
\begin{equation}
\label{eqn:rel_sing_stack}
    \Sing(\ZZ/\ZZ \times \YY) \to \Sing(\ZZ \times \YY)_\ZZ
\end{equation}
associated to the graph of $f$
    \[
    \ZZ \to \ZZ \times \YY,
    \]
we obtain the following consequence of \cref{prop:rel_D_stack_0}:

\begin{construction}
\label{construction:rel_D_stack_prod}
The Morita equivalence (\ref{eqn:morita_ShCat}) for $\ZZ \times \YY$ exchanges the objects
    \[
    \ul\DMod(\Sing(\YY)_\ZZ)^\Rightarrow
    \qquad
    \text{and}
    \qquad
    \LOC(\HH_{\ZZ \to \YY})
    \]
\end{construction}

\subsection{Execution of Construction \ref{construction:rel_D_stack}}
\label{subsubsection:rel_D_stack_construction}

\subsubsection{} As indicated in \S\ref{subsubsection:rel_D_stack_explain}, the task at hand is more or less to construct the morphism (\ref{eqn:rel_D_stack_map}). This will be accomplished by  subjecting the functor of IndCoh $!$-pullback
\begin{equation}
\label{eqn:ICoh_pull_H}
    \ICoh(\YY) \to \ICoh(\YY^\wedge_\ZZ),
\end{equation}
viewed as a morphism of $\HH(\YY)$-modules, to a series of transformations.

\medskip

First we apply $\LOC$ to (\ref{eqn:ICoh_pull_H}), and obtain a morphism
\begin{equation}
\label{eqn:ICoh_pull_Loc}
    \LOC(\ICoh(\YY))
    \to
    \LOC(\ICoh(\YY^\wedge_\ZZ))
\end{equation}
in $\ShCat_\HH(\YY)$.

\medskip

Then (\ref{eqn:ICoh_pull_Loc}) may be viewed as a morphism
\begin{equation}
\label{eqn:ICoh_pull_morita}
    \ul\DMod(\Sing(\YY))^\Rightarrow
    \to
    \DD
\end{equation}
in $\ShCat_\SSh(\YY)$.

\subsubsection{Digression} To proceed further, we will need to introduce the category
    \[
    \ShCat_\SSh(\YY)_\ZZ
    =
    \lim_{S \to \YY} \DMod(\Sing(S)_{\ZZ_S})^\Rightarrow \MOD,
    \]
where the limit is taken over $(\Sm_{/\YY})^\op$ with respect to the transition functors of induction:
\begin{equation*}
    \IND : 
    \DMod(\Sing(S)_{\ZZ_S})^\Rightarrow \MOD
    \to
    \DMod(\Sing(T)_{\ZZ_T})^\Rightarrow \MOD
\end{equation*}

\medskip

Consider the adjunctions
    \[
    \IND_S :
    \DMod(\Sing(S))^\Rightarrow \MOD
    \rightleftharpoons
    \DMod(\Sing(S)_{\ZZ_S})^\Rightarrow \MOD
    : \OBLV_S.
    \]
Observe that the left adjoints $\IND_S$ assemble tautologically into a functor
\begin{equation}
\label{eqn:SSh_Z_ind}
    \IND :
    \ShCat_\SSh(\YY)
    \to
    \ShCat_\SSh(\YY)_\ZZ
\end{equation}
making the following diagrams commute:
\[
\begin{tikzcd}
    \ShCat_\SSh(\YY) \arrow[r, "\IND"] \arrow[d] & \ShCat_\SSh(\YY)_\ZZ \arrow[d] \\
    \DMod(\Sing(S))^\Rightarrow \MOD
    \arrow[r, "\IND_S"] & \DMod(\Sing(S)_{\ZZ_S})^\Rightarrow \MOD
\end{tikzcd}
\]
Since the map $\ZZ \to \YY$ is schematic, the Beck--Chevalley transformations
\begin{equation}
\label{eqn:SSh_Z_BC}
    \begin{tikzcd}
        \DMod(\Sing(S))^\Rightarrow \MOD \arrow[d, "\IND_{T \to S}" left] & \DMod(\Sing(S)_{\ZZ_S})^\Rightarrow \MOD \arrow[l, "\OBLV_S" above] \arrow[d, "\IND_{T \to S}"] \\
        \DMod(\Sing(T))^\Rightarrow \MOD &  \DMod(\Sing(T)_{\ZZ_T})^\Rightarrow \MOD \arrow[l, "\OBLV"] 
    \end{tikzcd}
\end{equation}
are isomorphisms. Therefore, we obtain by passing to the limit that (\ref{eqn:SSh_Z_ind}) admits a right adjoint
    \[
    \OBLV : 
    \ShCat_\SSh(\YY)_\ZZ
    \to
    \ShCat_\SSh(\YY).
    \]
making the following diagrams commute:
\[
\begin{tikzcd}
    \ShCat_\SSh(\YY) \arrow[d] & \ShCat_\SSh(\YY)_\ZZ \arrow[d] \arrow[l, "\OBLV" above] \\
    \DMod(\Sing(S))^\Rightarrow \MOD
    & \DMod(\Sing(S)_{\ZZ_S})^\Rightarrow \MOD \arrow[l, "\OBLV_S"] 
\end{tikzcd}
\]

\medskip

Note that we may identify
\begin{equation}
\label{eqn:oblv_ind_DSing}
    \OBLV \circ \IND(\ul\DMod(\Sing(\YY))^\Rightarrow)
    \simeq
    \ul\DMod(\Sing(\YY)_\ZZ)^\Rightarrow.
\end{equation}

\subsubsection{} 
\label{subsubsection:S_description}
Note that the value of (\ref{eqn:ICoh_pull_Loc}) at an object $S \in \Sm_{/\YY}$ is given by the functor of IndCoh $!$-pullback
\begin{equation}
\label{eqn:ICoh_pull_S}
    \ICoh(S)
    \to
    \ICoh(S^\wedge_{\ZZ_S}),
\end{equation}
regarded as a morphism of $\HH(S)$-modules. 

\medskip

We may view (\ref{eqn:ICoh_pull_S}) as the functor
    \[
    \ICoh(S)
    \to
    \ICoh(S) 
    \underset{\DMod(\Sing(S))^\Rightarrow}{\otimes}
    \DMod(\Sing(S)_{\ZZ_S})^\Rightarrow,
    \]
induced by the D-module $!$-pullback
\begin{equation}
\label{eqn:DMod_pull_S}
    \DMod(\Sing(S))^\Rightarrow
    \to
    \DMod(\Sing(S)_{\ZZ_S})^\Rightarrow.
\end{equation}
In other words, the value (\ref{eqn:ICoh_pull_morita}) at $S$ identifies with the morphism (\ref{eqn:DMod_pull_S}).

\medskip

In particular, this lifts $\DD_S$ to an object
    \[
    \tilde{\DD}_S 
    \in 
    \DMod(\Sing(S)_{\ZZ_S})^\Rightarrow \MOD.
    \]
We claim:

\begin{proposition}
\label{prop:D_S_lift}
The data of $\tilde{\DD}_S$ for $S \in \Sm_{/\YY}$ determine a lift of $\DD$ to an object
    \[
    \tilde{\DD}
    \in
    \ShCat_\SSh(\YY)_\ZZ.
    \]
\end{proposition}

\subsubsection{} Note that this assertion is a proposition and not a construction as \cref{prop:rel_crys_locally_ff} (which is applicable by the running assumption that $\ZZ \to \YY$ is schematic) implies that the functors $\OBLV_S$ are locally fully faithful. Namely, it amounts to the assertion that for every arrow $T \to S$ in $\Sm_{/\YY}$, the map of $\DMod(\Sing(T))^\Rightarrow$-modules
\begin{equation}
\label{eqn:D_tilde_map}
\begin{split}
    \OBLV_T \circ \IND_{T \to S} (\tilde{\DD}_S)
    \simeq
    \IND_{T \to S} \circ \OBLV_S (\tilde{\DD}_S)
    & \simeq
    \IND_{T \to S} (\DD_S) \\
    & \qquad \simto
    \DD_T
    \simeq
    \OBLV_T(\tilde{\DD}_T)
\end{split}
\end{equation}
is actually linear over $\DMod(\Sing(T)_{\ZZ_T})^\Rightarrow$.

\subsubsection{} We now take \cref{prop:D_S_lift} as given and resume the construction of (\ref{eqn:rel_D_stack_map}). Its proof will appear at the end of this subsection, \S\ref{subsubsection:proof_prop_D_S_lift}.

\medskip

By \cref{prop:D_S_lift}, we may perceive (\ref{eqn:ICoh_pull_morita}) as a morphism
\begin{equation}
\label{eqn:ICoh_pull_oblv}
    \ul\DMod(\Sing(\YY))^\Rightarrow
    \to
    \OBLV(\tilde{\DD})
\end{equation}
in $\ShCat_\SSh(\YY)$. 

\medskip

The transpose of (\ref{eqn:ICoh_pull_oblv}) under the $(\IND, \OBLV)$ adjunction is a morphism
\begin{equation}
\label{eqn:ICoh_pull_ind}
    \IND(\ul\DMod(\Sing(\YY))^\Rightarrow)
    \to
    \tilde{\DD}
\end{equation}
in $\ShCat_\SSh(\YY)_\ZZ$.

\medskip

After applying $\OBLV$ to (\ref{eqn:ICoh_pull_ind}), we obtain a morphism
\begin{equation}
\label{eqn:ICoh_pull_ind_oblv}
    \OBLV \circ \IND(\ul\DMod(\Sing(\YY))^\Rightarrow)
    \to
    \DD
\end{equation}
in $\ShCat_\SSh(\YY)$. 

\medskip

Finally, after making the identification (\ref{eqn:oblv_ind_DSing}), we find that (\ref{eqn:ICoh_pull_ind_oblv}) defines the sought-after morphism
\begin{equation}
\label{eqn:ICoh_pull_final}
    \ul\DMod(\Sing(\YY)_\ZZ)^\Rightarrow
    \to
    \DD
\end{equation}
in $\ShCat_\SSh(\YY)$. 

\subsubsection{} To finish \cref{construction:rel_D_stack}, it remains to verify that (\ref{eqn:ICoh_pull_final}) is an isomorphism. Indeed, by the description of (\ref{eqn:ICoh_pull_morita}) given in \S\ref{subsubsection:S_description}, the value of (\ref{eqn:ICoh_pull_final}) at an object $S \in \Sm_{/\YY}$ is the canonical arrow
    \[
    \DMod(\ZZ_S) 
    \underset{\DMod(S)}{\otimes}
    \DMod(\Sing(S))^\Rightarrow
    \simto
    \DMod(\Sing(S)_{\ZZ_S})^\Rightarrow.
    \]

\qed[\cref{construction:rel_D_stack}]

\subsubsection{Proof of Proposition \ref{prop:D_S_lift}} 
\label{subsubsection:proof_prop_D_S_lift}
Consider the diagram
\[
\begin{tikzcd}[column sep={5.5em, between origins}]
    & \HH(\ZZ) \MOD \arrow[dd] \arrow[rr] & & \HH(\ZZ_S) \MOD \arrow[dd] \arrow[rr] \arrow[dl] & & \HH(\ZZ_T) \MOD \arrow[dd] \arrow[dl] \\
    & & \DMod(\Sing(S)_{\ZZ_S})^\Rightarrow \MOD \arrow[dd] \arrow[rr] & & \DMod(\Sing(T)_{\ZZ_T})^\Rightarrow \MOD \arrow[dd] & \\
    & \HH(\YY) \MOD \arrow[dl] \arrow[rr] & & \HH(S) \MOD \arrow[rr] \arrow[dl] & & \HH(T) \MOD \arrow[dl] \\
    \ShCat_\SSh(\YY) \arrow[rr] & & \DMod(\Sing(S))^\Rightarrow \MOD \arrow[rr] & & \DMod(\Sing(T))^\Rightarrow \MOD
\end{tikzcd}
\]
where the upper back-to-front maps are the composites
\[\HH(\ZZ_S) \MOD \totext{\sim} D(\Sing(\ZZ_S))^\Rightarrow \MOD \totext{\IND} D(\Sing(S)_{\ZZ_S})^\Rightarrow \MOD.\]
The back two faces are commutative by the functoriality of $\HH$, the front face is commutative by the Beck--Chevalley transformation (\ref{eqn:SSh_Z_BC}), the lower two faces are commutative by the construction of the Morita equivalence (\ref{eqn:morita_ShCat}), and the two vertical faces are commutative by \cref{construction:morita_transfer}.

\medskip

At the moment we have \emph{not} filled in the upper face of the cube. Observe that the isomorphism (\ref{eqn:D_tilde_map}) is encoded by the commutativity of the preceding diagram, when evaluated on the object
    \[
    \ICoh(\ZZ) \in \HH(\ZZ) \MOD,
    \]
and that a completion of the right half of the diagram into a commutative cube would lift (\ref{eqn:D_tilde_map}) to an isomorphism in $\DMod(\Sing(T)_{\ZZ_T})^\Rightarrow$. Let us construct this completion.

\medskip

All five faces of the diagram
\[
\begin{tikzcd}[column sep={5.5em, between origins}]
    & \HH(\ZZ_S) \MOD \arrow[dd] \arrow[dl] & & \HH(\ZZ_T) \MOD \arrow[dd] \arrow[dl] \arrow[ll] \\
    \DMod(\Sing(S)_{\ZZ_S})^\Rightarrow \MOD \arrow[dd] & & \DMod(\Sing(T)_{\ZZ_T})^\Rightarrow \MOD \arrow[dd] \arrow[ll] & \\
    & \HH(S) \MOD \arrow[dl] & & \HH(T) \MOD \arrow[dl]  \arrow[ll] \\
    \DMod(\Sing(S))^\Rightarrow \MOD & & \DMod(\Sing(T))^\Rightarrow \MOD \arrow[ll]
\end{tikzcd}
\]
obtained by taking right adjoints of the horizontal arrows remains commutative, so it is equivalent to complete this new diagram into a commutative cube. To do so, apply \cref{construction:transfer_bimod_composite} to the composites
    \[
    \ZZ_T \to \ZZ_S \to S
    \qquad
    \ZZ_T \to T \to S
    \]
and paste together the two resulting prisms. 

\qed[\cref{prop:D_S_lift}]

\subsubsection{Porism} 
\label{subsubsection:rel_D_stack_pullbacks}
\cref{construction:rel_D_stack} comes with an additional compatibility: the Morita equivalence (\ref{eqn:morita_ShCat}) exchanges the morphism
    \[
    \ul\DMod(\Sing(\YY))^\Rightarrow 
    \to
    \ul\DMod(\Sing(\YY)_\ZZ)^\Rightarrow
    \]
induced by D-module $!$-pullback with the morphism
    \[
    \LOC \Big( 
    \ICoh(\YY)
    \to \ICoh(\YY^\wedge_\ZZ)
    \Big)
    \]
induced by IndCoh $!$-pullback.

\medskip

Combining this with \cref{prop:rel_D_stack_0}, we see that \cref{construction:rel_D_stack_prod} comes with an additional compatibility: the Morita equivalence (\ref{eqn:morita_ShCat}) for $\ZZ \times \YY$ exchanges the morphism
    \[
    \ul\DMod(\Sing(\ZZ \times \YY))^\Rightarrow
    \to
    \ul\DMod(\Sing(\ZZ/\ZZ \times \YY))^\Rightarrow
    \]
induced by D-module $!$-pullback with the morphism
    \[
    \LOC \Big(
    \ICoh(\ZZ \times \YY)
    \to
    \ICoh((\ZZ \times \YY)^\wedge_\ZZ)
    \to
    \HH_{\ZZ \to \YY}
    \Big)
    \]
obtained by composing the functor of IndCoh $!$-pullback with the right adjoint of the colocalization defining $\HH_{\ZZ \to \YY}$.

\subsection{Sheared crystals of categories on the stack of singularities}\label{s:5-sheared-crystals}

\subsubsection{} The goal of this somewhat digressive subsection is to study the relationship between $\DMod(\Sing(\YY))^\Rightarrow \MOD$ and $\ShCat_\SSh(\YY)$. We will see that they are connected by a pair of adjoint functors
    \[
    \LOC : 
    \DMod(\Sing(\YY))^\Rightarrow \MOD
    \rightleftharpoons
    \ShCat_\SSh(\YY) :
    \GAMMA
    \]
in $\twoDGCAT$, and that this adjunction is in fact a colocalization. Applying the Morita equivalence (\ref{eqn:morita_ShCat}) and the H-affineness equivalence, we obtain that $\ICoh(\YY)$ is naturally a $(\HH(\YY), \DMod(\Sing(\YY))^\Rightarrow)$-bimodule defining a colocalization
    \[
    \DMod(\Sing(\YY))^\Rightarrow
    \mon
    \HH(\YY) \MOD.
    \]

\subsubsection{} Observe that the family of functors
    \[
    \IND :
    \DMod(\Sing(\YY))^\Rightarrow \MOD
    \to 
    \DMod(\Sing(S))^\Rightarrow \MOD.
    \]
assemble into a functor
    \[
    \LOC : \DMod(\Sing(\YY))^\Rightarrow\MOD \to \ShCat_{\mathbf{S}^\Rightarrow}(\YY).
    \]
We claim:

\begin{proposition}
\label{prop:GAMMA_2DGCAT}
The functor $\LOC$ admits a right adjoint in $\twoDGCAT$.
\end{proposition}

\subsubsection{} It is easy to see that $\LOC$ admits a right adjoint
    \[
    \GAMMA :
    \ShCat_\SSh(\YY) 
    \to 
    \DMod(\Sing(\YY))^\Rightarrow \MOD,
    \]
computed by the formula
    \[
    \CC
    \mapsto 
    \lim_{S \to \YY} \GAMMA(S; \CC),
    \]
where the limit is taken over $\Sm_{/\YY}$. The content of Proposition \ref{prop:GAMMA_2DGCAT} is that this functor is continuous and commutes with the operation of tensoring with a plain dg category.

\subsubsection{Proof of Proposition \ref{prop:GAMMA_2DGCAT}} In light of the Morita equivalence (\ref{eqn:morita_ShCat}) and the H-affineness theorem, it is equivalent to show that the functor
    \[
    \DMod(\Sing(\YY))^\Rightarrow \MOD
    \overset{\LOC}{\too} 
    \ShCat_\SSh(\YY)
    \simto
    \ShCat_\HH(\YY)
    \simto 
    \HH(\YY)\MOD
    \]
admits a right adjoint in $\twoDGCAT$. 

\medskip

Now observe that we have the adjunction
    \[
    \IND : 
    \DGCAT 
    \rightleftharpoons 
    \DMod(\Sing(\YY))^\Rightarrow\MOD
    : \OBLV
    \]
in $\twoDGCAT$, where the right adjoint $\OBLV$ is conservative. Therefore, it is enough to show that the composite
    \[
    \DGCAT 
    \overset{\IND}{\too}
    \DMod(\Sing(\YY))^\Rightarrow \MOD
    \to \HH(\YY)\MOD
    \]
admits a continuous right adjoint.

\medskip

This functor is the one defined by $\ICoh(\YY)$, as an $\HH(\YY)$-module. Therefore, it is equivalent to show that $\ICoh(\YY)$ is dualizable as an $\HH(\YY)$-module. This follows from the fact that $\ICoh(\YY)$ is dualizable as a plain dg category and that $\HH(\YY)$ is rigid.

\qed[Proposition \ref{prop:GAMMA_2DGCAT}]

\begin{proposition}
\label{prop:LOC_ff}
The functor $\LOC$ is fully faithful.
\end{proposition}

\subsubsection{Proof} Since $(\LOC, \GAMMA)$ is an adjunction in $\twoDGCAT$, it suffices to check that the functor
    \[
    \DMod(\Sing(\YY))^\Rightarrow
    \to
    \End_{\ShCat_\SSh(\YY)}
    \big(
    \LOC(\DMod(\Sing(\YY))^\Rightarrow)
    \big)
    \]
is an equivalence. But this functor rewrites as
    \[
    \DMod(\Sing(\YY))^\Rightarrow
    \to
    \lim_{S \to \YY} \DMod(\Sing(S))^\Rightarrow,
    \]
where the limit is taken over $\Sm_{/\YY}$, so the result follows from smooth descent for D-modules.

\qed[\cref{prop:LOC_ff}]

\subsubsection{Remark} Note that we have a commutative diagram
\[
\begin{tikzcd}
    \DMod(\YY)\MOD \arrow[d, "\RES"] \arrow[r, "\LOC"] & \ShCat(\YY_\dR) \arrow[d] \\
    \DMod(\Sing(\YY))^\Rightarrow \MOD \arrow[r, "\LOC"] & \ShCat_\SSh(\YY)
\end{tikzcd}
\]
in which the vertical arrows are fully faithful. Therefore, \cref{prop:LOC_ff} confirms \cite[Conjecture 2.6.6]{gaitsgory2015sheaves} for stacks which are nil-isomorphic to a quasi-smooth stack.

\section{Trace computations}

The outline of this section is as follows:
\begin{itemize}
    \item In \S\ref{s:6-cat-trace-review} we review the basic formalism of categorical traces and prove some simple results about the behavior of traces with respect to duality. The results of this section are a categorification of the fact that the trace of a matrix equals the trace of its transpose. 
    \item In \S\ref{s:6-marcel-duchamp} we use \Cref{thm:morita_D_H} to calculate the trace of $\HH(\YY) \MOD$ with respect to the bimodule $\HH_{\phi:\YY \to \YY}$ for any schematic endomorphism of a quasi-smooth stack. 
    \item In \S\ref{s:6-rel-enh-trace} we relate the Chern class of $\ICoh(\YY)$ as a ``Weil-object'' of $\HH(\YY) \MOD$ with the enhanced trace of \cite{gaitsgory2022toy}.
    \item In \S\ref{s:6-sing-supp} we calculate how taking traces interacts with the singular support filtration on $\ICoh(\YY)$. 
\end{itemize}

\subsection{Review of categorical trace}\label{s:6-cat-trace-review}

\subsubsection{} Recall that for any endomorphism $\phi : \CC \to \CC$ of a dualizable object in an ambient symmetric monoidal category $\mathbf{O}$, it makes sense to form the trace
    \[
    \tr(\phi; \CC) \in \End(\id_\mathbf{O}),
    \]
given definitionally by the composite
    \[
    \id_\bfO 
    \overset{\coev}{\longrightarrow}
    \CC \otimes \CC^\vee
    \overset{\phi \otimes \id}{\longrightarrow}
    \CC \otimes \CC^\vee
    \overset{\ev}{\longrightarrow}
    \id_\bfO.
    \]
Under the identification $\id_\bfO \simeq \id_\bfO^\vee$ the following diagram always naturally commutes
\[
\begin{tikzcd}
\id_\bfO \arrow[r, "\psi"] \arrow[d, "\sim"] & \id_\bfO \arrow[d, "\sim"] \\
\id_\bfO^\vee \arrow[r, "\psi^\vee"] & \id_\bfO^\vee
\end{tikzcd}
\]
for any $\psi$ in $\End(\id_\bfO)$. This gives rise to a canonical identification
\begin{equation}
\label{eqn:tr_of_dual}
    \tr(\phi; \CC) \simeq \tr(\phi^\vee; \CC^\vee).
\end{equation}

\subsubsection{} Suppose now that $\mathbf{O}$ is a symmetric monoidal $(\infty, 2)$-category. Then a natural transformation $\alpha : \phi \to \psi$ of endomorphisms induces a morphism
    \[
    \tr(\phi; \CC) \to \tr(\psi; \CC).
    \]
The formation of this arrow is compatible with duality in the sense that the morphism
    \[
    \tr(\phi^\vee; \CC^\vee)
    \to
    \tr(\psi^\vee; \CC^\vee)
    \]
induced by the natural transformation $\alpha^\vee : \phi^\vee \to \psi^\vee$ fits into a commutative diagram:
\begin{equation}
\label{eqn:tr_dual}
    \begin{tikzcd}
        \tr(\phi; \CC) \arrow[r] \arrow[d, "\sim"] & \tr(\psi; \CC) \arrow[d, "\sim"] \\
        \tr(\phi^\vee; \CC^\vee) \arrow[r] & \tr(\psi^\vee; \CC^\vee)
    \end{tikzcd}
\end{equation}

\subsubsection{} More generally, suppose that $\DD$ is another dualizable object equipped with an endomorphism, which we shall denote also by $\phi$. For any right adjointable functor $F : \CC \to \DD$ equipped with a natural transformation
\begin{equation}
\label{eqn:endo_lax}
    \begin{tikzcd}
        \CC \arrow[d, "F" left] \arrow[r, "\phi"] & \CC \arrow[d, "F"] \\
        \DD \arrow[r, "\phi" below] & \DD
        \arrow[Rightarrow, from=1-2, to=2-1]
    \end{tikzcd}
    \qquad \qquad
    \alpha : F \circ \phi \to \phi \circ F
\end{equation}
one obtains a morphism
\begin{equation}
\label{eqn:tr_endo_lax}
    \cl(F, \alpha) : \tr(\phi; \CC) \to \tr(\phi; \DD),
\end{equation}
given definitionally by the composite
\begin{equation*}
\begin{split}
    \tr(\phi; \CC)
    \overset{\mathrm{unit}}{\longrightarrow}
    & \tr(\phi \circ F^R \circ F; \CC) \\
    & \qquad \simeq
    \tr(F \circ \phi \circ F^R; \DD)
    \overset{\alpha}{\longrightarrow}
    \tr(\phi \circ F \circ F^R; \DD)
    \overset{\mathrm{counit}}{\longrightarrow}
    \tr(\phi; \DD),
\end{split}
\end{equation*}
where middle equivalence the cyclicity of trace (see for instance \cite[Definition 3.24]{benzvi2019nonlineartraces} \cite[\S 0.2.4]{gaitsgory2022toy}).

\subsubsection{} 
\label{subsubsection:tr_funct_dual}
The $2$-categorical functoriality of traces enjoys a further compatibility with duality. Namely, consider the natural transformation
    \[
    \alpha^\vee : 
    \phi^\vee \circ F^\vee
    \to
    F^\vee \circ \phi^\vee
    \]
obtained from $\alpha$ by passing to dual functors. Since $F$ is right adjointable, $F^\vee$ is left adjointable with $(F^\vee)^L \simeq (F^R)^\vee$. Therefore, we obtain from $\alpha^\vee$ a natural transformation
    \[
    \beta : 
    (F^\vee)^L \circ \phi^\vee
    \to
    \phi^\vee \circ (F^\vee)^L.
    \]
We will obtain the following compatibility by a routine diagram chase:

\begin{computation}
\label{computation:adj_dual_tr}
There is a commutative diagram:
    \[
    \begin{tikzcd}
        \tr(\phi; \CC) \arrow[r, "{\cl(F, \alpha)}"] \arrow[d, "\sim"] & \tr(\phi; \DD) \arrow[d, "\sim"] \\
        \tr(\phi^\vee; \CC^\vee) \arrow[r, "{\cl((F^\vee)^L, \beta)}"] & \tr(\phi^\vee; \DD^\vee)
    \end{tikzcd}
    \]
\end{computation}

\subsubsection{} Applying (\ref{eqn:tr_dual}) to every arrow appearing in the definition of $\cl((F^\vee)^L, \beta)$ identifies the lower circuit of this diagram with the composite
\begin{equation}
\label{eqn:adj_dual_tr_1}
\begin{split}
    & \tr(\phi; \CC) 
    \overset{\unit}{\longrightarrow}
    \tr(F^R \circ F \circ \phi; \CC)
    \simeq
    \tr(F \circ \phi \circ F^R; \DD) \\
    & \qquad \qquad \overset{\beta^\vee}{\longrightarrow}
    \tr(F \circ F^R \circ \phi; \DD)
    \overset{\counit}{\longrightarrow}
    \tr(\phi; \DD)
    \simeq
    \tr(\phi^\vee; \DD^\vee),
\end{split}
\end{equation}
where the second arrow is induced by the natural transformation
    \[
    \beta^\vee : \phi \circ F^R \to F^R \circ \phi. 
    \]

\medskip

Using the cyclic property of the trace, we rewrite the first line of (\ref{eqn:adj_dual_tr_1}) to replace it with
\begin{equation}
\label{eqn:adj_dual_tr_2}
\begin{split}
    & \tr(\phi; \CC) 
    \overset{\unit}{\longrightarrow}
    \tr(\phi \circ F^R \circ F; \CC)
    \simeq
    \tr(F \circ \phi \circ F^R; \DD) \\
    & \qquad \qquad \overset{\beta^\vee}{\longrightarrow}
    \tr(F \circ F^R \circ \phi; \DD)
    \overset{\counit}{\longrightarrow}
    \tr(\phi; \DD)
    \simeq
    \tr(\phi^\vee; \DD^\vee).
\end{split}
\end{equation}
Comparing this expression with the definition of $\cl(\alpha; \phi)$, we see that it will suffice to produce a commutative diagram:
\begin{equation}
\label{eqn:adj_dual_tr_dgrm}
\begin{tikzcd}
    F \circ \phi \circ F^R \arrow[r, "\beta^\vee"] \arrow[d, "\alpha"] & F \circ F^R \circ \phi \arrow[d, "\counit"] \\
    \phi \circ F \circ F^R \arrow[r, "\counit"] & \phi
\end{tikzcd}
\end{equation}

\medskip

By definition, the upper circuit of (\ref{eqn:adj_dual_tr_dgrm}) is given by the upper circuit of:
\begin{equation}
\label{eqn:adj_dual_tr_dgrm_2}
\begin{tikzcd}
    F \circ \phi \circ F^R \arrow[d, "\unit"] & & \\
    F \circ F^R \circ F \circ \phi \circ F^R \arrow[r, "\alpha"] \arrow[d, "\counit"] & F \circ F^R \circ \phi \circ F \circ F^R \arrow[r, "\counit"] \arrow[d, "\counit"] & F \circ F^R \circ \phi \arrow[d, "\counit"] \\
    F \circ \phi \circ F^R \arrow[r, "\alpha"] & \phi \circ F \circ F^R \arrow[r, "\counit"] & \phi
\end{tikzcd}
\end{equation}
The triangular identity implies that the left vertical composite of (\ref{eqn:adj_dual_tr_dgrm_2}) is the identity, so its lower circuit agrees with the lower circuit of (\ref{eqn:adj_dual_tr_dgrm}). This completes Computation \ref{computation:adj_dual_tr}.

\qed[\cref{computation:adj_dual_tr}]

\subsubsection{} Let $\AA$ be a monoidal category and let $\Phi$ be an $(\AA, \AA)$-bimodule, perceived as an endofunctor
    \[
    \Phi : \AA \MOD \to \AA \MOD
    \qquad
    \MM \mapsto \Phi \otimes_\AA \MM.
    \]
It makes sense to take the trace of this endofunctor when $\Phi$ is right dualizable; this is the dg category
    \[
    \tr(\Phi; \AA \MOD) = \AA \underset{\AA \otimes \AA^\rev}{\otimes} \Phi.
    \]

\medskip

Now let $(\MM, \alpha)$ be an $\AA$-module equipped with lax equivariance data for $\Phi$; that is, $\MM$ is an $\AA$-module and $\alpha$ is a morphism of $\AA$-modules
    \[
    \alpha : \MM \to \Phi \otimes_\AA \MM.
    \]
If $\MM$ admits a right dual as an $\AA$-module, then one may form the Chern class
    \[
    \cl(\MM, \alpha) \in \tr(\Phi; \AA \MOD),
    \]
given by the composite
\begin{equation}
\label{eqn:chern_class}
\begin{split}
    \id
    \overset{\coev}{\longrightarrow} 
    \MM^\vee \otimes_\AA \MM
    & \overset{\id \otimes \alpha}{\longrightarrow}
    \MM^\vee \otimes_\AA (\Phi \otimes_\AA \MM) \\
    & \qquad \qquad \simeq
    (\MM \otimes \MM^\vee) \underset{\AA \otimes \AA^\rev}{\otimes} \Phi
    \overset{\ev \otimes \id}{\longrightarrow}
    \AA \underset{\AA \otimes \AA^\rev}{\otimes} \Phi.
\end{split}
\end{equation}
Let us note for future reference that this composite admits an alternative description. Since $\MM$ is dualizable, we may view $\alpha$ as a morphism of $(\AA, \AA)$-bimodules
    \[
    \alpha' : \MM \otimes \MM^\vee \to \Phi.
    \]
Then (\ref{eqn:chern_class}) identifies with the composite
\begin{equation}
\label{eqn:chern_class_rewrite}
\begin{split}
    \id
    & \overset{\coev}{\longrightarrow} 
    \MM^\vee \otimes_\AA \MM \\
    & \qquad \qquad \simeq 
    \AA \underset{\AA \otimes \AA^\rev}{\otimes} (\MM \otimes \MM^\vee)
    \to
    \AA \underset{\AA \otimes \AA^\rev}{\otimes} \Phi.
\end{split}
\end{equation}

\subsection{Carr\'e du champ}\label{s:6-marcel-duchamp}

\subsubsection{} Let $\phi : \YY \to \YY$ be a schematic endomorphism of a quasi-smooth stack. Consider the $(\HH(\YY), \HH(\YY))$-bimodule
    \[
    \HH_\phi = \HH_{\YY \to \YY}
    \]
associated to the morphism $\phi$ by \S\ref{subsubsection:rel_H_stack}. We will use the same symbol to denote the corresponding endomorphism:
    \[
    \HH_\phi : \HH(\YY)\MOD \to \HH(\YY)\MOD.
    \]

\medskip

Observe that $\ICoh(\YY)$ is naturally equipped with a lax equivariance datum
    \[
    \alpha :
    \ICoh(\YY) \to \HH_\phi \underset{\HH(\YY)}{\otimes} \ICoh(\YY)
    \]
for $\HH_\phi$, corresponding under the self-duality of $\ICoh(\YY)$ as an $\HH(\YY)$-module to the composite
\begin{equation}
\label{eqn:alpha_definition}
    \ICoh(\YY \times \YY) 
    \to
    \ICoh((\YY \times \YY)^\wedge_\YY)
    \to
    \HH_\phi
\end{equation}
of IndCoh $!$-pullback along the formal completion of graph of $\phi$ with the right adjoint of the colocalization defining $\HH_\phi$. The goal of this subsection is to give a description of the object
    \[
    \cl(\ICoh(\YY), \alpha)
    \in
    \tr(\HH_\phi; \HH(\YY)\MOD).
    \]

\subsubsection{} Recall that a stack $\XX$ is said to be quasi-quasi-smooth if its cotangent complex has tor amplitude $[2, -1]$. When this happens, $\Cot_\XX[-2]$ is a perfect complex of tor amplitude $[0, -3]$ and we may form the algebraic stack
    \[
    \Sing_2(\XX) 
    = \Spec_\XX \left( \Sym (\Cot_\XX[-2])^\vee \right).
    \]
We think of $\Sing_2(\XX)$ as measuring the extent to which $\XX$ fails to be quasi-smooth.

\subsubsection{} The reason we need to consider quasi-quasi-smooth stacks is that the stack of fixed points $\YY^\phi$ is not usually quasi-smooth. Instead, it is quasi-quasi-smooth, being a fiber product of quasi-smooth stacks:
\[
\begin{tikzcd}
    \YY^\phi \arrow[r] \arrow[d] & \YY \arrow[d, "\id \times \phi"] \\
    \YY \arrow[r, "\Delta"] & \YY \times \YY
\end{tikzcd}
\]
One has that
    \[
    \Sing_2(\YY^\phi)
    \simeq
    \Sing(\YY/\YY \times \YY) \underset{\Sing(\YY \times \YY)}{\times} \Sing_\phi(\YY/ \YY \times \YY),
    \]
where
    \[
    \Sing_\phi(\YY/ \YY \times \YY) \to \Sing(\YY) \times \Sing(\YY).
    \]
is the map (\ref{eqn:rel_sing_stack}) for the morphism $\phi$.

\begin{computation}
\label{computation:bimod_trace}
We have
\begin{equation}
\label{eqn:trace_computation}
    \tr(\HH_\phi(\YY); \HH(\YY) \MOD)
    \simeq
    \DMod \left( \Sing_2(\YY^\phi)
    \right)^\Rightarrow,
\end{equation}
and this identification takes the object
\begin{equation}
\label{eqn:ICoh_class}
    \cl(\ICoh(\YY), \alpha)
    \in
    \tr(\HH_\phi; \HH(\YY)\MOD)
\end{equation}
to the object
    \[
    (\omega_{\Sing_2(\YY^\phi)})^\Rightarrow 
    \in
    \DMod(\Sing_2(\YY^\phi))^\Rightarrow.
    \]
\end{computation}

\subsubsection{} Let us begin with the identification
\begin{equation}
\label{eqn:tr_H_phi}
\begin{split}
    \tr(\HH_\phi; \HH(\YY)\MOD)
    & \simeq \HH(\YY) \underset{\HH(\YY) \otimes \HH(\YY)^\rev}{\otimes} \HH_\phi \\
    & \qquad \qquad \simeq
    \Maps_{(\HH(\YY),\HH(\YY))}(\HH(\YY), \HH_\phi).
\end{split}
\end{equation}
Since the coevaluation
    \[
    \Vect 
    \to 
    \ICoh(\YY) \underset{\HH(\YY)}{\otimes} \ICoh(\YY)
    \simeq
    \Maps_{(\HH(\YY),\HH(\YY))}(\HH(\YY), \ICoh(\YY \times \YY))
    \]
classifies the map of $(\HH(\YY),\HH(\YY))$-bimodules
    \[
    \HH(\YY) 
    \mon 
    \ICoh((\YY \times \YY)^\wedge_\YY)
    \overset{(-)_*}{\longrightarrow}
    \ICoh(\YY \times \YY),
    \]
the description (\ref{eqn:chern_class_rewrite}) of Chern classes implies that the isomorphism (\ref{eqn:tr_H_phi}) takes the class of $\ICoh$ (\ref{eqn:ICoh_class}) to the object of
    \[
    \Maps_{(\HH(\YY),\HH(\YY))}(\HH(\YY), \HH_\phi)
    \]
given by the composite
\begin{equation}
\label{eqn:H_phi_composite}    
    \HH(\YY) 
    \to \ICoh(\YY \times \YY)
    \to \HH_\phi.
\end{equation}

\medskip

After applying the equivalence
    \[
    (\HH(\YY), \HH(\YY))\BIMOD
    \simeq
    \HH(\YY \times \YY) \MOD
    \simeq
    \ShCat_\HH(\YY \times \YY)
    \overset{\text{(\ref{eqn:morita_ShCat})}}{\longrightarrow}
    \ShCat_\SSh(\YY \times \YY),
    \]
we find that (\ref{eqn:tr_H_phi}) identifies with
\begin{equation}
\label{eqn:tr_H_SSh}
    \Maps_{\ShCat_\SSh(\YY \times \YY)}\left(
    \ul\DMod(\Sing(\YY/\YY \times \YY))^\Rightarrow, 
    \ul\DMod(\Sing_\phi(\YY/ \YY \times \YY))^\Rightarrow 
    \right).
\end{equation}
By \S\ref{subsubsection:rel_D_stack_pullbacks}, this identification sends (\ref{eqn:H_phi_composite}) to the composite
\begin{equation}
\label{eqn:SSh_phi_composite}
    \ul\DMod(\Sing(\YY/\YY \times \YY))^\Rightarrow
    \to
    \ul\DMod(\Sing(\YY \times \YY))^\Rightarrow
    \to
    \ul\DMod(\Sing_\phi(\YY/\YY \times \YY))^\Rightarrow
\end{equation}
in $\ShCat_\SSh(\YY \times \YY)$.

\medskip

We finish the computation by writing (\ref{eqn:tr_H_SSh}) as the limit of the categories
\begin{equation}\label{eq:Sing_2-mapping-space}
\begin{split}
    & \lim_{S \to \YY \times \YY} \Maps_{\DMod(\Sing(S))^\Rightarrow} \left(
    \DMod(\Sing(\YY/\YY \times \YY)_S)^\Rightarrow,
    \DMod(\Sing_\phi(\YY/\YY \times \YY)_S)^\Rightarrow
    \right) \\
    & \qquad \qquad 
    \simeq
    \lim_{S \to \YY \times \YY}     
    \DMod \left(
    \Sing_2(\YY^\phi)_S
    \right)^\Rightarrow
    \simeq
    \DMod(\Sing_2(\YY^\phi))^\Rightarrow
\end{split}
\end{equation}
taken over $\Sm_{/\YY \times \YY}$. By base-change for the Cartesian diagram 
\[
\begin{tikzcd}
\Sing_2(\YY^\phi) \arrow[r] \arrow[d] & \Sing(\YY / \YY \times \YY) \arrow[d] \\
\Sing_\phi(\YY / \YY \times \YY) \arrow[r] & \Sing(\YY \times \YY)
\end{tikzcd}
\]
the composite arrow (\ref{eqn:SSh_phi_composite}) is sent to $(\omega_{\Sing_2(\YY^\phi)})^\Rightarrow$ under (\ref{eq:Sing_2-mapping-space}).

\qed[\Cref{computation:bimod_trace}]

\subsubsection{Remark} In the case where $\phi$ is the identity morphism, this computation produces a equivalence of $\mathbf{E}_2$-categories:
    \[
    \ZZ_{\mathbf{E}_2}(\HH(\YY))
    =
    \End_{\HH(\YY) \BIMOD}(\HH(\YY))
    \simeq
    \DMod(\Sing_2(L\YY))^\Rightarrow.
    \]
On the other hand, the first author has exhibited (see \cite{beraldo2021center}) an equivalence of plain dg categories
    \[
    \ZZ_{\mathbf{E}_2}(\HH(\YY))
    \simeq
    \DMod^\der(L\YY).
    \]
The corresponding $\mathbf{E}_2$-monoidal structure on the dg category $\DMod^\der(L \YY)$ is expected to be the one induced by the ``pair of pants" construction for $L \YY$. To confirm this expectation, one needs the establish the functoriality of derived D-modules.

\subsection{Relation to the enhanced trace}\label{s:6-rel-enh-trace}

\subsubsection{} In this subsection we will explain the relation between the object
\begin{equation}
\label{eqn:enh_tr_SSh}
    \cl(\ICoh(\YY), \alpha) 
    \in 
    \tr(\HH_\phi; \HH(\YY) \MOD)
    \simeq
    \DMod(\Sing_2(\YY^\phi))^\Rightarrow
\end{equation}
and the so-called \emph{enhanced trace}
\begin{equation}
\label{eqn:enh_tr_QCoh}
    \cl(\ICoh(\YY), \phi^!)
    \in
    \tr(\RES_\phi; \QCoh(\YY) \MOD)
    \simeq
    \QCoh(\YY^\phi)
\end{equation}
of \cite{gaitsgory2022toy}.

\medskip

Let us recall the definition of the latter object. Consider the natural action of $\QCoh(\YY)$ on $\ICoh(\YY)$. Then the functor of IndCoh $!$-pullback
    \[
    \phi^! : \ICoh(\YY) \to \ICoh(\YY)
    \]
defines a lax equivariance datum 
\begin{equation}
\label{eqn:phi_pull}
    \ICoh(\YY)
    \to
    \RES_\phi(\ICoh(\YY))
\end{equation}
for the endofunctor
    \[
    \RES_\phi : \QCoh(\YY) \to \QCoh(\YY),
    \]
whose class is the aforementioned enhanced trace.

\subsubsection{} The monoidal functor
    \[
    \QCoh(\YY)
    \overset{\Upsilon}{\too}
    \ICoh(\YY)
    \overset{\Delta_*}{\too}
    \ICoh((\YY \times \YY)^\wedge_\YY).
    \]
factors through a monoidal functor
\begin{equation}
\label{eqn:QCoh_H}
    \QCoh(\YY)
    \to
    \HH(\YY).
\end{equation}
As we shall see, the relation between the objects (\ref{eqn:enh_tr_SSh}) and (\ref{eqn:enh_tr_SSh}) is mediated by the functor of restriction
    \[
    \OBLV :
    \HH(\YY) \MOD
    \to
    \QCoh(\YY) \MOD
    \]
along (\ref{eqn:QCoh_H}). 

\medskip

The aforementioned restriction functor enjoys the following compatibility with respect to a morphism $\ZZ \to \YY$ of quasi-smooth stacks:
\begin{equation}
\label{eqn:H_QCoh_oblv_pullback}
    \begin{tikzcd}
        \HH(\YY) \MOD \arrow[r, "\OBLV"] \arrow[d, "\HH_{\ZZ \to \YY}" left] & \QCoh(\YY) \MOD \arrow[d, "\IND"] \\
        \HH(\ZZ) \MOD \arrow[r, "\OBLV"] & \QCoh(\ZZ) \MOD
    \end{tikzcd}
\end{equation}
This is verified by a straightforward monad computation (see \cite[Proposition 3.4.1]{beraldo2021center} for a more general assertion).

\subsubsection{} Let us apply (\ref{eqn:H_QCoh_oblv_pullback}) to our endomorphism $\phi : \YY \to \YY$. We claim that the functor
    \[
    \OBLV : \HH(\YY) \MOD \to \QCoh(\YY) \MOD
    \]
is right adjointable, and hence the diagram 
\begin{equation}
\label{eqn:H_Q_phi}
    \begin{tikzcd}
        \HH(\YY) \MOD \arrow[r, "\HH_\phi"] \arrow[d, "\OBLV"] & \HH(\YY) \MOD \arrow[d, "\OBLV"] \\
        \QCoh(\YY) \MOD \arrow[r, "\IND_\phi"] & \QCoh(\YY) \MOD
    \end{tikzcd}
\end{equation}
induces a functor
\begin{equation}
\label{eqn:H_Q_phi_tr}
    \tr(\HH_\phi; \HH(\YY) \MOD)
    \to
    \tr(\IND_\phi; \QCoh(\YY) \MOD).
\end{equation}

\medskip

Indeed, the functor $\OBLV$ is given by the bimodule
    \[
    \HH(\YY)
    \in
    (\QCoh(\YY), \HH(\YY)),
    \]
and $\HH(\YY)$ is dualizable as a $\QCoh(\YY)$-module because it is dualizable as a plain dg category. Therefore, the existence of the claimed right adjoint is a consequence of the following general assertion: 

\begin{proposition}
\label{prop:A_dualizable}
Let $\AA$ and $\BB$ be monoidal categories, and let $\MM$ be an $(\AA, \BB)$-bimodule. Suppose that $\MM$ admits a right dual $\MM^\vee$ as an $\AA$-module. Then $\MM^\vee$ is naturally a right dual of $\MM$ as an $(\AA, \BB)$-bimodule.
\end{proposition}

\subsubsection{Proof} We need to show that the adjunction
    \[
    \MM \otimes_\BB - : 
    \BB \MOD \rightleftharpoons \AA \MOD
    : \Maps_\AA(\MM, -)
    \]
actually takes place in $\twoDGCAT$. That is, we have to see that the right adjoint is continuous and commutes with the operation of tensoring with a plain dg category. 

\medskip

It suffices to check this condition after composing with the forgetful functor 
    \[
    \DGCAT \leftarrow \BB \MOD
    : \OBLV
    \]
Our assumption on $\MM$ says that this composite participates in an adjunction in $\twoDGCAT$
    \[
    \MM :
    \DGCAT 
    \rightleftharpoons
    \AA \MOD
    \rightleftharpoons
    \BB \MOD
    : \MM^\vee
    \]
so the result follows.

\qed[\cref{prop:A_dualizable}]

\subsubsection{} Now observe that $\QCoh(\YY)\MOD$ is canonically self dual because $\QCoh(\YY)$ is symmetric monoidal, and that this duality identifies the endofunctor $(\IND_\phi)^\vee$ with $\RES_\phi$. Therefore, by (\ref{eqn:tr_of_dual}), we obtain an equivalence
\begin{equation}
\label{eqn:Q_ind_res_tr}
    \tr(\IND_\phi; \QCoh(\YY) \MOD)
    \simeq
    \tr(\RES_\phi; \QCoh(\YY) \MOD).
\end{equation}
We claim:

\begin{computation}
\label{computation:tr_oblv}
The composite
\begin{equation}
\label{eqn:tr_oblv}
\begin{split}
    \tr(\HH_\phi; \HH(\YY) \MOD)
    & \overset{(\ref{eqn:H_Q_phi_tr})}{\longrightarrow}
    \tr(\IND_\phi; \QCoh(\YY) \MOD) \\
    & \qquad \qquad \simeq
    \tr(\RES_\phi; \QCoh(\YY) \MOD)
\end{split}
\end{equation}
maps $\cl(\ICoh(\YY), \alpha)$ to $\cl(\ICoh(\YY), \phi^!)$.
\end{computation}

\subsubsection{} Note that after composing with (\ref{eqn:H_Q_phi}), $\alpha$ defines a lax equivariance datum in $\QCoh(\YY) \MOD$, which we shall denote by the same symbol:
    \[
    \alpha : \ICoh(\YY) \to \IND_\phi(\ICoh(\YY)).
    \]
The functor (\ref{eqn:H_Q_phi_tr}) takes
    \[
    \cl(\ICoh(\YY), \alpha) \in \tr(\HH_\phi; \HH(\YY) \MOD)
    \]
to
    \[
    \cl(\ICoh(\YY), \alpha) \in \tr(\IND_\phi; \QCoh(\YY) \MOD).
    \]

\medskip

Let us apply Computation \ref{computation:adj_dual_tr} to $\alpha$, viewed as a natural transformation
    \[
    \begin{tikzcd}
        \DGCAT \arrow[d, "\ICoh(\YY)" left] \arrow[r, "\id"] & \DGCAT \arrow[d, "\ICoh(\YY)"] \\
        \QCoh(\YY) \arrow[r, "\IND_\phi" below] & \QCoh(\YY)
        \arrow[Rightarrow, from=1-2, to=2-1]
    \end{tikzcd}
    \]
The self-duality of $\ICoh(\YY)$ as a module over $\QCoh(\YY)$ identifies the functor
    \[
    \ICoh(\YY) : \DGCAT \to \QCoh(\YY)
    \]
with the left adjoint of its dual. We obtain that the isomorphism (\ref{eqn:Q_ind_res_tr}) takes $\cl(\ICoh(\YY), \alpha)$ to the class of the lax equivariance datum
    \[
    \beta : \ICoh(\YY) \to \RES_\phi(\ICoh(\YY))
    \]
obtained by adjunction from the map
\begin{equation}
\label{eqn:alpha_dual}
\begin{split}
    \IND_\phi(\ICoh(\YY))
    & \simeq
    \IND_\phi(\ICoh(\YY)^\vee) \\
    & \qquad \simeq
    \IND_\phi(\ICoh(\YY))^\vee
    \overset{\alpha^\vee}{\longrightarrow}
    \ICoh(\YY)^\vee
    \simeq
    \ICoh(\YY).
\end{split}    
\end{equation}
It remains to see that $\beta$ is canonically identified with (\ref{eqn:phi_pull}). 

\medskip

It is equivalent to identify (\ref{eqn:alpha_dual}) with the transpose
\begin{equation}
\label{eqn:phi_pull_transpose}
    \IND_\phi(\ICoh(\YY)) 
    \to
    \ICoh(\YY)
\end{equation}
of (\ref{eqn:phi_pull}). This is a formal consequence of the commutative diagram of $(\QCoh(\YY), \HH(\YY))$-bimodules:
    \[
    \begin{tikzcd}
        \IND_\phi(\ICoh(\YY)) \otimes \ICoh(\YY) \arrow[r, "\text{(\ref{eqn:alpha_definition})}"] \arrow[d, "\text{(\ref{eqn:phi_pull_transpose})} \, \otimes \id"] & \IND_\phi(\HH(\YY)) \arrow[d, "\text{(\ref{eqn:H_Q_phi})}"] \\
        \ICoh(\YY) \otimes \ICoh(\YY) \arrow[r, "(\IND_\phi \otimes \id) \circ \, \text{(\ref{eqn:alpha_definition})}"] & \HH_\phi
    \end{tikzcd}
    \]

\qed[\cref{computation:tr_oblv}]

\subsubsection{Remark} We do not currently have a good description of the composite appearing in Computation \ref{computation:tr_oblv} as a map
    \[
    \DMod(\Sing_2(\YY^\phi))^\Rightarrow
    \to
    \QCoh(\YY^\phi).
    \]
However, we expect that such a description would contain interesting geometric content.

\subsection{Variant with singular support}\label{s:6-sing-supp}

\subsubsection{} Consider a closed conical subset $\Nilp$ of $\Sing(\YY)$. When we regard $\Nilp$ as living over $\Sing(\YY) \times \Sing(\YY)$, we will do so via the map
    \[
    \Nilp 
    \mon
    \Sing(\YY) 
    \simeq 
    \Sing(\YY/\YY \times \YY)
    \to
    \Sing(\YY \times \YY).
    \]

\subsubsection{} The category $\ShCat_\SSh(\YY)$, being a limit of symmetric monoidal categories under symmetric monoidal functors, carries a natural symmetric monoidal structure. We let
    \[
    \pr_\Nilp : \ShCat_\SSh(\YY) \to \ShCat_\SSh(\YY)
    \]
denote the idempotent endofunctor of $\ShCat_\SSh(\YY)$ given by tensoring with $\ul\DMod(\Nilp)^\Rightarrow$. One has:
    \[
    \GAMMA(S; \pr_\Nilp(\CC))
    \simeq
    \DMod(\Nilp_S)^\Rightarrow 
    \underset{\DMod(\Sing(S))^\Rightarrow}{\otimes}
    \GAMMA(S; \CC).
    \]

\medskip

Since $\ul\DMod(\Nilp)^\Rightarrow$ is a colocalization of $\ul\DMod(\Sing(\YY))^\Rightarrow$, we obtain that $\pr_\Nilp$ is a colocalization of the identity endofunctor of $\ShCat_\SSh(\YY)$
\begin{equation}
\label{eqn:pr_Nilp_colocalization}
    \pr_\Nilp
    \rightleftharpoons
    \id.
\end{equation}
Since $\pr_\Nilp$ is idempotent, this implies $\pr_\Nilp$ is canonically adjoint to itself, with unit and counit
    \[
    \eta :
    \pr_\Nilp \circ \pr_\Nilp \simeq \pr_\Nilp 
    \rightleftharpoons
    \id 
    : \epsilon.
    \]
In particular, $\pr_\Nilp$ admits a right adjoint.
\begin{proposition}
\label{proposition:morita_singular_support}
The Morita equivalence (\ref{eqn:morita_ShCat}) exchanges the colocalizations
    \[
    \ul\DMod(\Nilp)^\Rightarrow
    \qquad
    \text{and}
    \qquad
    \LOC(\ICoh_\Nilp(\YY))
    \]
\end{proposition}

\subsubsection{Proof} Consider the image of the colocalization
    \[
    \ul\DMod(\Nilp)^\Rightarrow
    \rightleftharpoons
    \ul\DMod(\Sing(\YY))^\Rightarrow
    \]
under the Morita equivalence (\ref{eqn:morita_ShCat}). According to \cref{prop:ICoh_N_compat}, its value at an object $S \in \Sm_{/\YY}$ is the colocalization
\begin{equation}
\label{eqn:ICoh_Nilp_S_colocalization}
    \ICoh_{\Nilp_S}(S)
    \rightleftharpoons
    \ICoh(S),
\end{equation}
so its image under $\GAMMA$ is the colocalization of $\ICoh(\YY)$ obtained from (\ref{eqn:ICoh_Nilp_S_colocalization}) by passing to the limit over $\Sm_{/\YY}$. This is definitionally the colocalization
    \[
    \ICoh_\Nilp(\YY)
    \rightleftharpoons
    \ICoh(\YY),
    \]
so the result now follows from H-affineness.

\qed[\cref{proposition:morita_singular_support}]

\subsubsection{} After applying the Morita equivalence (\ref{eqn:morita_ShCat}) and H-affineness, we obtain from $\pr_\Nilp$ an endomorphism
    \[
    \pr_\Nilp : 
    \HH(\YY) \MOD
    \to
    \HH(\YY) \MOD,
    \]
to be denoted by the same symbol.

\medskip

Let
    \[
    \HH_\Nilp
    \in
    (\HH(\YY), \HH(\YY)) \BIMOD
    \]
denote the bimodule responsible for the endofunctor $\pr_\Nilp$. In these terms, the colocalization (\ref{eqn:pr_Nilp_colocalization}) manifests as a colocalization
    \[
    \HH_\Nilp
    \rightleftharpoons
    \HH(\YY)
    \]
of $(\HH(\YY), \HH(\YY))$-bimodules. The idempotence of $\pr_\Nilp$ says that
\begin{equation}
\label{eqn:proj_Nilp_bimod_idempotence}
    \HH_\Nilp \underset{\HH(\YY)}{\otimes} \HH_\Nilp \simeq \HH_\Nilp
\end{equation}
as full subcategories of $\HH(\YY)$ stable under the left and right actions of $\HH(\YY)$.

\subsubsection{} We obtain a map
    \[
    \pr_\Nilp \circ \HH_\phi
    \to
    \pr_\Nilp \circ \HH_\phi \circ \pr_\Nilp
    \to
    \HH_\phi \circ \pr_\Nilp
    \]
encoding a diagram:
\begin{equation}
\label{eqn:proj_Nilp_lax}
    \begin{tikzcd}
        \HH(\YY) \MOD \arrow[r, "\HH_\phi"] \arrow[d, "\pr_\Nilp"] & \HH(\YY) \MOD \arrow[d, "\pr_\Nilp"] \\
        \HH(\YY) \MOD \arrow[r, "\HH_\phi"] & \HH(\YY) \MOD
        \arrow[Rightarrow, from=1-2, to=2-1]
    \end{tikzcd}
\end{equation}
Therefore, $\pr_\Nilp$ induces an endomorphism of the trace:
\begin{equation}
\label{eqn:Nilp_proj_trace}
    \tr(\HH_\phi; \HH(\YY)\MOD)
    \to
    \tr(\HH_\phi; \HH(\YY)\MOD).
\end{equation}
We claim:

\begin{computation}
\label{computation:proj_Nilp_trace}
Under the identification (\ref{eqn:trace_computation}), this endomorphism identifies with the operation
    \[
    (\omega_{\Nilp^\phi})^\Rightarrow \otimes (-) :
    \DMod(\Sing_2(\YY^\phi))^\Rightarrow
    \to
    \DMod(\Sing_2(\YY^\phi))^\Rightarrow,
    \]
where
    \[
    \Nilp^\phi = \Sing_\phi(\YY/\YY \times \YY) \underset{\Sing(\YY \times \YY)}{\times} \Nilp.
    \]
\end{computation}

\subsubsection{} By the idempotence (\ref{eqn:proj_Nilp_bimod_idempotence}) of $\HH_\Nilp$ as a bimodule, (\ref{eqn:Nilp_proj_trace}) is obtained from the endomorphism
    \[
    \HH(\YY)
    \to
    \HH_\Nilp
    \to
    \HH(\YY)
    \]
by tensoring up along
    \[
    - \underset{\HH(\YY) \otimes \HH(\YY)^\rev}{\otimes} \HH_\phi.
    \]

\medskip

Now observe that
    \[
    \HH_\Nilp
    \simeq
    \pr_{\Nilp \times \Nilp}(\HH(\YY))
    \]
as objects of
    \[
    (\HH(\YY), \HH(\YY)) \BIMOD
    \simeq
    \HH(\YY \times \YY) \MOD.
    \]
Therefore, it follows from \cref{prop:rel_D_stack_0} that (\ref{eqn:Nilp_proj_trace}) corresponds the endomorphism of (\ref{eqn:tr_H_SSh}) induced by the composite
    \[
    \ul\DMod(\Sing(\YY/\YY \times \YY))^\Rightarrow
    \to
    \ul\DMod(\Nilp)^\Rightarrow
    \to
    \ul\DMod(\Sing(\YY/\YY \times \YY)
    \]

\qed[\cref{computation:proj_Nilp_trace}]

\subsubsection{} Consider the pullback of $\Nilp$ along the projection $\Sing_\phi(\YY) \to \Sing(\YY)$, to be denoted $\Nilp_\phi$. We impose the assumption:
\begin{equation}
\label{eqn:Nilp_phi_assumption}
\tag{$*$}
    \text{The image of $\Nilp \underset{\Sing(\YY)}{\times} \Sing_\phi(\YY)$ under $\Sing(\phi)$ is contained in $\Nilp$.}
\end{equation}
In other words, we require the existence of the dashed arrow in the diagram:
    \[
    \begin{tikzcd}
        \Nilp \underset{\Sing(\YY)}{\times} \Sing_\phi(\YY) \arrow[rr, dashed] \arrow[dd] \arrow[dr] & & \Nilp \arrow[dr] & \\
        & \Sing_\phi(\YY) \arrow[rr] \arrow[dd] & & \Sing(\YY) \\
        \Nilp \arrow[dr] & & & \\ 
        & \Sing(\YY) & &
    \end{tikzcd}
    \] 
The reason for imposing the condition (\ref{eqn:Nilp_phi_assumption}) is so that the following assertion holds true:

\begin{proposition}
\label{prop:Nilp_phi_containment}
Under the assumption (\ref{eqn:Nilp_phi_assumption}), the map
    \[
    \pr_\Nilp \circ \HH_\phi \circ \pr_\Nilp
    \to
    \HH_\phi \circ \pr_\Nilp.
    \]
is an isomorphism.
\end{proposition}

\subsubsection{Proof} Note that \cref{prop:Nilp_phi_containment} asserts that the essential image of $\HH_\phi \circ \pr_\Nilp$ is contained in that of $\pr_\Nilp$. In other words, we need to show for every object $S \in \Sm_{/\YY}$ that the essential image of the composite
\begin{equation}
\label{eqn:proj_Nilp_ev}
\begin{split}
    & \HH(\YY) \MOD 
    \overset{\pr_\Nilp}{\too}
    \HH(\YY) \MOD
    \overset{\HH_\phi}{\too}
    \HH(\YY) \MOD \\
    & \qquad \overset{\GAMMA(\YY_S; -)}{\too}
    \HH(\YY_S) \MOD
    \overset{\ICoh(\YY_S)^\tau}{\too}
    \DMod(\Sing(\YY_S))^\Rightarrow \MOD
\end{split}
\end{equation}
is contained in
    \[
    \DMod(\Nilp_{\YY_S})^\Rightarrow \MOD
    \subset
    \DMod(\Sing(\YY_S))^\Rightarrow \MOD.
    \]

\medskip

By the functoriality of the $\HH$-theory, (\ref{eqn:proj_Nilp_ev}) identifies with
\begin{equation}
\label{eqn:proj_Nilp_ev_2}
\begin{split}
    & \HH(\YY) \MOD 
    \overset{\pr_\Nilp}{\too}
    \HH(\YY) \MOD
    \overset{\GAMMA(S; -)}{\too}
    \HH(S) \MOD \\
    & \qquad \overset{\HH_{\YY_S \to S}}{\too}
    \HH(\YY_S) \MOD
    \overset{\ICoh(\YY_S)}{\too}
    \DMod(\Sing(\YY_S))^\Rightarrow \MOD.
\end{split}
\end{equation}
Applying (\ref{eqn:morita_pull}) rewrites this composite as
\begin{equation}
\begin{split}
    & \HH(\YY) \MOD 
    \overset{\pr_\Nilp}{\too}
    \HH(\YY) \MOD
    \overset{\GAMMA(S; -)}{\too}
    \HH(S) \MOD \\
    & \qquad \overset{\ICoh(S)}{\too}
    \DMod(\Sing(S))^\Rightarrow \MOD
    \overset{\DMod(\Sing(S)_{\YY_S})^\Rightarrow}{\too}
    \DMod(\Sing(\YY_S))^\Rightarrow \MOD,
\end{split}
\end{equation}
and applying the definition of $\pr_\Nilp$ further rewrites our functor as 
\begin{equation}
\begin{split}
    & \HH(\YY) \MOD 
    \overset{\GAMMA(S; -)}{\too}
    \HH(S) \MOD
    \overset{\ICoh(S)}{\too}
    \DMod(\Sing(S))^\Rightarrow \MOD \\
    & \qquad \overset{\DMod(\Nilp_S)^\Rightarrow}{\too}
    \DMod(\Sing(S))^\Rightarrow \MOD
    \overset{\DMod(\Sing(S)_{\YY_S})^\Rightarrow}{\too}
    \DMod(\Sing(\YY_S))^\Rightarrow \MOD.
\end{split}
\end{equation}
The desired containment now follows from (\ref{eqn:Nilp_phi_assumption}), base changed to $S$.

\qed[\cref{prop:Nilp_phi_containment}]

\subsubsection{} By Propositions \ref{proposition:morita_singular_support} and \ref{prop:Nilp_phi_containment}, $\alpha$ descends to a lax equivariance datum $\alpha_\Nilp$ fitting into the diagram:
\begin{equation}
\label{eqn:alpha_Nilp}
\begin{tikzcd}
    \ICoh(\YY) \arrow[d] \arrow[r, "\alpha"] & \HH_\phi \underset{\HH(\YY)}{\otimes} \ICoh(\YY) \arrow[d] \\
    \ICoh_\Nilp(\YY) \arrow[r, dashed, "\alpha_\Nilp"] & \HH_\phi \underset{\HH(\YY)}{\otimes} \ICoh_\Nilp(\YY)
\end{tikzcd}
\end{equation}

\medskip

Observe that $\alpha_\Nilp$ encodes the concatenation of the lax equivariance datum
\[
\begin{tikzcd}
    \DGCAT \arrow[r] \arrow[d, "\ICoh(\YY)" left] & \DGCAT \arrow[d, "\ICoh(\YY)"] \\
    \HH(\YY) \arrow[r, "\HH_\phi" below] & \HH(\YY)
    \arrow[Rightarrow, from=1-2, to=2-1]
\end{tikzcd}
\]
of $\alpha$ with (\ref{eqn:proj_Nilp_lax}). Therefore, the combination of Computations \ref{computation:bimod_trace} and \ref{computation:proj_Nilp_trace} shows that the equivalence (\ref{eqn:trace_computation}) takes
\begin{equation}
\label{eqn:cl_alpha_Nilp}
    \cl(\ICoh_\Nilp(\YY), \alpha_\Nilp)
    \in
    \tr(\HH_\phi; \HH(\YY) \MOD)
\end{equation}
to the object
    \[
    (\omega_{\Nilp^\phi})^\Rightarrow 
    \in \DMod(\Sing_2(\YY^\phi))^\Rightarrow.
    \]

\medskip

We obtain from (\ref{eqn:alpha_Nilp}) that (\ref{eqn:H_QCoh_oblv_pullback}) restricts to a lax equivariance datum $(\phi^!)_\Nilp$ fitting into the diagram:
\[
\begin{tikzcd}
    \ICoh_\Nilp(\YY) \arrow[r, "(\phi^!)_\Nilp"] \arrow[d, hook] & \RES_\phi(\ICoh_\Nilp(\YY)) \arrow[d, hook] \\
    \ICoh(\YY) \arrow[r, "\phi^!"] & \RES_\phi(\ICoh(\YY))
\end{tikzcd}
\]
As a result, \cref{computation:tr_oblv} shows that (\ref{eqn:tr_oblv}) maps (\ref{eqn:cl_alpha_Nilp}) to
    \[
    \cl(\ICoh_\Nilp(\YY), (\phi^!)_\Nilp)
    \in
    \tr(\RES_\phi; \QCoh(\YY) \MOD).
    \]

\subsubsection{Remark} This computation implies \cite[Conjecture 24.6.9]{arinkin2020stack}.

\appendix

\section{Locally fully faithful functors}

In this appendix we will formulate and prove a lifting criterion \cref{aproposition:lff_lifting} for locally fully faithful functors.

\subsubsection{} Recall that a functor of $(\infty, 2)$-categories
    \[
    F : \CC \to \DD
    \]
is said to be \emph{locally fully faithful} if, for every pair of objects $c, c' \in \CC$, the functor
    \[
    \Maps_\CC(c, c')
    \to
    \Maps_\DD(F(c), F(c'))
    \]
is fully faithful.

\medskip

Let us say right away that \cref{aproposition:lff_lifting} is concerned with functors out of an $(\infty, 1)$-category. Since the condition of being locally fully faithful is preserved by passage to $(\infty, 1)$-cores, the ambit of this appendix is contained entirely in the world of $(\infty, 1)$-categories.

\subsubsection{} Recall that for any category $\CC$, one associates its twisted arrows category $\Tw(\CC)$, whose objects are arrows $c \to c'$ in $\CC$ and whose morphisms from $c \to c'$ to $d \to d'$ are\footnote{\emph{Warning}: the literature does not appear to be consistent as to whether the `twisted arrows category' is what we have denoted $\Tw(\CC)$ or rather its opposite $\Tw(\CC)^\op$.} commutative diagrams:
\[
\begin{tikzcd}
    c \arrow[d] & d \arrow[l] \arrow[d] \\
    c' \arrow[r] & d'
\end{tikzcd}
\]
The formation of the twisted arrows category is natural in $\CC$; every functor $F : \CC \to \DD$ induces a functor $\Tw(F) : \Tw(\CC) \to \Tw(\DD)$.

\medskip

One has the twisted arrows formula for the space of natural transformations between a pair of functors $\CC \to \DD$:
\begin{equation}
\label{eqn:tw_arr_formula}
    \Maps(F, G) \simeq \lim_{\Tw(\CC)} \Maps_\DD(F(c), G(c')).
\end{equation}

\subsection{Lax surjective functors}

\subsubsection{} We will say that a functor $F: \CC \to \DD$ is \emph{lax surjective} if for every object $d \in \DD$ there is an object $c \in \CC$ and a map $F(c) \to d$. 

\begin{aproposition}
\label{prop:tw-ar-lax-surj}
If $F: \CC \to \DD$ is essentially surjective, then $\Tw(F) : \mathrm{Tw}(\CC) \to \mathrm{Tw}(\DD)$ is lax surjective. 
\end{aproposition}

\subsubsection{Proof} Let $f : d \to d'$ be an object of $\mathrm{Tw}(\DD)$. Since $F$ is essentially surjective, we may choose $c \in \CC$ along with an isomorphism $\alpha: F(c) \simto d$. Then the diagram
\[
\begin{tikzcd}
    F(c) \arrow[d, "\id"] & d \arrow[l, "\alpha^{-1}"'] \arrow[d, "f"]\\
    F(c) \arrow[r, "f \circ \alpha"] & d'
\end{tikzcd}
\]
defines a map from $F(c \to c)$ to $d \to d'$ in $\mathrm{Tw}(\DD)$. 

\qed[\cref{prop:tw-ar-lax-surj}]

\subsubsection{} Let $\Prop \subset \Cat$ denote the full subcategory of of propositions---that is, $\Prop$ is the full subcategory of $\Cat$ spanned by $* = \top$ and $\varnothing = \bot$. This subcategory is closed under limits and the limit of a diagram of propositions
    \[
    F : I \to \Prop
    \]
is the proposition
    \[
    (\forall i \in I)(F(i)).
    \]

\begin{aproposition}\label{prop:lax-surj-prop}
Let $F: J \to \Prop$ be a diagram of propositions. For any lax surjective functor $\alpha : I \to J$, the natural morphism
    \[
    \lim_{j \in J} F(j) \to \lim_{i \in I} F(\alpha(i))
    \]
is an isomorphism.
\end{aproposition}

\subsubsection{Proof} We need to establish the reverse implication:
    \[
    (\forall i \in I)(F(\alpha(i)))
    \to
    (\forall j \in J)(F(j)).
    \]
Since $\alpha$ is lax surjective, we have
    \[
    (\forall j \in J)(\exists i \in I)(F(\alpha(i)) \to F(j)),
    \]
so the conclusion follows.
\qed[Proposition \ref{prop:lax-surj-prop}]

\subsubsection{} Let $\Poset \subset \Cat$ denote the full subcategory of partially ordered sets. Observe that this full subcategory is closed under limits, and that a category $\CC$ is a poset if and only if the functor $\CC \to *$ is locally fully faithful.

\begin{aproposition}
\label{prop:lax-surj-lifting-criterion}
Let $F : J \to \Poset$ be a diagram of posets. For any lax surjective functor $\alpha : I \to J$, the natural morphism
\begin{equation}
\label{eqn:poset_lim}
    \lim_{j \in J} F(j) \to \lim_{i \in I} F(\alpha(i))
\end{equation}
is fully faithful. Furthermore, the essential image of this morphism consists of objects $(x_i : i \in I)$ satisfying the property:

\begin{equation}
\tag{$*$}
\begin{minipage}{28em}
    For every tuple $(j, i, i', \phi, \phi')$ determining a diagram of the form
    \[
    \begin{tikzcd}[column sep=1em, row sep=1em]
        \alpha(i) \arrow[dr, "\phi" below left] & & \alpha(i') \arrow[dl, "\phi'"] \\
        & j &
    \end{tikzcd}
    \]
    we have
        \[
        F_\phi(x_i) = F_{\phi'}(x_{i'}).
        \]
\end{minipage}    
\end{equation}
\end{aproposition}

\subsubsection{Proof} The mapping space between two objects of a poset is always a proposition, so the assertion about fully faithfulness follows from Propositions \ref{prop:lax-surj-prop}.

\medskip

It is clear that every object in the image of the morphism (\ref{eqn:poset_lim}) satisfies ($*$). It remains to check that every object $(x_i)$ of $\lim F(\alpha(i))$ satisfying ($*$) admits a lift $(x_j)$ to $\lim F(j)$. For $j \in J$, set 
    \[
    x_j = F_\phi(x_i)
    \]
where $\phi$ is any morphism $\alpha(i) \to j$. Such a morphism exists because $\alpha$ is lax surjective, and ($*$) says precisely that the value of $F_\phi(x_i)$ is independent of the choice of $\phi$.

\qed[\cref{prop:lax-surj-lifting-criterion}]

\subsection{Lifting along locally fully faithful functors}

\subsubsection{} For any category $\CC$ we let 
    \[
    \CC^\obj = \pi_0(\CC^{\,\simeq})
    \]
denote the category with the same objects as $\CC$ but where the only maps are the identity endomorphisms. So $\CC^\obj$ is a set and we may identify $\Tw(\CC^\mathrm{\obj}) \simeq \CC^\mathrm{\obj}$

\medskip

We choose a basepoint for every connected component of $\CC^{\,\simeq}$. This defines a (noncanonical) functor:
    \[
    \eta : \CC^\obj \to \CC^{\,\simeq} \to \CC.
    \]
In particular, we obtain a map
\begin{equation}
\label{eqn:poset_lim_tw_arr}
    \lim_{\Tw(\CC)} F(c \to c')
    \to
    \lim_{\Tw(\CC^\obj)} F(\eta(d \overset{\id}{\to} d))
    \simeq
    \underset{d \in \CC^\obj}{\sqcap} F(d),
\end{equation}
where by a slight abuse of notation we have denoted
    \[
    F(d) = F(\eta(d \overset{\id}{\to} d)).
    \]

\begin{aproposition}
\label{prop:limit-twisted-arrows-poset}
Let $F: \Tw(\CC) \to \Poset$ be a diagram of posets. The map (\ref{eqn:poset_lim_tw_arr}) is fully faithful. Furthermore, its essential image consists of those objects $(x_d)$ with the property: 
\begin{equation}
\label{eqn:tw_poset_property}
\tag{$**$}
    \text{For every arrow $c \to c'$ in $\CC$, the images of $x_c$ and $x_{c'}$ in $F(c \to c')$ are equal.}
\end{equation}
\end{aproposition}

\subsubsection{Proof} The functor $\Tw(\CC^\obj) \to \Tw(\CC)$ is lax surjective by \cref{prop:tw-ar-lax-surj}, so we may apply \cref{prop:lax-surj-lifting-criterion}. We need to show for any pair of diagrams
\begin{equation}
\label{eqn:tw_arr_diagrams}
    \begin{tikzcd}
        d \arrow[d, "\id"] & c \arrow[d, "f"] \arrow[l] \\
        d \arrow[r] & c'
    \end{tikzcd}
        \qquad \qquad
    \begin{tikzcd}
        d' \arrow[d, "\id"] & c \arrow[d, "f"] \arrow[l] \\
        d' \arrow[r] & c'
    \end{tikzcd}
\end{equation}
that the images of $x_d$ and $x_{d'}$ in $F(c \to c')$ are equal. 

\medskip

We are free to assume that $d' = c'$ and that the lower horizontal arrow in the right square of (\ref{eqn:tw_arr_diagrams}) is the identity. Then the two squares appearing in (\ref{eqn:tw_arr_diagrams}) may be rewritten as the composites:
    \[
    \begin{tikzcd}
        d \arrow[d, "\id"] & d \arrow[l] \arrow[d] & c \arrow[d, "f"] \arrow[l] \\
        d \arrow[r] & c' \arrow[r] & c'
    \end{tikzcd}
    \qquad \qquad
    \begin{tikzcd}
        c' \arrow[d, "\id"] & d \arrow[l] \arrow[d] & c \arrow[d, "f"] \arrow[l] \\
        c' \arrow[r] & c' \arrow[r] & c'
    \end{tikzcd}
    \]
and the condition (\ref{eqn:tw_poset_property}) implies that the images of $x_d$ and $x_{c'}$ in $F(d \to c')$ are equal.

\qed[\cref{prop:limit-twisted-arrows-poset}]

\begin{aproposition} 
\label{aproposition:lff_lifting}
Let $\mathscr{S}$ be a category, and let
    \[
    \begin{tikzcd}
        & \ZZ \arrow[d, "G"] \\
        \XX \arrow[r, "F"] & \YY
    \end{tikzcd}
    \]
be a diagram in $\Fun(\mathscr{S}, \Cat)$. Suppose that $G$ is pointwise fully faithful, in the sense that for every object $s$ of $\mathscr{S}$, the functor
    \[
    G(s) : \ZZ(s) \to \YY(s)
    \]
is locally fully faithful. Then the functor
\begin{equation}
\label{eqn:lift_ptwise}
\begin{split}
    & \Maps(\XX, \ZZ) \underset{\Maps(\XX, \YY)}{\times} \{F\} \\
    & \qquad \qquad \to
    \underset{s \in \mathscr{S}^\obj}{\sqcap} 
    \left(
    \Fun(\XX(s), \ZZ(s)) \underset{\Fun(\XX(s), \YY(s))}{\times} \{F(s)\}
    \right)
\end{split}
\end{equation}
is fully faithful, and a collection $(\tilde{F}_s)$ belongs to its essential image if and only if: for every arrow $s \to s'$ in $\mathscr{S}$, the images of
    \[
    \tilde{F}_s 
    \in 
    \Fun(\XX(s), \ZZ(s))
    \underset{\Fun(\XX(s), \YY(s))}{\times}
    \{F(s)\}
    \]
and
    \[
    \tilde{F}_{s'} 
    \in 
    \Fun(\XX(s'), \ZZ(s'))
    \underset{\Fun(\XX(s'), \YY(s'))}{\times}
    \{F(s')\}
    \]
in
    \[
    \Fun(\XX(s), \ZZ(s')) \underset{\Fun(\XX(s), \YY(s'))}{\times} \{F\}
    \]
are equal.
\end{aproposition}

\subsubsection{Remark} This result says that a lift of $F$ to a functor $\tilde{F} : \XX \to \ZZ$ amounts to a collection of lifts
    \[
    \begin{tikzcd}
        & \ZZ(s) \arrow[d, "G"] \\
        \XX(s) \arrow[r, "F"] \arrow[ur, "\tilde{F}_s", dashed] & \YY(s)
    \end{tikzcd}
    \]
for $s \in \mathscr{S}^\obj$ with the property that for every arrow $s \to s'$ in $\mathscr{S}$, the diagram
    \[
    \begin{tikzcd}[row sep=1em, column sep=1.5em]
        \XX(s) \arrow[dd] \arrow[dr, dashed] \arrow[drrr] & & & \\
        & \ZZ(s) \arrow[rr] \arrow[dd] & & \YY(s) \arrow[dd] \\
        \XX(s') \arrow[dr, dashed] \arrow[drrr] & & & \\
        & \ZZ(s') \arrow[rr] & & \YY(s')
    \end{tikzcd}
    \]
commutes. The assumption that $G$ is locally fully faithful gaurantees that the commutativity of this diagram is a property and not additional data.

\subsubsection{Proof} By the twisted arrows formula, we may perceive (\ref{eqn:lift_ptwise}) as the restriction morphism
\begin{equation}
\label{eqn:tw_arr_lift}
\begin{split}
    & \lim_{\Tw(\mathscr{S})} 
    \left(
    \Fun(\XX(s), \ZZ(s')) 
    \underset{\Fun(\XX(s), \YY(s'))}{\times} 
    \{F\}
    \right) \\
    & \qquad \qquad \to
    \underset{s \in \mathscr{S}^\obj}{\sqcap} 
    \left( 
    \Fun(\XX(s), \ZZ(s)) 
    \underset{\Fun(\XX(s), \YY(s))}{\times} 
    \{F(s)\}
    \right).
\end{split}
\end{equation}
The functors
    \[
    G(s') \circ (-) :
    \Fun(\XX(s), \ZZ(s')) 
    \to
    \Fun(\XX(s), \YY(s'))
    \]
are locally fully faithful because $G(s')$ is locally fully faithful, so the LHS of (\ref{eqn:tw_arr_lift}) is a limit of posets and the result follows from \cref{prop:limit-twisted-arrows-poset}.

\qed[Proposition \cref{aproposition:lff_lifting}]

\bibliographystyle{alpha}
\bibliography{ref}

\end{document}